\pdfoutput=1
\RequirePackage{ifpdf}
\ifpdf % We~are running pdfTeX in pdf mode
\documentclass[pdftex]{sigma}
\else
\documentclass{sigma}
\fi

\numberwithin{equation}{section}

\newtheorem{Theorem}{Theorem}[section]
\newtheorem*{Theorem*}{Theorem}
\newtheorem{Corollary}[Theorem]{Corollary}
\newtheorem{Lemma}[Theorem]{Lemma}
\newtheorem{Proposition}[Theorem]{Proposition}

\theoremstyle{definition}
\newtheorem{Definition}[Theorem]{Definition}

\newtheorem{Example}[Theorem]{Example}
\newtheorem{Remark}[Theorem]{Remark}

\usepackage{ytableau}
\usepackage{mathrsfs}
\usepackage{bbm}
\DeclareMathOperator{\Par}{\mathbb{Y}}
\DeclareMathOperator{\RYT}{\mathrm{RYT}_{\geq 0}}
\DeclareMathOperator{\SYT}{\mathrm{SYT}}
\DeclareMathOperator{\PSYT}{\mathrm{PSYT}_{\geq 0}}
\DeclareMathOperator{\RSSYT}{\mathrm{RSSYT}_{\geq 0}}
\DeclareMathOperator{\sH}{\mathscr{H}}
\DeclareMathOperator{\sA}{\mathscr{A}}
\DeclareMathOperator{\sD}{\mathscr{D}}

\DeclareMathOperator{\Ind}{Ind}
\DeclareMathOperator{\Res}{Res}

\DeclareMathOperator{\sort}{sort}
\DeclareMathOperator{\revsort}{revsort}
\DeclareMathOperator{\Stab}{Stab}

\DeclareMathOperator{\MacD}{\mathfrak{P}}
\DeclareMathOperator{\sE}{\mathcal{E}}
\DeclareMathOperator{\Inv}{\mathrm{Inv}}
\DeclareMathOperator{\inv}{inv}
\DeclareMathOperator{\rk}{\mathrm{rk}}

\DeclareMathOperator{\Top}{\mathrm{top}}
\DeclareMathOperator{\Min}{\mathrm{min}}
\DeclareMathOperator{\I}{\mathrm{I}}
\DeclareMathOperator{\Hilb}{\mathrm{Hilb}}

\DeclareMathOperator{\End}{\mathrm{End}}

\begin{document}

\allowdisplaybreaks

\newcommand{\arXivNumber}{2405.00756}

\renewcommand{\PaperNumber}{076}

\FirstPageHeading

\ShortArticleName{Murnaghan-Type Representations for the Positive Elliptic Hall Algebra}

\ArticleName{Murnaghan-Type Representations\\ for the Positive Elliptic Hall Algebra}

\Author{Milo BECHTLOFF WEISING}

\AuthorNameForHeading{M.~Bechtloff Weising}

\Address{Department of Mathematics, Virginia Tech, 460~McBryde Hall, 225~Stanger Street,\\
Blacksburg, VA 24061-1026, USA}
\Email{\mail{milojbw@vt.edu}, \mail{milojbw@gmail.com}}
\URLaddress{\url{https://sites.google.com/view/milobechtloffweising/}}

\ArticleDates{Received July 01, 2025, in final form July 27, 2026; Published online August 14, 2026}

\Abstract{We construct a~new family of graded representations \smash{$\widetilde{W}_{\lambda}$} for the positive elliptic Hall algebra $\mathcal{E}^{+}$ indexed by Young diagrams $\lambda$ which generalize the standard $\mathcal{E}^{+}$ action on symmetric functions. These representations have homogeneous bases of eigenvectors for the action of the Macdonald element $P_{0,1} \in \mathcal{E}^{+}$ with distinct $\mathbb{Q}(q,t)$-rational spectrum generalizing the symmetric Macdonald functions. The analysis of the structure of these representations exhibits interesting combinatorics arising from the limits of periodic standard Young tableaux. We find an explicit combinatorial rule for the action of the multiplication operators $e_r^{\bullet}$ generalizing the Pieri rule for symmetric Macdonald functions.}


\Classification{05E05; 05E10; 33D52}

\section{Introduction}
The space of symmetric functions $\Lambda$ is a~central object in algebraic combinatorics deeply connecting the fields of representation theory, geometry, and combinatorics. In his influential paper~\cite{MacDSLC}, Macdonald introduced a~special basis $P_{\lambda}$ for $\Lambda$ over $\mathbb{Q}(q,t)$ simultaneously generalizing many other important and well-studied symmetric function bases like the Schur functions~$s_{\lambda}$. These symmetric functions $P_{\lambda}$, called the symmetric Macdonald functions, exhibit many striking combinatorial properties and can be defined as the eigenvectors of a~certain operator~${\Delta\colon \Lambda \rightarrow \Lambda}$ called the Macdonald operator constructed using polynomial difference operators. It was discovered through the works of Bergeron, Garsia, Haiman, Tesler, and many others~\cite{BGSciFi,BGHT, haiman2000hilbert} that variants of the symmetric Macdonald functions called the modified Macdonald functions~$\widetilde{H}_{\lambda}(x;q,t)$ have deep ties to the geometry of the Hilbert schemes $\Hilb_n\bigl(\mathbb{C}^2\bigr)$. On the side of representation theory, it was shown first in full generality by Cherednik~\cite{C_2001} that one can recover the symmetric Macdonald functions by considering the representation theory of certain algebras called the spherical double affine Hecke algebras (DAHAs) in type ${\rm GL}_n$.

The elliptic Hall algebra (EHA), $\sE$, was introduced by Burban and Schiffmann~\cite{BS} through their study of the Hall algebra of the category of coherent sheaves on an elliptic curve over a~finite field. This algebra has connections to many areas of mathematics including, most importantly for the present paper, to Macdonald theory. The algebra $\sE$ may be described explicitly by generators~$P_{r,\ell}$ indexed by non-zero lattice points $(r,\ell)$ and relations \cite[Theorem~5.4]{BS}. The positive elliptic Hall algebra $\sE^{+}$, is the subalgebra generated by those elements $P_{r,\ell}$ where either~${r\geq 1}$ or $r=0$ and $\ell \neq 0$. In~\cite{SV2}, Schiffmann and Vasserot realize $\mathcal{E}^{+}$ as a~limit of the positive spherical DAHAs in type ${\rm GL}_n$.
 They show further that there is a~natural action of~$\mathcal{E}^{+}$ on $\Lambda$ recovering the spherical DAHA representations originally considered by Cherednik. In particular, the action of~$P_{0,1} \in \sE^{+}$ gives the Macdonald operator $\Delta$. The action of~$\sE^{+}$ on $\Lambda$ can be realized as the action of certain generalized convolution operators on the torus equivariant $K$-theory of the schemes $\Hilb_n\bigl(\mathbb{C}^2\bigr)$.

Dunkl and Luque in~\cite{DL_2011} introduced symmetric and non-symmetric vector-valued (vv) Macdonald polynomials. The term vector-valued here refers to polynomial-like objects of the form~${\sum_{\alpha}c_{\alpha}X^{\alpha}\otimes v_{\alpha}}$ for some scalars $c_{\alpha}$, monomials $X^{\alpha}$, and vectors $v_{\alpha}$ lying in some $\mathbb{Q}(q,t)$-vector space. The non-symmetric vv Macdonald polynomials are distinguished bases for certain DAHA representations built from the irreducible representations of the finite Hecke algebras in type A. These DAHA representations are indexed by Young diagrams and exhibit interesting combinatorial properties relating to periodic Young tableaux. The symmetric vv Macdonald polynomials are distinguished bases for the spherical (i.e., Hecke-invariant) subspaces of these DAHA representations. Naturally, the spherical DAHA acts on this spherical subspace with the special element $\xi_1+\dots + \xi_n$ of spherical DAHA acting diagonally on the symmetric vv Macdonald polynomials.

Dunkl and Luque in~\cite{DL_2011} (and in later work of Colmenarejo, Dunkl, and Luque~\cite{CDL_2022} and Dunkl~\cite{D_2019}) only consider the finite rank non-symmetric and symmetric vv Macdonald polynomials. It is natural to ask if there is an infinite-rank limit construction using the symmetric vv Macdonald polynomials to give generalized symmetric Macdonald functions and an associated representation of the positive elliptic Hall algebra $\mathcal{E}^{+}$. In this paper, we describe such a~construction (see Theorem~\ref{positive EHA module}).
\begin{Theorem}
 There exists a~family \smash{$\bigl(\widetilde{W}_{\lambda}\bigr)_{\lambda \in \Par}$} of graded mutually non-isomorphic cyclic $\mathcal{E}^+$-modules such that \smash{$\widetilde{W}_{\varnothing} \cong \Lambda_{q,t}$}.
\end{Theorem}
The spaces \smash{$\widetilde{W}_{\lambda}$} are limits of the spherical representations associated to the vector-valued DAHA representations $V_{\lambda}$ of Dunkl--Luque. The $V_{\lambda}$ are obtained by inflating irreducible Hecke--Specht modules (labeled by integer partitions) for the finite Hecke algebra to the affine Hecke algebra and then inducing to DAHA. Inside of each $V_{\lambda}$ is the spherical submodule $W_{\lambda}$ on which the spherical DAHA acts. We achieve our construction of the $\mathcal{E}^{+}$ representations \smash{$\widetilde{W}_{\lambda}$} by carefully choosing projection maps $W_{\lambda} \rightarrow W_{\lambda'}$ which are equivariant for the corresponding actions of spherical DAHAs (see Proposition~\ref{relations for connecting maps}). For combinatorial reasons, there is essentially a~unique natural way to obtain this construction. For any $\lambda$, we consider the increasing chains of Young diagrams $\lambda^{(n)} = (n-|\lambda|,\lambda)$ for $n \geq |\lambda| + \lambda_1$ to build the representations \smash{$\widetilde{W}_{\lambda}$}. These special sequences of Young diagrams are central to Murnaghan's theorem~\cite{M_1938} regarding the reduced Kronecker coefficients. As such, we refer to the $\mathcal{E}^{+}$-representations \smash{$\widetilde{W}_{\lambda}$} as Murnaghan-type.

Along with obtaining a~new family of graded $\mathcal{E}^{+}$-representations \smash{$\widetilde{W}_{\lambda}$} indexed by Young diagrams $\lambda$, we construct natural generalizations of the symmetric Macdonald functions $\MacD_{T}$ built as limits of the symmetric vv Macdonald polynomials indexed by certain labelings of infinite Young diagrams (see Definition~\ref{stable limit Definition} and Lemma~\ref{action of MacD operator}).

\begin{Theorem}
 Let $\lambda \in \Par$. There exists an explicit $\Delta$-weight basis $(\MacD_{T})_{T \in \Omega(\lambda)}$ of \smash{$\widetilde{W}_{\lambda}$}. In the case $\lambda = \varnothing$, the indexing set $\Omega(\varnothing)$ can be identified with $\Par$ and for all $\mu \in \Par$, $\MacD_{\mu}$ and $P_{\mu}\bigl(x;q^{-1},t\bigr)$ differ by a~non-zero scalar.
\end{Theorem}

We emphasize that we are interested in investigated both the representations \smash{$\widetilde{W}_{\lambda}$} and the bases \smash{$(\MacD_{T})_{T \in \Omega(\lambda)}$}. For this reason, we develop a~combinatorial understanding of the bases \smash{$(F_{\tau})_{\tau \in \PSYT(\lambda)}$}, \smash{$ (P_{T})_{T\in \RSSYT(\lambda)}$}, and $(\MacD_{T})_{T \in \Omega(\lambda)}$ while simultaneously constructing the modules $V_{\lambda}$, $W_{\lambda}$, and \smash{$\widetilde{W}_{\lambda}$}, respectively. Our goal is, in part, to show that the spaces \smash{$\widetilde{W}_{\lambda}$} along with the bases $(\MacD_{T})_{T \in \Omega(\lambda)}$ constitute a~natural generalization of the usual Macdonald theory picture from both a~representation-theoretic and a~combinatorial perspective. As evidence to our goal, we obtain a~Pieri rule for the action of the multiplication operators $e_r^{\bullet}$ on the generalized symmetric Macdonald function basis $\MacD_{T}$ (see Theorem~\ref{Pieri Rule for Finite vars} and Corollary~\ref{Pieri Rule}).

\begin{Theorem}
 Let $\lambda \in \Par$ and $r\geq 1$. There is an explicit combinatorial formula for the coefficients \smash{$\mathfrak{d}^{(r)}_{S,T}$} in the expansions
 \[
 e_r^{\bullet}\MacD_T = \sum_{S \in \Omega(\lambda)} \mathfrak{d}^{(r)}_{S,T} \MacD_{S}.
 \]

\end{Theorem}
After studying the particular case of the $e_1$-Pieri coefficients, we show that the modules \smash{$\widetilde{W}_{\lambda}$} are cyclically generated by their unique elements of minimal degree $\MacD_{T_{\lambda}^{\min}}$. This demonstrates that the combinatorial understanding for the modules \smash{$\widetilde{W}_{\lambda}$} developed throughout this paper is rich enough to prove non-trivial results that recover well-known results in the $\lambda = \varnothing$ case. The existence of these representations of the elliptic Hall algebra raises many questions about possible new relations between Macdonald theory and geometry. Other authors have constructed families of $\sE^{+}$ representations~\cite{FFJMM_2011,FJMM_2012}. Although there should exist a~relationship between the Murnaghan-type representations \smash{$\widetilde{W}_{\lambda}$} and those of other authors, the construction in this paper appears to be distinct from prior $\sE^{+}$-module constructions.

For technical reasons (regarding the misalignment of the spectrum of the Cherednik operators~$\xi_i$), we need to reprove many of the results of Dunkl and Luque in~\cite{DL_2011} using a~re-oriented version of the Cherednik operators $\theta_i$. Since the elements~$\theta_i$ are not uniformly conjugate to the~$\xi_i$ on the vector-valued polynomial spaces $V_{\lambda}$, we are not immediately able to use the results of Dunkl and Luque. This alternative choice of conventions greatly assists during the construction of the generalized Macdonald functions $\MacD_{T}$. The $\theta_i$ satisfy additional stability properties which the~$\xi_i$ fail to satisfy. The combinatorial model underpinning the non-symmetric vv Macdonald polynomials originally defined by Dunkl and Luque is similar but in this paper, we opt for a~tableau model for our combinatorics in order to agree with the usual Macdonald theory case. Our tableau model is well behaved when we limit from the finite-rank spherical DAHA representations $W_{\lambda}$ to the infinite-rank positive EHA representations \smash{$\widetilde{W}_{\lambda}$}.

Here we give a~brief overview of this paper. First, in Section~\ref{defs and nots}, we review relevant definitions and notations, define our combinatorial models as well as recall the limit spherical DAHA construction of Schiffmann--Vasserot. In Section~\ref{DAHA Modules from Young Diagrams}, we reprove many of the results of Dunkl--Luque but for the re-oriented Cherednik operators including describing the non-symmetric vv Macdonald polynomials $F_{\tau}$ and their associated Knop--Sahi relations (see Proposition~\ref{weight basis Proposition}). We define (see Definition~\ref{connecting maps def}) the DAHA modules $V_{\lambda^{(n)}}$ and connecting maps \smash{$\Phi^{(n)}_{\lambda}\colon V_{\lambda^{(n+1)}} \rightarrow V_{\lambda^{(n)}}$} which are used in the limiting process. Next in Section~\ref{Positive EHA Representations from Young Diagrams}, we describe the spherical subspaces~\smash{$W_{\lambda}^{(n)}$} of Hecke invariants of~\smash{$V_{\lambda}^{(n)}$} and the symmetric vv Macdonald polynomials $P_{T}$ including an explicit expansion of the $P_{T}$ into the $F_{\tau}$ (see Corollary~\ref{expansion of sym into nonsym}). Despite the non-symmetric vv Macdonald polynomials appearing in this paper not agreeing with those of Dunkl--Luque, once Hecke-symmetrized, they do agree (see Corollary~\ref{Corollary recovering sym vv Macd}). We use the connecting maps to define the limit spaces \smash{$\widetilde{W}_{\lambda}$} and show in Theorem~\ref{positive EHA module} that they possess a~graded action of $\mathcal{E}^{+}$ having a~distinguished basis of generalized symmetric Macdonald functions $\MacD_{T}$. In Section~\ref{Pieri Rule section}, we obtain a~Pieri formula (see Corollary~\ref{Pieri Rule}) for the action of $e_r^{\bullet}$ on the generalized Macdonald functions $\MacD$. We apply this rule in the $r=1$ situation to combinatorially verify that the $e_1$-Pieri coefficients are non-vanishing (see Corollary~\ref{non-vanishing of e_1 Pieri for limit MacD}) which shows that the $\mathcal{E}^+$-modules \smash{$\widetilde{W}_{\lambda}$} are cyclic.

The construction of the spaces \smash{$\widetilde{W}_{\lambda}$} is highly related to the Carlsson--Mellit algebra $\mathbb{A}_{q,t}$ \cite{CM_2015} and its variant $\mathbb{B}_{q,t}$ defined by Carlsson--Gorsky--Mellit~\cite{CGM_2017}. In the author's preprint~\cite{DDPA_BW}, it is shown that any sequence of DAHA modules satisfying conditions analogous to those in Proposition~\ref{relations for connecting maps} yields a~corresponding representation of $\mathbb{B}_{q,t}$ and thus of $\mathbb{A}_{q,t}$. For instance, that construction applies to the spaces $V_{\lambda^{(n)}}$ to yield a~family of $\mathbb{B}_{q,t}$-modules indexed by partitions~$\lambda$. The~${k=0}$ (i.e., fully spherical) component of these representations are just the \smash{$\widetilde{W}_{\lambda}$} in the current work. In the~${\lambda = \varnothing}$ case, this recovers the $\mathbb{A}_{q,t}$-module in Carlsson--Mellit's proof of the shuffle theorem~\cite{CM_2015}.

\section{Definitions and notations}\label{defs and nots}

\subsection{Some combinatorics}

We start with a~description of many of the combinatorial objects which we need for the remainder of this paper.

\begin{Definition}\label{Tableaux defs}
 A~\textit{partition} is a~(possibly empty) sequence of weakly decreasing positive integers. Denote by $\Par$ the set of all partitions. Given a~partition $\lambda = (\lambda_1,\dots, \lambda_r)$, we set $\ell(\lambda) := r$ and $|\lambda| := \lambda_1 + \cdots + \lambda_r$. For $\lambda = (\lambda_1,\dots, \lambda_r) \in \Par$ and $n \geq n_{\lambda}:= |\lambda| + \lambda_1$, we set $\lambda^{(n)}:= (n -|\lambda|, \lambda_1,\dots, \lambda_r)$. We identify partitions as defined above with \textit{Young diagrams} of the corresponding shape in English notation, i.e., justified up and to the left.

 Fix a~partition $\lambda$ with $|\lambda| = n$. We require each of the following combinatorial constructions for types of labeling of the Young diagram $\lambda$. If a~diagram $\lambda$ appears as the domain of a~labeling function, then we are referring to the set of boxes of $\lambda$ as the domain.
\begin{itemize}\itemsep=0pt
 \item A~non-negative \textit{reverse Young tableau} $\RYT(\lambda)$ is a~labeling $T\colon \lambda \rightarrow \mathbb{Z}_{\geq 0}$ which is weakly decreasing left to right across rows and top to bottom down columns.
 \item A~non-negative \textit{reverse semi-standard Young tableau} $\RSSYT(\lambda)$ is a~labeling $T\colon \lambda \rightarrow \mathbb{Z}_{\geq 0}$ which is weakly decreasing left to right across rows and strictly decreasing top to bottom down columns.
 \item A~\textit{standard Young tableau} $\SYT(\lambda)$ is a~labeling $\tau\colon \lambda \rightarrow \{1,\dots,n\}$ which is strictly increasing left to right across rows and top to bottom down columns.
 \item A~non-negative \textit{periodic standard Young tableau} $\PSYT(\lambda)$ is a~labeling $\tau\colon \lambda \rightarrow \big\{ jq^{b}\mid 1\leq j \leq n,\, b \geq 0\big\}$ in which each $1\leq j \leq n$ occurs in exactly one box of $\lambda$ and where the labeling is strictly increasing left to right across rows and top to bottom down columns. Here we order the formal products $jq^m$ by $jq^m < kq^{\ell}$ if $m > \ell$ or in the case that $m = \ell$ we have $j < k$. Note that $\SYT(\lambda) \subset \PSYT(\lambda)$.
\end{itemize}

Define $\tau^{\rm rs}_{\lambda}, \tau^{\rm cs}_{\lambda} \in \SYT(\lambda)$ to be the \textit{row-standard} and \textit{column-standard} labelings of $\lambda$, respectively. Namely, $\tau^{\rm rs}_{\lambda}$ is the unique labeling of $\lambda$ which, starting from the top left box, reads off each row left to right before moving onto the next row below. Similarly, $\tau^{\rm cs}_{\lambda}$ is the unique labeling of $\lambda$ which, starting from the top left box, reads off each column top to bottom before moving onto the next column to the right.

\end{Definition}

\begin{Example}\quad
\[
\ytableausetup{centertableaux, boxframe= normal, boxsize= 2.25em}
\begin{ytableau}
 17q^7 & 15q^5 & 16q^5 & 11q^3 & 7q^1 & 2q^0 \\
 14q^6 & 12q^4 & 13q^4 & 9q^2 & 8q^0 & \none \\
 10q^2 & 4q^1 & 5q^1 & 6q^1 & \none & \none \\
 3q^1 & 1q^0 & \none & \none & \none & \none \\
\end{ytableau} \in \PSYT(6,5,4,2).
\]
\end{Example}

\begin{Definition}\label{ordering and operations on tableaux defs}
Given a~box $\square$ in a~Young diagram $\lambda$ we define the \textit{content} of $\square$ as $c(\square) := a-b$ where $\square = (a,b)$ as drawn in the $\mathbb{N}\times \mathbb{N}$ grid (English notation). Let $\tau \in \PSYT(\lambda)$ and $1\leq i\leq n$. Whenever $\tau(\square) = iq^{m}$ for some box $\square \in \lambda$, we write
\begin{itemize}
\itemsep=0pt
 \item $c_{\tau}(i):= c(\square)$,
 \item $w_{\tau}(i):= m$.
\end{itemize}

We define the operations $s_j$ and $\Psi$ on the set $\PSYT(\lambda)$ as follows. Set $w_{\tau}:= (w_{\tau}(1),\dots,\allowbreak w_{\tau}(n)) \in \mathbb{Z}^{n}_{\geq0}$. Let $1\leq j \leq n-1$ and suppose that for some boxes $\square_1,\square_2 \in \lambda$, $\tau(\square_1) = jq^m$ and $\tau(\square_2) = (j+1)q^{\ell}$.
\begin{itemize}
\itemsep=0pt
 \item Let $\tau'$ be the labeling defined by $\tau'(\square_1) = (j+1)q^m$, $\tau'(\square_2) = jq^{\ell}$, and $\tau'(\square) = \tau(\square)$ for~${\square \in \lambda \setminus \{\square_1,\square_2\}}$. If $\tau' \in \PSYT(\lambda)$, then we write $s_j(\tau):= \tau'$.
 \item Let $\Psi(\tau) \in \PSYT(\lambda)$ be the labeling defined by whenever $\tau(\square) = kq^a$, then either $\Psi(\tau)(\square) = (k-1)q^a$ when $k \geq 2$ or $\Psi(\tau)(\square) = nq^{a+1}$ when $k = 1$.
\end{itemize}

We give the set $\PSYT(\lambda)$ the partial order $\geq$ defined by the following cover relations:
\begin{itemize}
\itemsep=0pt
 \item For all $\tau \in \PSYT(\lambda)$, $\Psi(\tau) > \tau$.
 \item If $w_{\tau}(i)<w_{\tau}(i+1)$, then $s_i(\tau) > \tau$.
 \item If $w_{\tau}(i) = w_{\tau}(i+1)$ and $c_{\tau}(i)-c_{\tau}(i+1) > 1$, then $s_i(\tau) > \tau$.
\end{itemize}

Lastly, define the map $\mathfrak{p}_{\lambda}\colon \PSYT(\lambda) \rightarrow \RYT(\lambda)$ by $\mathfrak{p}_{\lambda}(\tau)(\square) = b$ whenever $\tau(\square) = iq^b$. We write $\PSYT(\lambda;T)$ for the set of all $\tau \in \PSYT(\lambda)$ with $\mathfrak{p}_{\lambda}(\tau) = T \in \RYT(\lambda)$.

\end{Definition}

Note that the operation $\Psi\colon \PSYT(\lambda) \rightarrow \PSYT(\lambda)$ is indeed well defined. If $\tau \in \PSYT(\lambda)$ and $\tau(\square_1) < \tau(\square_2)$, then $\Psi(\tau)(\square_1) < \Psi(\tau)(\square_2)$ since the operation defined by $kq^a \mapsto (k-1)q^a$ when $k \geq 2$ and $kq^a \mapsto nq^{a+1}$ when $k = 1$ is strictly decreasing for the ordering on formal products defined in Definition~\ref{Tableaux defs}.

\begin{Example}\quad
\[\Psi \left( \ytableausetup{centertableaux, boxframe= normal, boxsize= 2.25em}
\begin{ytableau}
 1q^7& 3q^5 & 5q^5 & 8q^2 & 12q^1 & 17q^0 \\
 2q^6& 4q^5 & 6q^5 & 14q^0 & 16q^0 & \none \\
 7q^2& 10q^1 & 11q^1 & 15q^0 & \none & \none \\
 9q^1& 13q^0 & \none & \none & \none & \none \\
\end{ytableau} \right) = \ytableausetup{centertableaux, boxframe= normal, boxsize= 2.25em}
\begin{ytableau}
 17q^8& 2q^5 & 4q^5 & 7q^2 & 11q^1 & 16q^0 \\
 1q^6& 3q^5 & 5q^5 & 13q^0 & 15q^0 & \none \\
 6q^2& 9q^1 & 10q^1 & 14q^0 & \none & \none \\
 8q^1& 12q^0 & \none & \none & \none & \none \\
\end{ytableau}\]
\end{Example}

We frequently require the next basic lemma regarding the ordering $\leq$ on $\PSYT(\lambda)$.

\begin{Lemma}\label{locally maximal periodic tableau}
 Let $\lambda \in \Par$ and $T \in \RYT(\lambda)$. There are unique $\Min(T), \Top(T) \in \PSYT(\lambda;T)$ such that for all $\tau \in \PSYT(\lambda)$ with $\mathfrak{p}_{\lambda}(\tau) = T$, $\Min(T) \leq \tau \leq \Top(T)$.
\end{Lemma}

\begin{proof}
 We construct the elements $\Top(T)$, $\Min(T)$ directly. Define $\Top(T)$ by first filling in the boxes $\square \in \lambda$ of $\lambda$ with the values $q^{T(\square)}$. Now we label these boxes with the values $\{1,\cdots, n \}$ by first decomposing $\lambda$ into skew diagrams where $T$ is constant on each sub-diagram. This gives us an increasing chain of Young diagrams $\lambda^{1}\subset \dots \subset \lambda^{r} = \lambda$. Next we fill each diagram $\lambda^{i}$ with the values $\big\{\big|\lambda^{1}\big|+\cdots + \big|\lambda^{i-1}\big|+1,\dots, \big|\lambda^{1}\big|+\dots + \big|\lambda^{i}\big| \big\} $ in column-standard order. This gives a~value~$iq^a$ in each box of $\lambda$.

 For $\Min(T)$, we proceed similarly by first filling in the boxes $\square \in \lambda$ of $\lambda$ with the values $q^{T(\square)}$. Then we decompose $\lambda$ into the same skew diagrams as before. Now we fill each diagram $\lambda^{i}$ with the values $\big\{n-\bigl(\big|\lambda^{1}\big|+\cdots + \big|\lambda^{i-1}\big|\bigr),\dots, n-\bigl(\big|\lambda^{1}\big|+\cdots + \big|\lambda^{i}\big|\bigr) \big\} $ in row-standard order. This gives a~value $iq^a$ in each box of $\lambda$.
\end{proof}

\begin{Example}\label{Example of labeling types}\quad
 Given
\[T = \ytableausetup{centertableaux, boxframe= normal, boxsize= 2.25em}
\begin{ytableau}
 7& 5 & 5 & 2 & 1 & 0 \\
 6& 5 & 5 & 0 & 0 & \none \\
 2& 1 & 1 & 0 & \none & \none \\
 1& 0 & \none & \none & \none & \none \\
\end{ytableau} \in \RYT(6,5,4,2),\]
we have that
\[\Min(T) = \ytableausetup{centertableaux, boxframe= normal, boxsize= 2.25em}
\begin{ytableau}
 17q^7& 12q^5 & 13q^5 & 10q^2 & 6q^1 & 1q^0 \\
 16q^6& 14q^5 & 15q^5 & 2q^0 & 3q^0 & \none \\
 11q^2& 7q^1 & 8q^1 & 4q^0 & \none & \none \\
 9q^1& 5q^0 & \none & \none & \none & \none \\
\end{ytableau}\]
and
\[\Top(T) = \ytableausetup{centertableaux, boxframe= normal, boxsize= 2.25em}
\begin{ytableau}
 1q^7& 3q^5 & 5q^5 & 8q^2 & 12q^1 & 17q^0 \\
 2q^6& 4q^5 & 6q^5 & 14q^0 & 16q^0 & \none \\
 7q^2& 10q^1 & 11q^1 & 15q^0 & \none & \none \\
 9q^1& 13q^0 & \none & \none & \none & \none \\
\end{ytableau}.\]

As in the proof of Lemma~\ref{locally maximal periodic tableau}, we correspond to $T$ the following increasing sequence of Young diagrams which fill with the labels of $T$ for reference:
\[
\ytableausetup{centertableaux, boxframe= normal, boxsize= 1em}
\begin{ytableau}
 7\\
\end{ytableau} \subset ~ \begin{ytableau}
 7 \\
 6 \\
\end{ytableau} \subset ~ \begin{ytableau}
 7 & 5& 5\\
 6 & 5 & 5 \\
\end{ytableau} \subset ~ \begin{ytableau}
 7 & 5& 5 & 2\\
 6 & 5 & 5 & \none \\
 2 & \none & \none & \none \\
\end{ytableau} \subset ~ \begin{ytableau}
 7 & 5& 5 & 2 & 1\\
 6 & 5 & 5 & \none & \none \\
 2 & 1 & 1 & \none & \none \\
 1 & \none & \none & \none & \none \\
\end{ytableau} \subset ~ \begin{ytableau}
 7 & 5& 5 & 2 & 1 & 0\\
 6 & 5 & 5 & 0 & 0 & \none \\
 2 & 1 & 1 & 0 & \none & \none \\
 1 & 0 & \none & \none & \none & \none \\
\end{ytableau}.\]

\end{Example}

\begin{Definition}\label{decomposing RYT into SYT and partition}
Let $\lambda \in \Par$ with $|\lambda| = n$ and $T \in \RYT(\lambda)$. Define $\nu(T) \in \mathbb{Z}_{\geq 0}^{n}$ to be the vector formed by listing the values of T in decreasing order, i.e., $\nu(T) = \sort(w_{\tau})$ for any $\tau \in \PSYT(\lambda;T)$. Define $S(T) \in \SYT(\lambda)$ by ordering the boxes of $\lambda$ according to $\square_1 \leq \square_2$ if and only if
\begin{itemize}
\itemsep=0pt
 \item $T(\square_1) > T(\square_2)$ or
 \item $T(\square_1) = T(\square_2)$ and $\square_1$ comes before $\square_2$ in the column-standard labeling of $\lambda$.
\end{itemize}
We often write as a~shorthand $\square_1 <_{T} \square_2$ whenever $S(T)(\square_1) < S(T)(\square_2)$.
Define the statistic~${b_T \in \mathbb{Z}_{\geq 0}}$ by
\[
b_T:= \sum_{i=1}^{n} \nu(T)_i( c_{S(T)}(i) + i-1).
\]

Lastly, define the composition $\mu(T)$ of $n$ as follows. Decompose $\lambda$ into maximal size horizontal strips $h_1,\dots, h_m$ such that $T$ is constant on each strip. We order these strips according to the filling $\Min(T)$, i.e., the $\Min(T)$ labels in $h_i$ are strictly less than those in $h_{i+1}$ for all $i$. We see that each of these horizontal strips $h_i$ has some length $r_i\geq 1$. Then $\mu(T)$ is given as $(r_1,\dots, r_m)$.
\end{Definition}

\begin{Remark}
 For every $T \in \RYT(\lambda)$, we can recover $T$ from the pair $(S(T),\nu(T))$ by labeling~$\lambda$ with the entries of $\nu(T)$ following the order of $S(T)$. Further, the standard Young tableau~$S(T)$ is the largest such tableau following the partial order defined in Definition~\ref{ordering and operations on tableaux defs}.
\end{Remark}

Below is an example calculation of the various data which we associate to $T \in \RYT(\lambda)$.

\begin{Example}
 For $T \in \RYT(6,5,4,2)$ as in Example~\ref{Example of labeling types}, we have that
\begin{gather*}
S(T) = \ytableausetup{centertableaux, boxframe= normal, boxsize= 2.25em}
\begin{ytableau}
 1& 3 & 5 & 8 & 12 & 17 \\
 2& 4 & 6 & 14 & 16 & \none \\
 7& 10 & 11 & 15 & \none & \none \\
 9& 13 & \none & \none & \none & \none \\
\end{ytableau} \in \SYT(6,5,4,2),\\
\nu(T) = (7,6,5,5,5,5,2,2,1,1,1,1,0,0,0,0,0)\in \mathbb{Z}_{\geq 0}^{17},\\
b_T = 0+0+15+15+30+30+8+20+5+8+10+15+0+0+0+0+0 = 156,\\
\text{and}\quad \mu(T) = (1,2,1,1,1,2,1,1,1,2,2,1,1).
\end{gather*}
\end{Example}

The next definition is crucial for many of the results in this paper.

\begin{Definition}\label{inversions Definition}
 Let $\lambda \in \Par$, with $|\lambda|= n$ and $\tau \in \PSYT(\lambda)$ with $T = \mathfrak{p}_{\lambda}(\tau)$. An ordered pair of boxes $(\square_1,\square_2) \in \lambda \times \lambda$ is called an \textit{inversion pair} of $\tau$ if $S(T)(\square_1) < S(T)(\square_2)$ and $i > j$ where $\tau(\square_1) = iq^a$, $\tau(\square_2) = jq^b$ for some $a,b \geq 0$.
 The set of all inversion pairs of $\tau$ is denoted by $\Inv(\tau)$. We use the shorthand $\I(T)$ for the set $\Inv(\Min(T))$.
\end{Definition}

\begin{Example}
 In the labeling
\[ \ytableausetup{centertableaux, boxframe= normal, boxsize= 2.25em}
\begin{ytableau}
 17q^7& 12q^5 & 13q^5 & 10q^2 & 6q^1 & 1q^0 \\
 16q^6& 14q^5 & 15q^5 & 2q^0 & 3q^0 & \none \\
 11q^2& 7q^1 & 8q^1 & 4q^0 & \none & \none \\
 9q^1& 5q^0 & \none & \none & \none & \none \\
\end{ytableau}\]
we have that the pairs $\bigl(17q^7, 12q^5\bigr)$, $\bigl(14q^5,13q^5\bigr)$, and $\bigl(5q^0,4q^0\bigr)$ are all inversions. Here we have referred to boxes according to their labels.
\end{Example}

In the following definition, our conventions for the Bruhat ordering differ from many other authors. These conventions are used to help properly state some triangularity properties later in the paper. However, one may obtain the below definition from the more standard conventions in~\cite[Section~2.1]{haglund2007combinatorial} by reversing the order of the entries of each vector $(a_1,\dots, a_n) \rightarrow (a_n,\dots, a_1)$ and rewriting their Bruhat ordering from this reversed perspective.

\begin{Definition}\label{Bruhat def}
Define the Bruhat ordering $\preceq$ on $\mathbb{Z}_{\geq 0}^n$ using the following cover relations for~${\lambda \in \mathbb{Z}_{\geq 0}^{n}}$:
\begin{itemize}
\itemsep=0pt
 \item If $i<j$ with $\lambda_i < \lambda_j$, then $\lambda \prec (i,j)\lambda$.
 \item If $i<j$ with $\lambda_i + 1< \lambda_j$, then $\lambda \succ \lambda + e_i - e_j$.
\end{itemize}
 Here $e_i$ denotes the $i$-th standard basis vector of $\mathbb{Z}^{n}$ and $(i,j) \in \mathfrak{S}_n$ denotes the permutation swapping $i$ and $j$ and leaving all other inputs fixed. For $\alpha = (\alpha_1,\alpha_2,\dots, \alpha_n) \in \mathbb{Z}^{n}_{\geq 0}$, we define~${\widetilde{\gamma}(\alpha):= (\alpha_2,\dots, \alpha_n, \alpha_1 +1)}$. We write $\sort(\alpha)$ for the vector formed by listing the entries of $\alpha$ in weakly decreasing order and $\revsort(\alpha)$ for the vector formed by listing the entries of~$\alpha$ in weakly increasing order. We define $\Stab(\alpha)$ to be the corresponding stabilizer subgroup of~$\mathfrak{S}_n$ for $\alpha$, i.e., the set of all $\sigma \in \mathfrak{S}_n$ with $\sigma(\alpha) = \alpha$.
\end{Definition}

We require the following simple lemma regarding the interplay between the map $\widetilde{\gamma}$ on $\mathbb{Z}^n_{\geq 0}$ and the ordering $\prec$.

\begin{Lemma}\label{gamma preserves Bruhat}
 If $\alpha, \beta \in \mathbb{Z}^n_{\geq 0}$ satisfy $\alpha \prec \beta$, then $\widetilde{\gamma}(\alpha) \prec \widetilde{\gamma}(\beta)$.
\end{Lemma}
\begin{proof}
 We show that if $\alpha, \beta \in \mathbb{Z}_{\geq 0}^{n}$ and $\beta$ covers $\alpha$ with respect to the Bruhat order, then $\widetilde{\gamma}(\alpha) \prec \widetilde{\gamma}(\beta)$. We proceed in cases. Let $\lambda \in \mathbb{Z}_{\geq 0}^{n}$.

 First, suppose $1<i<j$ and $\lambda_i < \lambda_j$. Then
\[\widetilde{\gamma}(\lambda) \prec (i-1,j-1)\widetilde{\gamma}(\lambda) = \widetilde{\gamma}( (i,j) \lambda).\]

 Now suppose $1<j$ and $\lambda_1 < \lambda_j$. Then
 \[
 \widetilde{\gamma}( (1,j)\lambda) \succ \widetilde{\gamma}( (1,j)\lambda) + e_j - e_n \succeq (j,n)(\widetilde{\gamma}( (1,j)\lambda) + e_j - e_n) = \widetilde{\gamma}(\lambda).
 \]

 If now we have that $1<i<j$ and $\lambda_i < \lambda_j -1$, then
 \[
 \widetilde{\gamma}(\lambda) \succ \widetilde{\gamma}(\lambda) + e_{i-1}-e_{j-1} = \widetilde{\gamma}(\lambda + e_{i}-e_{j}).
 \]

 Lastly, consider the case when $1 < j$ and $\lambda_1 < \lambda_j -1$. If $\lambda_1 + 2 = \lambda_j$, then
 \[
 \widetilde{\gamma}(\lambda) \succ (j-1,n)\widetilde{\gamma}(\lambda) = \widetilde{\gamma}(\lambda + e_1 - e_j).
 \]
 Instead if $\lambda_1 < \lambda_j -2$, then
 \[
 \widetilde{\gamma}(\lambda) \succ (j-1,n)\widetilde{\gamma}(\lambda) \succ (j-1,n) \widetilde{\gamma}(\lambda) + e_{j-1}-e_n = \widetilde{\gamma}(\lambda+e_1-e_j).\tag*{\qed}
 \] \renewcommand{\qed}{}
\end{proof}

Here we review some necessary details about the extended affine symmetric groups.

\begin{Definition}
 Define $\widehat{\mathfrak{S}}_n$ to be the extended affine symmetric group given by
$
 \widehat{\mathfrak{S}}_n:= \mathfrak{S}_n \ltimes \mathbb{Z}^{n}
$
 where $\mathfrak{S}_n$ acts on $\mathbb{Z}^{n}$ by coordinate permutations. Denote by $t_1,\dots, t_n$ the standard generators of $\mathbb{Z}^{n} \subset \widehat{\mathfrak{S}}_n$. Further, we define the special element $\widetilde{\gamma}_n \in \widehat{\mathfrak{S}}_n$ given by
$
 \widetilde{\gamma}_n:= t_ns_{n-1}\cdots s_1$.
 For any $\beta \in \mathbb{Z}^{n}$, we write
$
 t_{\beta}:= t_1^{\beta_1}\cdots t_n^{\beta_n}$.
 Define the positive submonoid of $\widehat{\mathfrak{S}}_n$, $\widehat{\mathfrak{S}}_n^{+}$, as the monoid generated by $\big\{s_1,\dots, s_{n-1},\widetilde{\gamma}_n \big\}$ (i.e., no $\widetilde{\gamma}_n^{-1}$s).

The length $\ell(\sigma)$ of $\sigma \in \widehat{\mathfrak{S}}_n$ is the minimal number of $s_i$ required to express $\sigma$ in terms of the generators $\big\{s_1,\dots,s_{n-1},\widetilde{\gamma}_n\big\}$.
In particular, multiplication by $\widetilde{\gamma}_n$ does not effect the length of elements of \smash{$\widehat{\mathfrak{S}}_n$}. An expression of the form $\sigma = w_1\cdots w_{r}$ for $w_i \in \{s_1,\dots,s_{n-1},\widetilde{\gamma}_n \}$ is a~\textit{reduced word} for $\sigma$ if $\#( i \mid w_i \in \{s_1,\dots,s_{n-1}\}) = \ell(\sigma)$.
 We denote by \smash{$\widehat{\mathfrak{S}}_n/\mathfrak{S}_n$} the set of minimal length left coset representatives of \smash{$\widehat{\mathfrak{S}}_n$} with respect to the subgroup $\mathfrak{S}_n$.
 We denote the set of positive minimal length coset representatives of $\widehat{\mathfrak{S}}_n$ with respect to the subgroup $\mathfrak{S}_n$ by~\smash{$\bigl(\widehat{\mathfrak{S}}_n/\mathfrak{S}_n \bigr)^{+}:= \bigl(\widehat{\mathfrak{S}}_n/\mathfrak{S}_n\bigr) \cap \widehat{\mathfrak{S}}_n^{+}$}.
If $\mu = (\mu_1,\dots, \mu_r)$ is a~composition of $n = \mu_1+\cdots + \mu_r$, then we define the Young subgroup $\mathfrak{S}_{\mu}$ of $\mathfrak{S}_n$ corresponding to $\mu$ as $\mathfrak{S}_{\mu}:= \mathfrak{S}_{\mu_1}\times \cdots \times \mathfrak{S}_{\mu_r} \subset \mathfrak{S}_n$. We write $\mathfrak{S}_n/\mathfrak{S}_{\mu}$ for the set of minimal length left coset representatives for $\mathfrak{S}_n$ with respect to the subgroup $\mathfrak{S}_{\mu}$.

 For $\beta \in \mathbb{Z}^n$, define $\sigma_{\beta} \in \widehat{\mathfrak{S}}_n$ by
$
 \sigma_{\beta}:= \sigma t_{\sort(\beta)}
$
 where $\sigma$ is the unique minimal length coset representative in $\mathfrak{S}_n/\mathfrak{S}_{\Stab(\sort(\beta))}$ such that $\sigma(\sort(\beta))= \beta$.

\end{Definition}

The next lemma summarizes some standard results in the theory of (extended) affine permutations.

\begin{Lemma}\label{min coset reps Lemma}\label{properties of min length rep}
 We have that
$
 \widehat{\mathfrak{S}}_n/\mathfrak{S}_n= \{ \sigma_{\beta} \mid \beta \in \mathbb{Z}^n \}
$
 and
$
 \bigl(\widehat{\mathfrak{S}}_n/\mathfrak{S}_n \bigr)^{+}= \{ \sigma_{\beta} \mid \beta \in \mathbb{Z}^n_{\geq 0} \}$.
 For all $\alpha \in \mathbb{Z}^{n}_{\geq 0}$, we have the following:
 \begin{itemize}\itemsep=0pt
 \item If $\alpha$ is weakly decreasing, then $\sigma_{\alpha} = t_{\alpha}$.
 \item If $s_i(\alpha) \succ \alpha$, then $\sigma_{s_i(\alpha)}= s_i\sigma_{\alpha}$.
 \item If $s_i(\alpha) = \alpha$, then $s_i\sigma_{\alpha} = \sigma_{\alpha} s_{\sigma^{-1}(i)}$, where $\sigma$ is the minimal length permutation with $\sigma(\sort(\alpha))= \alpha$.
 \item $\sigma_{\widetilde{\gamma}_n(\alpha)} = \widetilde{\gamma}_n (\sigma_{\alpha})$.
 \end{itemize}
\end{Lemma}
\begin{proof}
 This follows directly from \cite[Proposition~2.1]{Stembridge} applied to the setting of $\widehat{\mathfrak{S}}_n$ along with using Lemma~\ref{gamma preserves Bruhat} for the last statement.
\end{proof}

Recall that in Definition~\ref{ordering and operations on tableaux defs} we only defined $s_i(\tau)$ for $\tau \in \PSYT(\lambda)$ in the situation where swapping the $i$ and $i+1$ labels in the boxes of $\tau$ resulted in an element of $\PSYT(\lambda)$. We now generalize this notion to elements of \smash{$\widehat{\mathfrak{S}}_n^{+}$}.

\begin{Definition}\label{action of extended affine permutations of PSYT}
 Suppose $z_r\cdots z_1$ is a~reduced word in $\widehat{\mathfrak{S}}_n^{+}$ written in the generators $z_i \in \{s_1,\dots, s_{n-1},\widetilde{\gamma}_n\}$. We define inductively on $r \geq 1$ if $z_{r-1}\cdots z_1(\tau) \in \PSYT(\lambda)$ the element $z_r\cdots z_1(\tau)$ of $\PSYT(\lambda)$ as either
 \begin{itemize}
\itemsep=0pt
 \item $\Psi(z_{r-1}\cdots z_1(\tau))$ if $z_r = \widetilde{\gamma}_n$, or
 \item $s_i(z_{r-1}\cdots z_1(\tau))$ if $z_r = s_i$ and swapping the $i$ and $i+1$ labels in the boxes of $z_{r-1}\cdots z_1(\tau)$ results in an element of $\PSYT(\lambda)$.
 \end{itemize}
 Otherwise, we leave this symbol undefined. This definition is only dependent on the element~${z_r\cdots z_1}$ of $\widehat{\mathfrak{S}}_n^{+}$ in that if $z_r\cdots z_1 = z_r'\cdots z_1'$ is another reduced word, then $z_r\cdots z_1(\tau)$ is defined if and only if $z_r'\cdots z_1'(\tau)$ is defined in which case $z_r\cdots z_1(\tau) = z_r'\cdots z_1'(\tau)$. Thus, we write $\sigma(\tau) = z_r\cdots z_1(\tau)$ unambiguously in this situation whenever $\sigma = z_r\cdots z_1$.
\end{Definition}

Implicit in the above definition is that the following relations hold.

\begin{Lemma}
For any $\tau \in \PSYT(\lambda)$ in each of the following equalities if one of either the left-hand or right-hand labelings define valid elements of $\PSYT(\lambda)$ then both sides are elements of $\PSYT(\lambda)$ and they are equal:
\begin{itemize}
\itemsep=0pt
 \item For all $1\leq i \leq n-1$, $s_i(s_i(\tau)) = \tau$.
 \item If $|i-j|>1$, then $s_is_j(\tau) = s_js_i(\tau)$.
 \item For all $1\leq i \leq n-2$, $s_is_{i+1}s_i(\tau) = s_{i+1}s_is_{i+1}(\tau)$.
 \item For all $1\leq i \leq n-2$, $\Psi s_{i+1}(\tau) = s_i\Psi(\tau)$.
 \item $\Psi^2 s_1(\tau) = s_{n-1}\Psi^2(\tau)$.
\end{itemize}

\end{Lemma}
\begin{proof}
Let $\tau \in \PSYT(\lambda)$. We prove each relation separately.

If $s_i(\tau) \in \PSYT(\lambda)$, then either $w_{\tau}(i) \neq w_{\tau}(i+1)$ or otherwise $i$ and $i+1$ are distinct rows and columns. The same is true about $s_i(\tau)$ so $s_i(s_i(\tau)) \in \PSYT(\lambda)$ which is obviously equal to $\tau$.

If $s_is_j(\tau) \in \PSYT(\lambda)$, then we know that either $w_{\tau}(j) \neq w_{\tau}(j+1)$ or $j$ and $j+1$ are in distinct rows and columns and likewise for $i$ and $i+1$ since $|i-j| > 1$. Therefore, $s_js_i(\tau) \in \PSYT(\lambda)$ is also well defined and both are equal because $s_is_j = s_js_i$ as permutations.

Let $\square_i$, $\square_{i+1}$, $\square_{i+2}$ denote the boxes of $\lambda$ labeled $iq^{a_i}$, $(i+1)q^{a_{i+1}}$, $(i+2)q^{a_{i+2}}$ by $\tau$ respectively for some $a_i,a_{i+1},a_{i+2} \geq 0$. The condition guaranteeing $s_is_{i+1}s_i(\tau) \in \PSYT(\lambda)$ is the same as the one guaranteeing $s_{i+1}s_is_{i+1}(\tau) \in \PSYT(\lambda)$: If any pair $\square_j$, $\square_{j'}$ of the boxes $\square_i$, $\square_{i+1}$, $\square_{i+2}$ are in the same row or column, then necessarily $a_{j} \neq a_{j'}$. Therefore, $s_is_{i+1}s_i(\tau)\in \PSYT(\lambda)$ if and only if $s_{i+1}s_is_{i+1}(\tau) \in \PSYT(\lambda)$ and both are obviously equal in this case since $s_{i}s_{i+1}s_{i} = s_{i+1}s_is_{i+1}$ as permutations.

Suppose $\Psi s_{i+1}(\tau) \in \PSYT(\lambda)$. Then either $w_{\tau}(i+1) \neq w_{\tau}(i+2)$ or else $i+1$ and $i+2$ appear in distinct rows and columns in the labeling $\tau$. But, then either
\[w_{\Psi(\tau)}(i) = w_{\tau}(i+1) \neq w_{\tau}(i+2)= w_{\Psi(\tau)}(i+1)
\]
 or else $i$ and $i+1$ appear in distinct rows and columns in the labeling~$\Psi(\tau)$. Therefore, $s_i\Psi(\tau) \in \PSYT(\lambda)$ and by Definition~\ref{ordering and operations on tableaux defs} $\Psi s_{i+1}(\tau) = s_i\Psi(\tau)$. In particular, if $\tau(\square) = jq^a$, then
\[
 \Psi s_{i+1}(\tau)(\square)=s_i\Psi(\tau)(\square) = \begin{cases}
 nq^{a+1}, & j=1 ,\\
 (i+1)q^{a}, & j=i+1,\\
 iq^a ,& j =i+2, \\
 (j-1)q^a, & j\neq 1,i+1,i+2 .
\end{cases}
\]
 Similarly, if we suppose $s_i\Psi(\tau) \in \PSYT(\lambda)$, then by the same argument but in reverse, $\Psi s_{i+1}(\tau) \in \PSYT(\lambda)$ and the expressions agree.

Lastly, suppose $\Psi^{2}s_1(\tau) \in \PSYT(\lambda)$ then $w_{\tau}(1) \neq w_{\tau}(2)$ or else $1$ and $2$ appear in different rows and columns of the labeling $\tau$. But, then either $w_{\Psi^2(\tau)}(n-1) = w_{\tau}(1)+1 \neq w_{\tau}(2) +1= w_{\Psi^2(\tau)}(2)$ or else $n-1$ and $n$ appear in distinct rows and columns in the labeling $\Psi^2(\tau)$. Therefore, $s_{n-1}\Psi^2(\tau) \in \PSYT(\lambda)$ and $\Psi^{2}s_1(\tau)=s_{n-1}\Psi^2(\tau)$. In particular, if $\tau(\square) = jq^a$, then
\[
 \Psi^{2}s_1(\tau)(\square) = s_{n-1}\Psi^2(\tau)(\square) = \begin{cases}
 nq^{a+1}, & j =1 ,\\
 (n-1)q^{a+1} ,& j =2 ,\\
 (j-2)q^a ,& j\neq 1,2 .\\
\end{cases}
\]

Likewise, if we suppose $s_{n-1}\Psi^2(\tau) \in \PSYT(\lambda)$, then by the same argument but in reverse, $\Psi^{2}s_1(\tau) \in \PSYT(\lambda)$ and the expressions agree.
\end{proof}

We need the following result regarding generating elements of $\PSYT(\lambda)$ from $\SYT(\lambda)$ using special creation operators $\zeta_i$.

\begin{Lemma}\label{creation operator Lemma}
 For $T \in \RYT(\lambda)$, we have that
 \smash{$\Top(T) = \zeta_1^{\nu(T)_1 - \nu(T)_2}\cdots \zeta_n^{\nu(T)_{n}}(S(T))$}, where for all $1\leq i \leq n$,
 $\zeta_i:= (s_i\cdots s_{n-1}\Psi)^{i}$.
\end{Lemma}

\begin{proof}
 One may check by direct computation that if $T \in \RYT(\lambda)$ and $1\leq i \leq n$, then $\zeta_i(\Top(T)))$ is well defined according to Definition~\ref{action of extended affine permutations of PSYT} and in particular,
 $\zeta_i(\Top(T)) = \Top(T')$ where $T'(\square) = T(\square) + 1$ for $S(T)(\square) \leq i$ and $T'(\square) = T(\square)$ otherwise. Note that $S(T)= S(T')$ so applying $\zeta_i$ does not change the underlying diagram ordering corresponding to the labeling $T$. Thus, given any $T \in \RYT(\lambda)$ by applying each $\zeta_i$ one at a~time, we see that \smash{$\zeta_1^{\nu(T)_1 - \nu(T)_2}\cdots \zeta_n^{\nu(T)_{n}}(S(T))$} must equal $\Top(T)$.
\end{proof}

To define our combinatorial model for the representations we study in Section~\ref{DAHA Modules from Young Diagrams}, we need to identify an explicit bijection between $\PSYT(\lambda)$ and \smash{$\bigl(\widehat{\mathfrak{S}}_n/\mathfrak{S}_n \bigr)^{+} \times \SYT(\lambda)$}. We already have a~map \smash{$\PSYT(\lambda) \rightarrow \bigl(\widehat{\mathfrak{S}}_n/\mathfrak{S}_n \bigr)^{+}$} given by $\tau \rightarrow \sigma_{w_{\tau}}$. This is not bijective so we use elements of~$\SYT(\lambda)$ to refine this map to yield a~bijection. We now identify the correct choice of $\SYT(\lambda)$ for a~given $\tau \in \PSYT(\lambda)$.

\begin{Definition}\label{tableaux from PSYT def}
 For $\tau \in \PSYT(\lambda)$, we define $S(\tau) \in \SYT(\lambda)$ by the following recursion:
 \begin{itemize}\itemsep=0pt
 \item $S(\Top(T)) := S(T)$ as defined in Definition~\ref{decomposing RYT into SYT and partition}
 \item If $w_{\tau}(i) < w_{\tau}(i+1)$, then $S(s_i(\tau))= S(\tau)$.
 \item $S(\Psi(\tau)) = S(\tau)$.
 \item If $w_{\tau}(i) = w_{\tau}(i+1)$ and $c_{\tau}(i) - c_{\tau}(i+1)>1$, then $S(s_i(\tau))= s_jS(\tau)$, where $j = \sigma^{-1}(i)$ and $\sigma$ is the minimal length permutation with $\sigma(\sort(w_{\tau})) = w_{\tau}$.
 \end{itemize}
\end{Definition}

For any $\tau \in \PSYT(\lambda)$, we can recover $\tau$ from the data $(\sigma_{w_{\tau}}, S(\tau))$ as follows.

\begin{Proposition}\label{generating PSYT from SYT}
 For $\tau \in \PSYT(\lambda)$,
 $\tau = \sigma_{w_{\tau}}(S(\tau))$.
\end{Proposition}
\begin{proof}
 Using Lemmas~\ref{properties of min length rep} and~\ref{creation operator Lemma}, we see that for all $T \in \RYT(\lambda)$
 \begin{align*}
 \sigma_{w_{\Top(T)}}(S(\mathfrak{p}_{\lambda}(\Top(T))))& = \sigma_{\nu(T)}(S(T)) = t_{\nu(T)}(S(T)) = \zeta_1^{\nu(T)_1 - \nu(T)_2}\cdots \zeta_n^{\nu(T)_{n}}(S(T)) \\
 &= \Top(T).
 \end{align*}
 Let $\tau \in \PSYT(\lambda;T)$ and suppose for sake of induction that $\tau = \sigma_{w_{\tau}}(S(\tau))$. Now let $s_i(\tau) < \tau$. If $w_{\tau}(i) > w_{\tau}(i+1)$, then $S(s_i(\tau)) = S(\tau)$ and $\sigma_{w_{s_i(\tau)}}= s_i\sigma_{w_{\tau}}$ so that
 \[
 \sigma_{w_{s_i(\tau)}}(S(s_i(\tau))) = s_i\sigma_{w_{\tau}}(S(\tau)) = s_i(\tau).
 \]
 In the case instead that $w_{\tau}(i) = w_{\tau}(i+1)$ with $c_{\tau}(i+1)-c_{\tau}(i) > 1$, then $S(s_i(\tau)) = s_j(S(\tau))$ and~${\sigma_{w_{s_i(\tau)}} = \sigma_{w_{\tau}}}$ where $j = \sigma^{-1}(i)$ and $\sigma$ is the minimal length permutation with $\sigma(\sort(w_{\tau})) = w_{\tau}$.
 Then
 \[
 \sigma_{w_{s_i(\tau)}}(S(s_i(\tau)))=\sigma_{w_{\tau}}(s_jS(\tau))= (\sigma_{w_{\tau}}s_j)(S(\tau))= (s_i\sigma_{w_{\tau}})(S(\tau))= s_i(\tau).\tag*{\qed}
 \] \renewcommand{\qed}{}
\end{proof}

We may now obtain the desired bijection.

\begin{Proposition}\label{bijection between PSYT and pairs on min coset reps and SYT}
 The map \smash{$\Xi_{\lambda}\colon \PSYT(\lambda) \rightarrow \bigl(\widehat{\mathfrak{S}}_n/\mathfrak{S}_n \bigr)^{+} \times \SYT(\lambda)$} given by
 $\Xi_{\lambda}(\tau) := (\sigma_{w_{\tau}},\allowbreak S(\tau))$ is a~bijection.
\end{Proposition}
\begin{proof}
 It is immediate from Proposition~\ref{generating PSYT from SYT} that $\Xi_{\lambda}$ is injective. But, it is straightforward to check that given any $\sigma \in \bigl(\widehat{\mathfrak{S}}_n/\mathfrak{S}_n \bigr)^{+}$ and $S \in \SYT(\lambda)$, $\sigma(S)$ is a~well defined element of~$\PSYT(\lambda)$ in the sense of Definition~\ref{action of extended affine permutations of PSYT}. To see this, let $\mu$ be weakly decreasing. Using the same argument from the proof of Lemma~\ref{creation operator Lemma}, we see that $t_{\mu}(S) \in \PSYT(\lambda)$ is given by~${t_{\mu}(S)(\square) = S(\square) q^{\mu_{S(\square)}}}$. We proceed by induction on $\gamma \in \mathfrak{S}_n/\mathfrak{S}_{\text{Stab}(\mu)}$. The base case $\gamma = 1$ is trivial. Suppose that $\gamma(t_{\mu}(S))$ is a~well-defined element of $\PSYT(\lambda)$ for some $\gamma \in \mathfrak{S}_n/\mathfrak{S}_{\text{Stab}(\mu)}$ and $s_i\gamma > \gamma$ for some $1\leq i \leq n-1$, i.e.,
 $\mu_{\gamma^{-1}(i)} > \mu_{\gamma^{-1}(i+1)}$. It suffices to check that if~${\square_1,\square_2 \in \lambda}$ with $\square_1$ either above in the same column or to the left and in the same row as~$\square_2$ and $\mu_{S(\square_1)} = \mu_{S(\square_2)}$, then $s_i(\gamma(S(\square_1))) < s_i(\gamma(S(\square_2)))$. Since only the labels of the boxes labeled $i$ or $i+1$ change, we only need to consider the situation when either $\square_1$ or $\square_2$ are labeled by $i$ or $i+1$. Note that this excludes the situation where $\gamma(S(\square_1)) = i$ and $\gamma(S(\square_1)) = i+1$ since $\mu_{S(\square_1)} \neq \mu_{S(\square_2)}$ for these boxes by assumption. Now we have the following case work:
 \begin{itemize}
\itemsep=0pt
 \item If $ \gamma(S(\square_1)) = i$, then $ i+1 < \gamma(S(\square_2))$ so $s_i(\gamma(S(\square_1))) = i+1 < s_i(\gamma(S(\square_2)))$.
 \item If $\gamma(S(\square_1)) = i+1$, then $i+2 \leq \gamma(S(\square_2))$ so $s_i(\gamma(S(\square_1))) = i < i+2 \leq s_i(\gamma(S(\square_2)))$.
 \item If $\gamma(S(\square_2)) = i$, then $\gamma(S(\square_1))< i$ so $s_i(\gamma(S(\square_1)) = \gamma(S(\square_1)) < i < i+1 = s_i(\gamma(S(\square_2)))$.
 \item If $\gamma(S(\square_2)) = i+1$,
 then $\gamma(S(\square_1)) < i$ so $s_i(\gamma(S(\square_1)) = \gamma(S(\square_1)) < i = s_i(\gamma(S(\square_2)))$.
 \end{itemize}
 Therefore, \smash{$s_i\gamma t_{\mu}(S)$} is a~well-defined element of \smash{$\PSYT(\lambda)$}. From Lemma~\ref{min coset reps Lemma}, we know that every element \smash{$\sigma \in \bigl(\widehat{\mathfrak{S}}_n/\mathfrak{S}_n \bigr)^{+}$} can be written in the form \smash{$\sigma = \gamma t_{\mu}$} for some $\mu$ weakly decreasing and \smash{$\gamma \in \mathfrak{S}_n/\mathfrak{S}_{\text{Stab}(\mu)}$}. Thus~\smash{$\sigma(S) \in \PSYT(\lambda)$} is well defined for all~\smash{$\sigma \in \bigl(\widehat{\mathfrak{S}}_n/\mathfrak{S}_n \bigr)^{+}$} and~${S \in \SYT(\lambda)}$.
 This shows that $\Xi_{\lambda}$ is also surjective and thus bijective.
\end{proof}

\subsection{Finite Hecke algebra}
For references to the definitions of finite, affine, and double affine Hecke algebras, we refer the reader to Cherednik's book~\cite{CherednikBook}. Here we detail our conventions for the finite Hecke algebras in type $A_{n-1}$ occurring in this paper.

\begin{Definition}\label{finite hecke alg Definition}
 Define the finite Hecke algebra $\sH_n$ to be the $\mathbb{Q}(q,t)$-algebra generated by~${T_1,\dots, T_{n-1}}$ subject to the relations
 \begin{itemize}\itemsep=0pt
 \item $(T_i-1)(T_i+t) = 0$ for $1\leq i \leq n-1$,
 \item $T_iT_{i+1}T_1 = T_{i+1}T_iT_{i+1}$ for $1\leq i \leq n-2$,
 \item $T_iT_j = T_jT_i$ for $|i-j| > 1$.
 \end{itemize}
 Define the \textit{Jucys--Murphy} elements $\overline{\theta}_1,\dots, \overline{\theta}_n \in \sH_n$ by $\overline{\theta}_1 := 1$ and $\overline{\theta}_{i+1}:= tT_{i}^{-1}\overline{\theta}_{i}T_{i}^{-1}$ for $1\leq i \leq n-1$. Further, define $\overline{\varphi}_1,\dots, \overline{\varphi}_{n-1}$ by $\overline{\varphi}_i: = \bigl(tT_{i}^{-1}\bigr)\overline{\theta}_i - \overline{\theta}_{i}\bigl(tT_i^{-1}\bigr)$. For a~permutation~${\sigma \in \mathfrak{S}_n}$ and a~reduced expression $\sigma = s_{i_1}\cdots s_{i_r}$, we write $T_{\sigma} := T_{i_1}\cdots T_{i_r}$.
\end{Definition}

Note that the definition of $\sH_n$ in \cite[Section~2.1]{DL_2011} differs from our definition given above. To translate between our differing conventions we may identify as follows:
\begin{itemize}\itemsep=0pt
 \item $s \rightarrow t$ (parameters),
 \item $T_i \rightarrow tT_i^{-1}$ (Hecke elements).
\end{itemize}

\begin{Remark}
 There are natural algebra inclusions $\sH_{n} \rightarrow \sH_{n+1}$ given by $T_i \rightarrow T_i$ for $1\leq i \leq n-1$. Under this embedding $\overline{\theta}_i \rightarrow \overline{\theta}_i$ for $1 \leq i \leq n$, so we can naturally identify the copies of $\overline{\theta}_i$ in both $\sH_n$ and $\sH_{n+1}$.
\end{Remark}

We require the following list of relations. The formulas involving the elements $\overline{\varphi}_i$ we refer to as the intertwiner relations.

\begin{Proposition}\label{additional relations for finite hecke}
 The following relations hold:
 \begin{itemize}
\itemsep=0pt
 \item $\overline{\theta}_i = t^{i-1}T_{i-1}^{-1}\cdots T_{1}^{-1}T_{1}^{-1}\cdots T_{i-1}^{-1}$ for $1\leq i \leq n$,
 \item $\overline{\theta}_i\overline{\theta}_j = \overline{\theta}_j\overline{\theta}_i$ for $1\leq i,j \leq n$,
 \item $T_i\overline{\theta}_j = \overline{\theta_j}T_i$ for $j \notin \{i,i+1\}$,
 \item $\overline{\varphi}_i = tT_{i}^{-1}\bigl(\overline{\theta}_i - \overline{\theta}_{i+1}\bigr) + (t-1)\overline{\theta}_{i+1}$ for $1 \leq i \leq n-1$,
 \item $\overline{\varphi}_i\overline{\varphi}_{i+1}\overline{\varphi}_i=\overline{\varphi}_{i+1}\overline{\varphi}_i\overline{\varphi}_{i+1}$ for $1\leq i \leq n-1$,
 \item $\overline{\varphi}_i\overline{\varphi}_j = \overline{\varphi}_j\overline{\varphi}_i$ for $|i-j| >1$,
 \item $\overline{\varphi}_{i}\overline{\theta}_j = \overline{\theta}_{s_i(j)}\overline{\varphi}_{i}$ for $1\leq i \leq n-1$ and $1 \leq j \leq n$,
 \item $\overline{\varphi}_i^{2} = \bigl(t\overline{\theta}_i-\overline{\theta}_{i+1}\bigr)\bigl(t\overline{\theta}_{i+1}-\overline{\theta}_{i}\bigr)$.
 \end{itemize}
\end{Proposition}
\begin{proof}
This result follows directly from using the map $\rho_n$ defined in Definition~\ref{map from affine to finite hecke} in conjunction with Proposition~\ref{additional relations for affine hecke} which is independently proven later.
\end{proof}

The following definition gives a~description of the irreducible representations of $\sH_n$. There are many equivalent methods for defining these representations but we choose to specify eigenvectors for the Jucys--Murphy elements $\overline{\theta}_i$ directly as we require these eigenvectors throughout this paper.

\begin{Definition}\label{irreps for finite hecke Definition}
 Let $\lambda \in \Par$ with $|\lambda| = n$. Define $S_{\lambda}$ to be the $\sH_n$-module with basis $e_{\tau}$ for~${\tau \in \SYT(\lambda)}$ with action determined by the following relations:
 \begin{itemize}
\itemsep=0pt
 \item $\overline{\theta}_i(e_{\tau}) = t^{c_{\tau}(i)}e_{\tau}$.
 \item If $s_i(\tau) > \tau$, then $\overline{\varphi}_i(e_{\tau}) = \bigl(t^{c_{\tau}(i)}-t^{c_{\tau}(i+1)}\bigr)e_{s_i(\tau)}$.
 \item If the labels $i,i+1$ are in the same row in $\tau$, then $T_i(e_{\tau}) = e_{\tau}$.
 \item If the labels $i,i+1$ are in the same column in $\tau$, then $T_i(e_{\tau}) = -te_{\tau}$.
 \end{itemize}
\end{Definition}

Using the relations from Proposition~\ref{additional relations for finite hecke}, we can show the following more explicit form for the action of the $T_i$ on the $\SYT(\lambda)$ basis:
\begin{itemize}
\itemsep=0pt
 \item If $s_i(\tau)> \tau$, then
 \[
 T_i(e_{\tau}) = e_{s_i(\tau)} + \frac{(1-t)t^{c_{\tau}(i)}}{t^{c_{\tau}(i)} - t^{c_{\tau}(i+1)}} e_{\tau}.
 \]

 \item If $s_i(\tau) < \tau$, then
 \[
 T_i(e_{\tau}) = -\frac{\bigl(t^{c_{\tau}(i+1)+1}-t^{c_{\tau}(i)}\bigr)\bigl(t^{c_{\tau}(i)+1}-t^{c_{\tau}(i+1)}\bigr)}{\bigl(t^{c_{\tau}(i+1)}-t^{c_{\tau}(i)}\bigr)^2} e_{s_i(\tau)} + \frac{(1-t)t^{c_{\tau}(i)}}{t^{c_{\tau}(i)} - t^{c_{\tau}(i+1)}} e_{\tau}.
 \]

\end{itemize}

Implicit in the above definition is the following result.

\begin{Proposition}
 Definition~{\rm\ref{irreps for finite hecke Definition}} is well posed, i.e., the action of the operators $T_i$ on $S_{\lambda}$ define an irreducible $\sH_n$-module.
\end{Proposition}
\begin{proof}
 As this construction is standard, we only give an outline. It follows from standard theory~\cite[Theorem 3.9]{Murphy} for the finite Hecke algebra $\sH_n$ (analogous to that of the symmetric group~$\mathfrak{S}_n$ in characteristic $0$) that there exists an irreducible representation of $\sH_n$, $S_{\lambda}$, corresponding to the partition $\lambda$ with a~basis of weight vectors for the Jucys--Murphy elements~$\overline{\theta}_i$, $v_{\tau}$~say, indexed by $\tau \in \SYT(\lambda)$. Further, the weights are given by $\overline{\theta}_i(v_{\tau}) = t^{c_{\tau}(i)} v_{\tau}$. As these weights are all distinct, it follows that this basis is unique up to re-normalization by nonzero scalars.
 The presentation given in Definition~\ref{irreps for finite hecke Definition} fixes a~specific normalization given by choosing first~${e_{\tau^{\rm rs}_{\lambda}} = v_{\tau^{\rm rs}_{\lambda}}}$, and then building the full basis $e_{\tau}$ using the intertwiner $\overline{{\varphi_i}}$ relations in Proposition~\ref{additional relations for finite hecke} with the choice that whenever $s_i(\tau) > \tau$ we have that \smash{$\overline{{\varphi_i}}(e_{\tau}) = \bigl(t^{c_{\tau}(i)}-t^{c_{\tau}(i+1)}\bigr)e_{s_i(\tau)}$}. Up to an initial arbitrary choice for the scalar multiple of $e_{\tau^{\rm rs}_{\lambda}}$, this uniquely determines the rest of the coefficients of the $e_{\tau}$.
\end{proof}

\begin{Remark}
 The set $\{ \lambda \in \Par \mid |\lambda| = n\}$ gives a~full set of irreducible $\sH_n$-modules up to isomorphism. Note that for $\tau,\tau' \in \SYT(\lambda)$, the $\overline{\theta}$-weights of $e_{\tau}= e_{\tau'}$ are equal if and only if~${\tau = \tau'}$.
\end{Remark}

The modules $S_{\lambda}$ satisfy the following property which will serve as the basis for the main algebraic constructions in this paper.

\begin{Lemma}\label{finite hecke map}
 Let $\lambda \in \Par$ and $n \geq n_{\lambda}$. Let $\square_0$ be the unique square in the skew-diagram \smash{$\lambda^{(n+1)} / \lambda^{(n)}$}. Consider the map \smash{$\mathfrak{q}_{\lambda}^{(n)}\colon\lambda^{(n+1)}\rightarrow \lambda^{(n)}$} given for \smash{$\tau \in \SYT\bigl(\lambda^{(n+1)}\bigr)$} as
 \[
 \mathfrak{q}_{\lambda}^{(n)}(e_{\tau}) := \begin{cases}
 e_{\tau|_{\lambda^{(n)}}}, & \tau(\square_0) = n+1,\\
 0, & \tau(\square_0) \neq n+1.
 \end{cases}
 \]

 Then $\mathfrak{q}_{\lambda}^{(n)}$ is a~$\sH_{n}$-module map.
\end{Lemma}
\begin{proof}
Let $\tau \in \SYT\bigl(\lambda^{(n+1)}\bigr)$. If $\tau(\square_0) = n+1$, then by Definition~\ref{irreps for finite hecke Definition}, the action of each $T_i$ for $1\leq i \leq n-1$ is the same on both $\tau \in S_{\lambda^{(n+1)}}$ and $\tau\mid _{\lambda^{(n)}} \in S_{\lambda^{(n)}}$. If $\tau(\square_0) \neq n+1$, then for all $1\leq i \leq n-1$, $T_i(e_{\tau})$ is a~linear combination of $e_{\tau}$ and $e_{s_i(\tau)}$ (when defined). Note that, by construction, \smash{$\mathfrak{q}_{\lambda}^{(n)}(e_{\tau}) = 0$}. Necessarily, when defined $s_i(\tau)(\square_0) \neq n+1$, and so \smash{$\mathfrak{q}_{\lambda}^{(n)}(e_{s_i(\tau)}) = 0$}. Thus \smash{$\mathfrak{q}_{\lambda}^{(n)}(T_i(e_{\tau})) = 0 = T_i\mathfrak{q}_{\lambda}^{(n)}(e_{\tau})$}. In all cases, we have that \smash{$\mathfrak{q}_{\lambda}^{(n)}(T_i(e_{\tau})) = T_i\mathfrak{q}_{\lambda}^{(n)}(e_{\tau})$}. Hence, \smash{$\mathfrak{q}_{\lambda}^{(n)}$} is a~$\sH_n$-module map.
\end{proof}

\subsection{Affine Hecke algebra}
We are interested in the presentation of the AHAs of type ${\rm GL}_n$ which follows.
\begin{Definition}\label{first presentation of AHA}
 Define the affine Hecke algebra $\sA_n$ to be the $\mathbb{Q}(q,t)$-algebra generated by~${T_1,\dots, T_{n-1}}$ and $\theta_1^{\pm 1},\dots, \theta_n^{\pm 1}$ subject to the relations
 \begin{itemize}
\itemsep=0pt
 \item $\sH_n \subset \sA_n$,
 \item $\theta_i\theta_j = \theta_j\theta_i$ for all $1\leq i,j \leq n$,
 \item $\theta_{i+1} = tT_i^{-1}\theta_i T_{i}^{-1}$ for $1\leq i \leq n-1$,
 \item $T_i\theta_j = \theta_j T_i$ for $j \notin \{i,i+1\}$.
 \end{itemize}

 Define the special elements $\pi_n$ and $\varphi_1,\dots, \varphi_{n-1}$ of $\sA_n$ by
 \begin{itemize}
\itemsep=0pt
 \item $\pi_n:= t^{n-1}\theta_1T_{1}^{-1}\cdots T_{n-1}^{-1}$,
 \item $\varphi_i := \bigl(tT_{i}^{-1}\bigr)\theta_i - \theta_i \bigl(tT_{i}^{-1}\bigr)$.
 \end{itemize}

 We denote by $\Theta^{(n)}$ the commutative subalgebra of $\sA_n$ generated by $\theta_1^{(n)},\dots, \theta_n^{(n)}$.
\end{Definition}

It is important to note that when converting between the AHA conventions in this paper and those in Dunkl--Luque \cite[Section~2.1]{DL_2011} the standard Cherednik elements $\xi_i$ of Dunkl--Luque do \textit{not} agree with the $\theta_i$ above. In particular, after the appropriate translation into our conventions we have that $\xi_i$ are given by $\xi_i = t^{n-i+1}T_{i-1}\cdots T_1 \pi_n T_{n-1}^{-1}\cdots T_i^{-1}$ as opposed to~$\smash{\theta_i = t^{-(n-i)}T_{i-1}^{-1}\cdots T_1^{-1}\pi_n T_{n-1}\cdots T_{i}}$. The distinction between the standard Cherednik elements $\xi_i$ and the reversed orientation Cherednik elements $\theta_i$ is important in this paper since the latter yield operators with additional stability properties which the $\xi_i$ do not satisfy.

\begin{Remark}
 We use the notation \smash{$\theta_i^{(n)}$} and \smash{$\theta_i^{(m)}$} to differentiate between the copies of~$\theta_i$ in~$\sA_n$ and $\sA_m$ for $n \neq m$.
\end{Remark}

The following proposition is required at many points throughout this paper when we study $\Theta^{(n)}$-weight vectors. We include the proofs of these relations for completeness and to emphasize that we may use intertwiner theory for AHAs with the $\theta_i$ elements instead of the standard $\xi_i$ with only slight differences.

\begin{Proposition}\label{additional relations for affine hecke}
 The following relations hold:
 \begin{itemize}
\itemsep=0pt
 \item $\varphi_i = tT_{i}^{-1}(\theta_i - \theta_{i+1}) + (t-1)\theta_{i+1} = (\theta_{i+1}-\theta_i)tT_i^{-1} + (1-t)\theta_{i+1}$ for $1 \leq i \leq n-1$,
 \item $\varphi_{i}\theta_j = \theta_{s_i(j)}\varphi_{i}$ for $1\leq i \leq n-1$ and $1 \leq j \leq n$,
 \item $\varphi_i^{2} = (t\theta_i-\theta_{i+1})(t\theta_{i+1}-\theta_{i})$,
 \item $\varphi_i\varphi_{i+1}\varphi_i=\varphi_{i+1}\varphi_i\varphi_{i+1}$ for $1\leq i \leq n-2$,
 \item $\varphi_i\varphi_j = \varphi_j\varphi_i$ for $|i-j| >1$.
 \end{itemize}
\end{Proposition}

\begin{proof}
We proceed by proving each of these relations in the order in which they appear above.

 Let $1 \leq i \leq n-1$. Then
\begin{align*}
 \varphi_i &= tT_i^{-1}\theta_i - \theta_i\bigl(tT_i^{-1}\bigr) = tT_i^{-1}\theta_i - T_i\theta_{i+1} = tT_i^{-1}\theta_i - \bigl(tT_i^{-1}+1-t\bigr)\theta_{i+1}\\
 &= tT_i^{-1}(\theta_i - \theta_{i+1}) + (t-1)\theta_{i+1}.
\end{align*}

By a~similar calculation, we also get
$ \varphi_i = (\theta_{i+1}-\theta_i)tT_i^{-1} + (1-t)\theta_{i+1}$.
This can also be written as
$ \varphi_i = (\theta_{i+1}-\theta_i)T_i + (1-t)\theta_i
$
which we need later in this proof.

Now we see
\begin{align*}
 \varphi_i\theta_i &= tT_i^{-1}(\theta_i-\theta_{i+1})\theta_i + (t-1)\theta_{i+1}\theta_i = tT_i^{-1}\theta_i(\theta_i-\theta_{i+1}) + (t-1)\theta_{i+1}\theta_i \\
 & = \theta_{i+1}T_i(\theta_i -\theta_{i+1}) + (t-1)\theta_{i+1}\theta_i
 = \theta_{i+1}( T_i(\theta_i-\theta_{i+1}) + (t-1)\theta_i ) \\
 & = \theta_{i+1}\bigl( \bigl(tT_i^{-1} + 1-t\bigr)(\theta_i -\theta_{i+1}) + (t-1)\theta_i \bigr) = \theta_{i+1} \bigl( tT_i^{-1}(\theta_i -\theta_{i+1}) + (t-1)\theta_{i+1}\bigr) \\
&= \theta_{i+1}\varphi_i
\end{align*}
and
\begin{align*}
 \varphi_i\theta_{i+1} &= tT_i^{-1}(\theta_i - \theta_{i+1})\theta_{i+1}+(t-1)\theta_{i+1}^{2} = (T_i+t-1)\theta_{i+1}(\theta_i-\theta_{i+1}) + (t-1)\theta_{i+1}^2 \\
 & = T_i\theta_{i+1}(\theta_i-\theta_{i+1}) + (t-1)\left( \theta_{i+1}(\theta_i-\theta_{i+1})+\theta_{i+1}^2 \right) \\
 &= t\theta_iT_i^{-1}(\theta_i - \theta_{i+1}) +(t-1)\theta_i\theta_{i+1} = \theta_i\left( tT_i^{-1}(\theta_i-\theta_{i+1}) +(t-1)\theta_{i+1} \right) = \theta_i\varphi_i.
\end{align*}

For any $j \notin \{i,i+1\}$, it follows since $\theta_j$ commutes with both $\theta_i$ and $T_i$, that
$
 \varphi_i\theta_j = \theta_j \varphi_i$.
Thus for any $1\leq j \leq n$,
$
 \varphi_i\theta_j = \theta_{s_i(j)}\varphi_i$.

Now we have that
\begin{align*}
 \varphi_i^2 &= \bigl(tT_i^{-1}\theta_i - \theta_i tT_i^{-1}\bigr)^2 = t^2T_i^{-1}\theta_iT_i^{-1}\theta_i - t^2 T_i^{-1}\theta_i^2 T_i^{-1} - t^2\theta_iT_i^{-2}\theta_i + t^2\theta_iT_i^{-1}\theta_iT_i^{-1} \\
& = t\theta_{i+1}\theta_{i} - t\theta_{i+1}T_i \theta_i T_i^{-1} -t\theta_i \bigl(1+(t-1)T_i^{-1}\bigr) \theta_i + t\theta_i\theta_{i+1} \\
& = 2t\theta_i\theta_{i+1} - t\theta_{i+1}\bigl(tT_i^{-1}+1-t\bigr)\theta_iT_i^{-1} - t\theta_i^2 +t(1-t)\theta_iT_i^{-1}\theta_i \\
& = 2t\theta_i\theta_{i+1} -t^2 \theta_{i+1}T_i^{-1}\theta_iT_i^{-1} +t(t-1)\theta_i\theta_{i+1}T_i^{-1}-t\theta_i^2 + (1-t)\theta_i\theta_{i+1} T_i \\
& = 2t\theta_i\theta_{i+1} - t\theta_{i+1}^2 -t\theta_{i}^2 + (1-t)\theta_i \theta_{i+1}\bigl(T_i - tT_i^{-1}\bigr)\\
& = 2t\theta_i\theta_{i+1}-t\theta_{i+1}^{2} - t\theta_i^2 +(1-t)^2 \theta_i \theta_{i+1} = \bigl(1+t^2\bigr)\theta_i\theta_{i+1} - t\theta_{i+1}^{2}-t\theta_i^{2}\\
& = (t\theta_i - \theta_{i+1})(t\theta_{i+1}-\theta_i).
\end{align*}

Now suppose $1\leq i \leq n-2$. By expanding each of the $\varphi_j$ from right to left using $\varphi_j = (\theta_{j+1}-\theta_j)T_j + (1-t)\theta_j$ and repeatedly applying the relation $\varphi_{j}\theta_k = \theta_{s_j(k)}\varphi_{j}$, we find
\begin{align*}
 \varphi_i\varphi_{i+1}\varphi_i ={}& (\theta_{i+2}-\theta_{i+1})(\theta_{i+2}-\theta_i)(\theta_{i+1}-\theta_i)T_iT_{i+1}T_i \\
& + (1-t)\theta_i (\theta_{i+2}-\theta_{i+1})(\theta_{i+2}-\theta_{i})T_{i+1}T_i \\
& + (1-t)\theta_{i+1}(\theta_{i+2}-\theta_i)(\theta_{i+1}-\theta_i)T_iT_{i+1} + (1-t)^{2}\theta_i \theta_{i+1} (\theta_{i+2}-\theta_i) T_i \\
& + (1-t)^2\theta_i \theta_{i+1} (\theta_{i+2}-\theta_i) T_{i+1} \\
 & + \bigl( t(1-t)\theta_i(\theta_{i+2}-\theta_{i+1})(\theta_{i+1}-\theta_i) + (1-t)^3 \theta_i^2 \theta_{i+1} \bigr).
\end{align*}

Using the same method, we also see that
\begin{align*}
 \varphi_{i+1}\varphi_{i}\varphi_{i+1}= {}&(\theta_{i+1}-\theta_i)(\theta_{i+2}-\theta_i)(\theta_{i+2}-\theta_{i+1})T_{i+1}T_iT_{i+1} \\
&+ (1-t)\theta_{i+1}(\theta_{i+1}-\theta_i)(\theta_{i+2}-\theta_{i})T_iT_{i+1} \\
& + (1-t)\theta_i (\theta_{i+2}-\theta_{i})(\theta_{i+2}-\theta_{i+1})T_{i+1}T_i + (1-t)^2 \theta_{i}\theta_{i+1}(\theta_{i+2}-\theta_i)T_{i+1} \\
& + (1-t)^2 \theta_{i}\theta_{i+1}(\theta_{i+2}-\theta_i)T_{i} \\
& + \bigl( t(1-t)\theta_i(\theta_{i+1}-\theta_i)(\theta_{i+2}-\theta_{i+1}) + (1-t)^3 \theta_{i}^2 \theta_{i+1} \bigr).
\end{align*}
From here, we may use the braid relation \smash{$T_iT_{i+1}T_i = T_{i+1}T_iT_{i+1}$} and some rearrangement of terms to see $\varphi_i\varphi_{i+1}\varphi_i = \varphi_{i+1}\varphi_{i}\varphi_{i+1}$.

Lastly, consider $|i - j| >1$. Since $T_iT_j = T_jT_i$, $T_i\theta_j = \theta_j T_i$, and $\theta_i\theta_j= \theta_j \theta_i$, we readily find that $\varphi_i\varphi_j = \varphi_j\varphi_i$.
\end{proof}

In this paper, we are interested in AHA modules which are \textit{pulled back} from irreducible finite Hecke representations. To do this we need to define algebra surjections $\sA_n \rightarrow \sH_n$. There are many such choices for these maps but we choose the maps $\rho_n$ defined below carefully so that the AHA modules we consider in this paper satisfy nontrivial stability conditions.

\begin{Definition}\label{map from affine to finite hecke}
 Define the $\mathbb{Q}(q,t)$-algebra homomorphism $\rho_n\colon\sA_n \rightarrow \sH_n$ by
 \begin{itemize}
\itemsep=0pt
 \item $\rho_n(T_i) = T_i$ for $1\leq i \leq n-1$,
 \item $\rho_n(\theta_1) = 1$.
 \end{itemize}
 For a~$\sH_n$-module $V$, we denote by $\rho_n^{*}(V)$ the $\sA_n$-module with action defined for $v \in V$ and~${X \in \sA_n}$ by $X(v) := \rho_n(X)(v)$.
\end{Definition}

\begin{Remark}
 Note that $\rho_n(\pi_n) = t^{n-1}T_1^{-1}\cdots T_{n-1}^{-1}$ and for all $1\leq i \leq n$, $\rho_n(\theta_i) = \overline{\theta}_i$. Thus~$\rho_n$ maps the reversed orientation Cherednik elements $\theta_i \in \sA_n$ to the Jucys--Murphy elements~${\overline{\theta}_i \in \sH_n}$. Further, for $\lambda \in \Par$ with $|\lambda| = n$, $\rho_n^{*}(\lambda)$ is an irreducible $\sA_n$-module with a~basis of $\theta$-weight vectors $\{e_{\tau}\}_{\tau \in \SYT(\lambda)}$ with distinct weights.
\end{Remark}

\subsection{Positive double affine Hecke algebra}
We use the following presentation for the positive DAHAs in type ${\rm GL}_n$.
\begin{Definition}\label{daha def}
 Define the positive double affine Hecke algebra $\sD_n$ to be the $\mathbb{Q}(q,t)$-algebra generated by $T_1,\dots,T_{n-1}$, $\theta_1^{\pm 1},\dots, \theta_n^{\pm 1}$, and $X_1,\dots, X_n$ subject to the relations
 \begin{itemize}
\itemsep=0pt
 \item $T_1,\dots,T_{n-1}$ and $\theta_1^{\pm 1},\dots, \theta_n^{\pm 1}$ satisfy the relations of $\sA_n$ in Definition~\ref{first presentation of AHA},
 \item $X_iX_j = X_jX_i$ for $1\leq i,j \leq n$,
 \item $X_{i+1} = tT_i^{-1}X_iT_i^{-1}$ for $1\leq i \leq n-1$,
 \item $T_iX_j = X_jT_i$ for $1\leq i\leq n-1$ and $1\leq j \leq n$ with $j \notin \{i,i+1\}$,
 \item $\pi_nX_i\pi_n^{-1} = X_{i+1}$ for $1\leq i \leq n-1$,
 \item $\pi_nX_n\pi_n^{-1} = qX_1$.
 \end{itemize}
 Here $\pi_n$ is the same as in Definition~\ref{first presentation of AHA}. Define the special element $\gamma_n:= X_nT_{n-1}\cdots T_1$.
\end{Definition}

The element $\gamma_n$ satisfies some nice intertwining relations which we need for later when we study $\Theta^{(n)}$-weight vectors.

\begin{Lemma}\label{intertwiner Lemma}
 The following hold:
 \begin{itemize}
\itemsep=0pt
 \item $\theta_i\gamma_n = \gamma_n \theta_{i+1}$ for $1\leq i \leq n-1$,
 \item $\theta_n \gamma_n = \gamma_n q\theta_1$.
 \end{itemize}
\end{Lemma}

\begin{proof}
 Let $1 \leq i \leq n-1$. We find that
 \begin{align*}
 \theta_i\gamma_n& = t^{-(n-i)}T_{i-1}^{-1}\cdots T_1^{-1} \pi_n T_{n-1}\cdots T_i X_n T_{n-1}\cdots T_1 \\
 & = t^{-(n-i)}T_{i-1}^{-1}\cdots T_1^{-1} \pi_n T_{n-1}X_n T_{n-2}\cdots T_i T_{n-1} \cdots T_1 \\
 & = t^{-(n-i)}T_{i-1}^{-1}\cdots T_1^{-1}\pi_n tX_{n-1}T_{n-1}^{-1} T_{n-2}\cdots T_i T_{n-1} \cdots T_1 \\
 & = t^{-(n-(i+1))} T_{i-1}^{-1}\cdots T_1^{-1} X_n \pi_n T_{n-1}^{-1} T_{n-2}\cdots T_i T_{n-1}\cdots T_1 \\
 & = t^{-(n-(i+1))} X_n T_{i-1}^{-1}\cdots T_1^{-1} \pi_n T_{n-1}^{-1} T_{n-2}\cdots T_i (T_{n-1}\cdots T_1).
 \end{align*}
 From the braid relations, we see that for all $1\leq j \leq n-2$, $T_j(T_{n-1}\cdots T_1) = (T_{n-1}\cdots T_1)T_{j+1}$
 and hence
 \begin{gather*}
 t^{-(n-(i+1))} X_n T_{i-1}^{-1}\cdots T_1^{-1} \pi_n T_{n-1}^{-1} T_{n-2}\cdots T_i (T_{n-1}\cdots T_1) \\
 \qquad = t^{-(n-(i+1))} X_n T_{i-1}^{-1}\cdots T_1^{-1} \pi_n T_{n-1}^{-1}(T_{n-1}\cdots T_1)T_{n-1}\cdots T_{i+1} \\
 \qquad= t^{-(n-(i+1))} X_n T_{i-1}^{-1}\cdots T_1^{-1} \pi_n T_{n-2} \cdots T_1 T_{n-1}\cdots T_{i+1} \\
 \qquad= t^{-(n-(i+1))} X_n T_{i-1}^{-1}\cdots T_1^{-1} T_{n-1} \cdots T_2 \pi_n T_{n-1}\cdots T_{i+1}\\
 \qquad= t^{-(n-(i+1))} X_n T_{i-1}^{-1}\cdots T_1^{-1} T_{n-1} \cdots T_2 T_1T_1^{-1}\pi_n T_{n-1}\cdots T_{i+1}\\
 \qquad= t^{-(n-(i+1))} X_nT_{n-1} \cdots T_1T_{i}^{-1}\cdots T_2^{-1}T_1^{-1}\pi_n T_{n-1}\cdots T_{i+1}\\
 \qquad= (X_nT_{n-1} \cdots T_1) \bigl(t^{-(n-(i+1))}T_{i}^{-1}\cdots T_1^{-1}\pi_n T_{n-1}\cdots T_{i+1}\bigr)
 = \gamma_n \theta_{i+1}.
 \end{gather*}

 Now we consider the last case
 \begin{align*}
 \theta_n \gamma_n& = T_{n-1}^{-1}\cdots T_1^{-1} \pi_n T_{n-1}\cdots T_1 = T_{n-1}^{-1}\cdots T_1^{-1} qX_1 \pi_n T_{n-1}\cdots T_1 \\
 & = t^{-(n-1)}X_n T_{n-1}\cdots T_1 q\pi_n T_{n-1}\cdots T_1= (X_nT_{n-1} \cdots T_1) \bigl(qt^{-(n-1)}\pi_n T_{n-1}\cdots T_1\bigr) \\
 & = \gamma_n q\theta_1.\tag*{\qed}
 \end{align*} \renewcommand{\qed}{}
\end{proof}

Recall the definition of the intertwiner elements $\varphi_i$ in Definition~\ref{first presentation of AHA}. We use the elements~${\{\varphi_1,\dots, \varphi_{n-1}, \gamma_n\}}$ to define intertwiner operators corresponding to elements of $\widehat{\mathfrak{S}}_n^{+}$.

\begin{Definition}\label{intertwiners def}
 For any $\sigma \in \widehat{\mathfrak{S}}_n^{+}$ with $
 \sigma = w_1\dots w_r$ written minimally in terms of the generators $w_i \in \{s_1,\dots, s_{n-1},\widetilde{\gamma}_n \}$, define $\varphi_{\sigma} \in \sD_n$ by
$ \varphi_{\sigma}:= \varpi_1\cdots \varpi_r \in \sD_n
$, where $\varpi_i:= \varphi_j$ if~${w_i = s_j}$ and $\varpi_{i} := \gamma_n$ if $w_i = \widetilde{\gamma}_n$.

In particular, we have that $\varphi_{s_i} = \varphi_{i}$ and $\varphi_{\widetilde{\gamma}_n} = \gamma_n$.
\end{Definition}

The main utility of considering the intertwiner operators $\varphi_{\sigma}$ comes from the next lemma.

\begin{Lemma}\label{weight vector intertwiners Lemma}
 If $v$ is a~$\Theta^{(n)}$-weight vector in some $\sD_n$-module with weight $\alpha= (\alpha_1,\dots,\alpha_n)$ and $\sigma \in \mathfrak{S}_n$ with $\varphi_{\sigma}(v) \neq 0$, then $\varphi_{\sigma}(v)$ is a~$\Theta^{(n)}$-weight vector with weight $\alpha^{\sigma}$ given by the following recursion:
 \begin{itemize}
\itemsep=0pt
 \item $\alpha^{s_i} = (\alpha_1,\dots,\alpha_{i+1},\alpha_i,\dots\alpha_n)$,
 \item $\alpha^{\widetilde{\gamma}_n} = (\alpha_2,\dots, \alpha_{n}, q\alpha_1)$,
 \item $(\alpha^{\sigma_2})^{\sigma_1} = \alpha^{\sigma_1\sigma_2}$.
 \end{itemize}
\end{Lemma}
\begin{proof}
 This result follows using induction on $\widehat{\mathfrak{S}}_n^{+}$ and the relations in Proposition~\ref{additional relations for affine hecke} and Lemma~\ref{intertwiner Lemma}. For reference, see \cite[Section~3.3.3]{CherednikBook}.
\end{proof}

We are primarily interested in modules for the \textit{spherical} subalgebra of $\sD_n$.

\begin{Definition}\label{spherical DAHA Definition}
 Let $\epsilon^{(n)} \in \sH_n$ denote the trivial idempotent given by
 \[
 \epsilon^{(n)}:= \frac{1}{[n]_{t}!}\sum_{\sigma \in \mathfrak{S}_n} t^{{n\choose 2} - \ell(\sigma)} T_{\sigma},
 \]
 where
 \smash{$[n]_{t}!:= \prod_{i=1}^{n}\bigl(\frac{1-t^i}{1-t}\bigr)$}. The positive spherical double affine Hecke algebra $\sD_n^{\mathrm{sph}}$ is the (non-unital) subalgebra of $\sD_n$ given by \smash{$\sD_n^{\mathrm{sph}}:= \epsilon^{(n)}\sD_n \epsilon^{(n)}$}.
\end{Definition}

The element
\[\epsilon^{(n)}:= \frac{1}{[n]_{t}!}\sum_{\sigma \in \mathfrak{S}_n} t^{{n\choose 2} - \ell(\sigma)} T_{\sigma} \in \sH_n\] is uniquely determined by the following properties:
\begin{itemize}
\itemsep=0pt
 \item $\epsilon^{(n)} \neq 0$ (non-zero),
 \item $\bigl(\epsilon^{(n)}\bigr)^2 = \epsilon^{(n)}$ (idempotent),
 \item $\epsilon^{(n)}T_i = T_i\epsilon^{(n)}$ for all $1\leq i\leq n-1$ (central),
 \item $T_i\epsilon^{(n)} = \epsilon^{(n)}$ (trivial-like).
\end{itemize}

We use without proof that $\epsilon^{(n)}$ as defined in Definition~\ref{spherical DAHA Definition} satisfies these properties but it is straightforward to check this using the defining relations of \smash{$\sH_n$}. Since \smash{$\bigl(\epsilon^{(n)}\bigr)^2 = \epsilon^{(n)}$}, we see that \smash{$\sD_n^{\mathrm{sph}}$} is a~unital algebra with unit \smash{$\epsilon^{(n)}$}. The algebra \smash{$\sD_n^{\mathrm{sph}}$} contains both of the subalgebras~$\smash{\mathbb{Q}(q,t)[X_1,\dots,X_n]^{\mathfrak{S}_n}\epsilon^{(n)}}$ and \smash{$\mathbb{Q}(q,t)\big[\theta_1^{\pm 1},\dots, \theta_n^{\pm 1}\big]^{\mathfrak{S}_n}\epsilon^{(n)}$}.

We may use $\epsilon^{(n)}$ to generate modules for the spherical DAHA. Given any $\sD_n$-module V, the space $\epsilon^{(n)}(V)$ is naturally a~$\sD_n^{\mathrm{sph}}$-module. In the standard picture of Cherednik theory \cite[Section 3.8]{CherednikBook}, the polynomial representation of $\sD_n$ on $\mathbb{Q}(q,t)[x_1,\dots, x_n]$ is symmetrized using~$\epsilon^{(n)}$ to yield the symmetric polynomial representation of $\sD_n^{\mathrm{sph}}$ on $\mathbb{Q}(q,t)[x_1,\dots, x_n]^{\mathfrak{S}_n}$.

\begin{Remark}
 It is a~standard result \cite[Theorem 3.2.1]{CherednikBook} that $\sD_n$ is a~\textit{free} right $\sA_n$-module with basis \smash{$\{X^{\alpha}\}_{\alpha \in \mathbb{Z}^n_{\geq 0}}$}. Importantly, for our purposes, this implies that for any $\sA_n$-module $V$ with~$\mathbb{Q}(q,t)$-basis $\{v_i \}_{i\in I}$, the induced module
\smash{$\Ind^{\sD_n}_{\sA_n} V := \sD_n \otimes_{\sA_n} V$}
 has $\mathbb{Q}(q,t)$-basis $\{X^{\alpha}\otimes v_i | \alpha \in \mathbb{Z}^n_{\geq 0}, i \in I \}$.
\end{Remark}

\subsection{Elliptic Hall algebra}

Here we give a~brief description of the positive elliptic Hall algebra which will be a~sufficient introduction for the purposes of this paper. For a~more complete description of the elliptic Hall algebra, we direct the reader to the original paper of Burban--Schiffmann~\cite{BS} introducing EHA and the subsequent papers of Schiffmann--Vasserot~\cite{SV2,SV} which develop more of the relations between EHA and Macdonald theory.

\begin{Definition}
 For $\ell \in \mathbb{Z}\setminus \{0\}$, $r > 0$, define the special elements \smash{$P_{0, \ell}^{(n)}, P_{r,0}^{(n)} \in \sD_n^{\mathrm{sph}}$} by
 \begin{itemize}
\itemsep=0pt
 \item \smash{$P_{0,\ell}^{(n)} = \epsilon^{(n)}\bigl( \sum_{i=1}^{n} \theta_i^{\ell} \bigr) \epsilon^{(n)}$},
 \item \smash{$P_{r,0}^{(n)} = q^{r} \epsilon^{(n)} \bigl( \sum_{i=1}^{n} X_i^{r} \bigr) \epsilon^{(n)}$}.
 \end{itemize}
\end{Definition}

\begin{Theorem}[{\cite[Theorem 1.2]{SV} and \cite[Proposition 4.1]{SV2}}]
 The elements \smash{$P_{0, \ell}^{(n)}$}, \smash{$P_{r,0}^{(n)}$} for $\ell \in \mathbb{Z}\setminus \{0\}$, $r > 0$ generate $\sD_n^{\mathrm{sph}}$ as a~$\mathbb{Q}(q,t)$-algebra. There is a~unique $\mathbb{Z}_{\geq0}\times \mathbb{Z}$ grading on $\sD_n^{\mathrm{sph}}$ determined by \smash{$\deg\bigl(P_{0,\ell}^{(n)}\bigr) = (0,\ell)$} and \smash{$\deg\bigl(P_{r,0}^{(n)}\bigr) = (r,0)$}.
 There is a~graded algebra surjection \smash{$\sD_{n+1}^{\mathrm{sph}}\rightarrow \sD_n^{\mathrm{sph}}$} determined for $\ell \in \mathbb{Z}\setminus \{0\}$, $r > 0$ by
 \smash{$ P_{0,\ell}^{(n+1)} \rightarrow P_{0,\ell}^{(n)}$} and \smash{$ P_{r,0}^{(n+1)} \rightarrow P_{r,0}^{(n)}$}.
\end{Theorem}

The existence of the $\mathbb{Z}_{\geq 0}\times \mathbb{Z}$-graded algebra surjections $\sD_{n+1}^{\mathrm{sph}}\rightarrow \sD_n^{\mathrm{sph}}$ allows for the following definition.

\begin{Definition}[{\cite[Section~1.3]{SV}}]\label{EHA def}
 The \textit{positive elliptic Hall algebra} $\sE^{+}$ is the limit of the $\mathbb{Z}_{\geq0}\times \mathbb{Z}$-graded algebras \smash{$\sD_n^{\mathrm{sph}}$} with respect to the maps \smash{$\sD_{n+1}^{\mathrm{sph}}\rightarrow \sD_n^{\mathrm{sph}}$}. For $\ell \in \mathbb{Z}\setminus \{0\}$, $r > 0$, define the special elements of $\sE^{+}$, \smash{$P_{0,\ell}:= \lim_n P_{0,\ell}^{(n)}$} and \smash{$P_{r,0}:= \lim_n P_{r,0}^{(n)}$}.
\end{Definition}

The positive elliptic Hall algebra contains elements $P_{a,b}$ for $(a,b) \in \mathbb{Z}_{\geq0}\times \mathbb{Z} \setminus \{(0,0)\}$ which may be defined using repeated commutators of the elements $P_{0,\ell}$, $P_{r,0}$. For example, $P_{1,1} = [P_{0,1},P_{1,0}]$. We will not require an explicit description of these elements for the purposes of this paper. Further, we will not require knowledge of the full elliptic Hall algebra $\mathcal{E}$ which is obtained as the Drinfeld double of $\mathcal{E}^{+}$ with respect to a~certain Hopf pairing. In the standard Macdonald theory picture \cite[Corollary~1.5]{SV}, we can realize the action of the full EHA on the ring of symmetric functions $\Lambda$ using multiplication operators $p_{r}^{\bullet}$, skewing operators $p_{r}^{\perp}$, and Macdonald operators $p_{\ell}[\Delta]$ roughly corresponding to the elements $P_{r,0}$, $P_{-r,0}$, $P_{0,\ell}$, respectively.

\begin{Remark}
 We will be considering the $\mathbb{Z}_{\geq 0}$-grading on $\mathcal{E}^{+}$ obtained by the specialization~${(a,b) \rightarrow a}$, i.e., for $r > 0$ and $\ell \in \mathbb{Z}\setminus \{0\}$
 \begin{itemize}
\itemsep=0pt
 \item $\deg(P_{0,\ell}) = 0$,
 \item $\deg(P_{r,0}) = r$.
 \end{itemize} When we refer to a~$\mathcal{E}^{+}$-module $V$ as \textit{graded} we are referring to the $\mathbb{Z}_{\geq 0}$-grading on $\mathcal{E}^{+}$.
\end{Remark}

\section{DAHA modules from Young diagrams}\label{DAHA Modules from Young Diagrams}
\subsection[The D\_n-module V\_lambda]{The $\boldsymbol{\sD_n}$-module $\boldsymbol{V_{\lambda}}$}
 We begin by defining a~collection of DAHA modules indexed by Young diagrams $\lambda \in \Par$. These modules are the same as those appearing in \cite[Definition 3.1]{DL_2011} but we take the approach of using induction from $\sA_n$ to $\sD_n$ for their definition.

 \begin{Definition}
 Let $\lambda \in \Par$ with $|\lambda| = n$. Define the $\sD_n$-module $V_{\lambda}$ to be the induced module
 $V_{\lambda}:= \Ind_{\sA_n}^{\sD_n}\rho_n^{*}(S_{\lambda})$. The modules $V_{\lambda}$ naturally have the basis given by $X^{\alpha}\otimes e_{\tau}$, where $X^{\alpha}$ is a~monomial and $\tau \in \SYT(\lambda)$. We refer to this as the \textit{standard basis} of $V_{\lambda}$.
 \end{Definition}

Using the theory of intertwiners for DAHA and some combinatorics, we are able to show the following structural results. The $F_{\tau}$ appearing below are the version of the non-symmetric vv Macdonald polynomials of Dunkl--Luque following our conventions but they are \textit{not} the same polynomials (see Corollary~\ref{Corollary recovering sym vv Macd}). The result below is analogous to \cite[Theorem 4.12]{DL_2011}.

\begin{Proposition}\label{weight basis Proposition}
 There exists a~basis of $V_{\lambda}$ consisting of $\Theta^{(n)}$-weight vectors \[
 \{F_{\tau}\mid \tau \in \PSYT(\lambda)\}\] with distinct $\Theta^{(n)}$-weights such that the following hold:
 \begin{itemize}
\itemsep=0pt
 \item \smash{$\theta_i^{(n)}(F_{\tau}) = q^{w_{\tau}(i)}t^{c_{\tau}(i)}F_{\tau}$}.
 \item If $\tau \in \SYT(\lambda)$, then $F_{\tau} = 1\otimes e_{\tau}$.
 \item If $s_i(\tau) > \tau$, then
 \[
 \left( tT_{i}^{-1} + \frac{(t-1)q^{w_{\tau}(i+1)}t^{c_{\tau}(i+1)}}{q^{w_{\tau}(i)}t^{c_{\tau}(i)}-q^{w_{\tau}(i+1)}t^{c_{\tau}(i+1)}}\right)(F_{\tau}) = F_{s_i(\tau)}.
 \]
 \item $F_{\Psi(\tau)} = q^{w_1(\tau)}X_n\pi_n^{-1}(F_{\tau})$.
 \end{itemize}
\end{Proposition}

\begin{proof}
 Using the analogue of Mackey decomposition for induction between $\sA_n$ and $\sD_n$ \cite[Proposition~3.9.5]{CherednikBook}, we find
 \begin{align*}
 \Res_{\Theta^{(n)}}^{\sD_n}(V_{\lambda}) &= \Res_{\Theta^{(n)}}^{\sD_n} \Ind_{\sA_n}^{\sD_n} \rho_n^{*}(S_{\lambda})
 = \bigoplus_{\sigma \in \bigl( \widehat{\mathfrak{S}}_n/\mathfrak{S}_n \bigr)^{+}} \bigl( \Res_{\Theta^{(n)}}^{\sA_n} \rho_n^{*}(S_{\lambda}) \bigr)^{\sigma} \\
 &= \bigoplus_{\substack{\sigma \in \bigl( \widehat{\mathfrak{S}}_n/\mathfrak{S}_n \bigr)^{+}\\ \tau \in \SYT(\lambda)}} \mathbb{Q}(q,t) (\varphi_{\sigma}\otimes e_{\tau}).
 \end{align*}

 As a~consequence, we find that the set \smash{$\{\varphi_{\sigma}\otimes e_{\tau}\}_{(\sigma,\tau) \in ( \widehat{\mathfrak{S}}_n/\mathfrak{S}_n )^{+} \times \SYT(\lambda)}$} is a~generalized $\Theta^{(n)}$-weight basis for $V_{\lambda}$. Define scalars $g_{\tau}$ recursively as follows:
 \begin{itemize}
\itemsep=0pt
 \item $g_{\tau_{\lambda}^{\rm rs}}:= 1$.
 \item If $s_i(\tau)> \tau$ and $w_{\tau}(i) = w_{\tau}(i+1)$, then $g_{s_i(\tau)}:= q^{-w_{\tau}(i)}g_{\tau}$.
 \item If $s_i(\tau) > \tau$ and $w_{\tau}(i)\neq w_{\tau}(i+1)$, then \smash{$g_{s_i(\tau)}:= \left( q^{w_{\tau}(i)}t^{c_{\tau}(i)} - q^{w_{\tau}(i+1)}t^{c_\tau(i+1)} \right)^{-1} g_{\tau}$}.
 \item \smash{$g_{\Psi(\tau)}:= t^{1-n-c_{\tau}(1)}g_{\tau}$}.
 \end{itemize}
 For $\tau \in \PSYT(\lambda)$, define $F_{\tau} := g_{\tau} \varphi_{\sigma_{w_{\tau}}} \otimes e_{S(\tau)}$. We now show that the set $(F_{\tau})_{\tau \in \PSYT(\lambda)}$ constitutes a~basis of $V_{\lambda}$ satisfying the conditions of the proposition statement.

 By Proposition~\ref{bijection between PSYT and pairs on min coset reps and SYT} and since the scalars $g_{\tau}$ are non-zero by construction, it follows that the set \smash{$(F_{\tau})_{\tau \in \PSYT(\lambda)}$} is a~basis for $V_{\lambda}$ labeled by $\PSYT(\lambda)$. Further, by induction using Lemma~\ref{weight vector intertwiners Lemma} and Proposition~\ref{generating PSYT from SYT} we see that each $F_{\tau}$ is a~\smash{$\Theta^{(n)}$}-weight vector with $\smash{\theta_i^{(n)}(F_{\tau})} = q^{w_{\tau}(i)}t^{c_{\tau}(i)}F_{\tau}$.

 Suppose $s_i(\tau) > \tau$. By Proposition~\ref{additional relations for affine hecke},
 \begin{gather*}
 \left( tT_{i}^{-1} + \frac{(t-1)q^{w_{\tau}(i+1)}t^{c_{\tau}(i+1)}}{q^{w_{\tau}(i)}t^{c_{\tau}(i)}-q^{w_{\tau}(i+1)}t^{c_{\tau}(i+1)}}\right)(F_{\tau}) \\
 \qquad = \bigl( q^{w_{\tau}(i)}t^{c_{\tau}(i)}-q^{w_{\tau}(i+1)}t^{c_{\tau}(i+1)}\bigr)^{-1} \\
 \phantom{\qquad = }{} \times\bigl( tT_{i}^{-1}(q^{w_{\tau}(i)}t^{c_{\tau}(i)}-q^{w_{\tau}(i+1)}t^{c_{\tau}(i+1)}) + (t-1)q^{w_{\tau}(i+1)}t^{c_{\tau}(i+1)}\bigr)(F_{\tau}) \\
 \qquad = \bigl( q^{w_{\tau}(i)}t^{c_{\tau}(i)}-q^{w_{\tau}(i+1)}t^{c_{\tau}(i+1)}\bigr)^{-1} \bigl( tT_{i}^{-1}(\theta_{i}-\theta_{i+1}) + (t-1)\theta_{i+1}\bigr)(F_{\tau})\\
 \qquad = \bigl( q^{w_{\tau}(i)}t^{c_{\tau}(i)}-q^{w_{\tau}(i+1)}t^{c_{\tau}(i+1)}\bigr)^{-1} \varphi_i(F_{\tau}) \\
 \qquad = \bigl( q^{w_{\tau}(i)}t^{c_{\tau}(i)}-q^{w_{\tau}(i+1)}t^{c_{\tau}(i+1)}\bigr)^{-1}g_{\tau} \varphi_{i}\varphi_{\sigma_{w_{\tau}}} \otimes e_{S(\tau)}.
 \end{gather*}
 If $w_{\tau}(i) = w_{\tau}(i+1)$, then from the proof of Proposition~\ref{generating PSYT from SYT}, we know that $\varphi_i \varphi_{\sigma_{w_{\tau}}} = \varphi_{\sigma_{w_{\tau}}} \varphi_j$ where $j= \sigma^{-1}(i)$ and $\sigma$ is the minimal length permutation such that $\sigma( \sort(w_{\tau})) = w_{\tau}$. Combining this observation with Definitions~\ref{tableaux from PSYT def} and~\ref{irreps for finite hecke Definition}, we have
 \begin{gather*}
 \bigl( q^{w_{\tau}(i)}t^{c_{\tau}(i)}-q^{w_{\tau}(i+1)}t^{c_{\tau}(i+1)}\bigr)^{-1}g_{\tau} \varphi_{i}\varphi_{\sigma_{w_{\tau}}} \otimes e_{S(\tau)} \\
 \qquad =\bigl( q^{w_{\tau}(i)}t^{c_{\tau}(i)}-q^{w_{\tau}(i+1)}t^{c_{\tau}(i+1)}\bigr)^{-1}g_{\tau} \varphi_{\sigma_{w_{\tau}}}\varphi_{j} \otimes e_{S(\tau)} \\
 \qquad = \bigl( q^{w_{\tau}(i)}t^{c_{\tau}(i)}-q^{w_{\tau}(i+1)}t^{c_{\tau}(i+1)}\bigr)^{-1}g_{\tau} \varphi_{\sigma_{w_{\tau}}} \otimes \varphi_{j}e_{S(\tau)} \\
 \qquad = \bigl( q^{w_{\tau}(i)}t^{c_{\tau}(i)}-q^{w_{\tau}(i+1)}t^{c_{\tau}(i+1)}\bigr)^{-1}g_{\tau} \varphi_{\sigma_{w_{\tau}}} \otimes (t^{c_{S(\tau)}(j)}-t^{c_{S(\tau)}(j+1)}) e_{s_jS(\tau)} \\
 \qquad = q^{-w_{\tau}(i)}\bigl( t^{c_{\tau}(i)}-t^{c_{\tau}(i+1)}\bigr)^{-1}g_{\tau} \varphi_{\sigma_{w_{s_i(\tau)}}} \otimes (t^{c_{\tau}(i)}-t^{c_{\tau}(i+1)}) e_{S(s_i(\tau))} \\
 \qquad = q^{-w_{\tau}(i)}g_{\tau}\varphi_{\sigma_{w_{s_i(\tau)}}} \otimes e_{S(s_i(\tau))} = g_{s_i(\tau)}\varphi_{\sigma_{w_{s_i(\tau)}}} \otimes e_{S(s_i(\tau))} = F_{s_i(\tau)}.
 \end{gather*}

 If instead $w_{\tau}(i) \neq w_{\tau}(i+1)$, then $\varphi_{i}\varphi_{\sigma_{w_{\tau}}} = \varphi_{s_i\sigma_{w_{\tau}}} = \varphi_{\sigma_{w_{s_i(\tau)}}}$ and $S(s_i(\tau)) = S(\tau)$. Thus,
 \begin{gather*}
 \bigl( q^{w_{\tau}(i)}t^{c_{\tau}(i)}-q^{w_{\tau}(i+1)}t^{c_{\tau}(i+1)}\bigr)^{-1}g_{\tau} \varphi_{i}\varphi_{\sigma_{w_{\tau}}} \otimes e_{S(\tau)} \\
 \qquad = \bigl( q^{w_{\tau}(i)}t^{c_{\tau}(i)}-q^{w_{\tau}(i+1)}t^{c_{\tau}(i+1)}\bigr)^{-1}g_{\tau}\varphi_{\sigma_{w_{s_i(\tau)}}} \otimes e_{S(s_i(\tau))} \\
 \qquad = g_{s_i(\tau)}\varphi_{\sigma_{w_{s_i(\tau)}}} \otimes e_{S(s_i(\tau))} = F_{s_i(\tau)}.
 \end{gather*}

 Now we check the last relation
 \begin{align*}
 q^{w_{\tau}(1)}X_n\pi_{n}^{-1}(F_{\tau}) &
 = t^{1-n-c_{\tau}(1)}t^{n-1}q^{w_{\tau}(1)}t^{c_{\tau}(1)}X_n\pi_{n}^{-1}(F_{\tau}) = t^{1-n-c_{\tau}(1)} t^{n-1}X_n\pi_{n}^{-1}\theta_1(F_{\tau}) \\
 & = t^{1-n-c_{\tau}(1)} X_n\pi_{n}^{-1}\pi_nT_{n-1}\cdots T_1(F_{\tau}) =t^{1-n-c_{\tau}(1)}X_nT_{n-1}\cdots T_1(F_{\tau}) \\
 &= t^{1-n-c_{\tau}(1)}\gamma_n(F_{\tau}) = t^{1-n-c_{\tau}(1)} g_{\tau} \gamma_n\varphi_{\sigma_{w_{\tau}}} \otimes e_{S(\tau)} \\
 &= g_{\Psi(\tau)} \varphi_{\widetilde{\gamma}_n\sigma_{w_{\tau}}} \otimes e_{S(\Psi(\tau))} = g_{\Psi(\tau)} \varphi_{\sigma_{\widetilde{\gamma}_n (w_{\tau})}} \otimes e_{S(\Psi(\tau))} \\
 &= g_{\Psi(\tau)} \varphi_{\sigma_{w_{\Psi(\tau)}}} \otimes e_{S(\Psi(\tau))}
 = F_{\Psi(\tau)}. \tag*{\qed}
 \end{align*} \renewcommand{\qed}{}
\end{proof}

The above proposition allows for the following definition.

\begin{Definition}
 For $\tau \in \PSYT(\lambda)$, define the \textit{non-symmetric vv Macdonald polynomial} for $\tau$ to be $F_{\tau} \in V_{\lambda}$ as described in the proof of Proposition~\ref{weight basis Proposition}.
\end{Definition}

\begin{Example}

We have that
\[
 \ytableausetup{smalltableaux}
\ytableaushort{{1q}{2q},3}= s_2s_1 \Psi^2 \ytableausetup{smalltableaux}
\ytableaushort{{1}{2},3}\,.
\]
 By starting with
 \[
 F_{\ytableausetup{smalltableaux}
\ytableaushort{{1}{2},3}} = 1\otimes e_{\ytableausetup{smalltableaux}
\ytableaushort{{1}{2},3}}\]
 and using the recurrence relations in Proposition~\ref{weight basis Proposition}, we find that
\begin{align*}
 F_{~\ytableausetup{smalltableaux}
\ytableaushort{{1q}{2q},3}} ={}& t^{-2}X_1X_2\otimes e_{~\ytableausetup{smalltableaux}
\ytableaushort{12,3}} + t^{-2} \left( \frac{1-t}{1-qt^2} \right)X_2X_3 \otimes e_{~\ytableausetup{smalltableaux}
\ytableaushort{13,2}} \\
& + \frac{t^{-2}}{1+t} \left(\frac{1-t}{1-qt^2} \right) X_2X_3 \otimes e_{~\ytableausetup{smalltableaux}
\ytableaushort{12,3}} -t^{-3}\left( \frac{1-t}{1-qt^{2}} \right)X_1X_3\otimes e_{~\ytableausetup{smalltableaux}
\ytableaushort{13,2}} \\
&+ \frac{t^{-1}}{1+t} \left( \frac{1-t}{1-qt^2} \right) X_1X_3 \otimes e_{~\ytableausetup{smalltableaux}
\ytableaushort{12,3}}\,.
\end{align*}
\end{Example}

\begin{Remark}\label{action of gamma}
 Note that from Proposition~\ref{weight basis Proposition}, we get that
\smash{$ \gamma_n(F_{\tau}) = t^{n-1 + c_{\tau}(1)}F_{\Psi(\tau)}$}.
 By induction, we see that for $r \geq 1$,
\smash{$\gamma_n^{r}(F_{\tau}) = t^{r(n-1)}t^{c_{\tau}(1)+\cdots + c_{\tau}(r)}F_{\Psi^{r}(\tau)}$}.
\end{Remark}

Now that we have a~concrete combinatorial description for the module $V_{\lambda}$, we study its decomposition into $\sA_n$-submodules.

\begin{Proposition}\label{Mackey decomposition}
 The $\sD_n$-module $V_{\lambda}$ has the following decomposition into $\sA_n$-submodules:
 \[
 \Res^{\sD_n}_{\sA_n} V_{\lambda} = \bigoplus_{T \in \RYT(\lambda)} U_T,
 \]
 where $U_T:= \text{span}_{\mathbb{Q}(q,t)}\{F_{\tau}\mid \tau \in \PSYT(\lambda;T)\}$.
 Further, each $\sA_n$-module $U_T$ is irreducible.
\end{Proposition}

\begin{proof}
 Let $T \in \RYT(\lambda)$. Note that it follows immediately from Proposition~\ref{weight basis Proposition} that each~$U_T$ is a~$\sA_n$-submodule of $V_{\lambda}$. Further, trivially $U_T \cap U_{T'} = \varnothing$ for $T \neq T'$ since the $F_{\tau}$ are a~basis for $V_{\lambda}$ and the sets $\PSYT(\lambda;T)$ partition $\PSYT(\lambda)$. Therefore,
 \[
 \Res^{\sD_n}_{\sA_n} V_{\lambda} = \bigoplus_{T \in \RYT(\lambda)} U_T.
 \]
 Now, let $T \in \RYT(\lambda)$. If $U \subset U_{T}$ is a~nonzero $\sA_n$-submodule, then $U$ must contain some $\Theta^{(n)}$-weight vector as $U_T$ is spanned by $\Theta^{(n)}$-weight vectors. Thus there exists some ${\tau_0 \!\in\! \PSYT(\lambda;T)}$ with $F_{\tau_0} \in U$. But then it follows readily from Proposition~\ref{weight basis Proposition}, that by using intertwiner operators $\varphi_i$ given any $\tau \in \PSYT(\lambda)$ we may find $A \in \sA_n$ such that $A(F_{\tau_0}) = F_{\tau}$. Therefore, $U = U_T$ and hence $U_T$ is irreducible.
\end{proof}

\begin{Remark}\label{PSYT orbits using Bruhat}
 It follows by using Frobenius reciprocity and Proposition~\ref{Mackey decomposition} that in fact there are surjective $\sA_n$-module maps
\smash{$ \Ind_{\sA_{\mu(T)}}^{\sA_n} \chi_T \rightarrow U_T$},
 where $\chi_T$ is the $1$-dimensional representation of $\sA_{\mu(T)}$ determined by the $\Theta^{(n)}$-weight of $F_{\Min(T)}$ and $T_i \rightarrow 1$ for relevant $T_i$. Thus each~$U_T$ is a~quotient of an induced module from a~parabolic subalgebra of $\sA_n$. In the case of~${T \in \RSSYT(\lambda)}$, this map is an isomorphism. We may witness the implied bijection between~${\PSYT(\lambda;T)}$ and $\mathfrak{S}_n/\mathfrak{S}_{\mu(T)}$ combinatorially using the map $\sigma \rightarrow \sigma(\Min(T))$
 for~${\sigma \in \mathfrak{S}_n/\mathfrak{S}_{\mu(T)}}$. It is straightforward to check by decomposing~$\lambda$ into horizontal strip diagrams where $T$ is constant along rows that this map is actually an isomorphism of posets.
\end{Remark}

The following lemma exhibits triangularity for the $T_i^{-1}$ operators with respect to the reversed Bruhat order on $\mathbb{Z}_{\geq 0}^{n}$.

\begin{Lemma}\label{action of T inverse}
 For $1\leq i \leq n-1$ and $a \geq 0$,
 \[
 \bigl(tT_i^{-1}\bigr)X_{i+1}^{a} = X_i^{a}\bigl(tT_i^{-1}\bigr) + (t-1)X_{i+1}\frac{X_i^{a}-X_{i+1}^{a}}{X_i-X_{i+1}}.
 \]
 Further, every monomial occurring in the term \smash{$X_{i+1}\frac{X_i^{a}-X_{i+1}^{a}}{X_i-X_{i+1}}$} is strictly lower than $X_i^{a}$ with respect to the Bruhat ordering $\preceq$. Consequently, it follows that for any $\alpha \in \mathbb{Z}^n_{\geq 0}$ with $s_i(\alpha) \succeq \alpha$ the following expansion holds for some scalars $c_{\beta}$:
 \[
 \bigl(tT_i^{-1}\bigr)X^{\alpha} = X^{s_i(\alpha)}\bigl(tT_i^{-1}\bigr) + \sum_{\beta \prec s_i(\alpha)} c_{\beta}X^{\beta}.
 \]
\end{Lemma}

\begin{proof}
We start with
 \begin{align*}
 tT_i^{-1}X_{i+1}^{a}&= (T_i+t-1)X_{i+1}^a = T_iX_{i+1}^{a} + (t-1)X_{i+1}^a \\
 &= X_i^{a}T_i + (1-t)X_i\frac{X_{i+1}^{a}-X_i^{a}}{X_i-X_{i+1}} -(1-t)X_{i+1}^{a} \\
 & = X_i^a\bigl(tT_i^{-1}+1-t\bigr) + (1-t)X_i\frac{X_{i+1}^a - X_i^a}{X_i-X_{i+1}} -(1-t)X_{i+1}^{a} \\
 & = X_i^{a}tT_i^{-1} + (1-t)X_i^{a} -(1-t)X_{i+1}^{a} + (1-t)X_i\frac{X_{i+1}^{a}-X_i^{a}}{X_i-X_{i+1}}\\
 &= X_i^{a}tT_i^{-1} + (t-1)X_{i+1}\frac{X_{i+1}^{a}-X_i^{a}}{X_{i+1}-X_i}.
 \end{align*}
Further,
\[
 X_{i+1}\frac{X_{i+1}^{a}-X_i^{a}}{X_{i+1}-X_i} = X_{i+1}^{a} + X_{i+1}^{a-1}X_i + \cdots + X_{i+1}^2X_i^{a-2} + X_{i+1}X_i^{a-1}
\]
 so that
\[
 tT_i^{-1}X_{i+1}^{a} = X_i^{a}tT_i^{-1} + (t-1)\bigl(X_{i+1}^{a} + X_{i+1}^{a-1}X_i + \cdots + X_{i+1}^2X_i^{a-2} + X_{i+1}X_i^{a-1}\bigr).
\]

Now let $\alpha \in \mathbb{Z}_{\geq 0}^{n}$ with $s_i(\alpha) \succ \alpha$, i.e., $\alpha_i < \alpha_{i+1}$.
Then
\begin{align*}
 tT_i^{-1}X^{\alpha} ={}& tT_i^{-1}X_1^{\alpha_1}\cdots X_{i-1}^{\alpha_{i-1}}X_i^{\alpha_i}X_{i+1}^{\alpha_{i+1}}X_{i+2}^{\alpha_{i+2}}\cdots X_n^{\alpha_n} \\
 ={}& X_1^{\alpha_1}\cdots X_{i-1}^{\alpha_{i-1}}X_{i+2}^{\alpha_{i+2}}\cdots X_n^{\alpha_n} tT_i^{-1}X_i^{\alpha_i}X_{i+1}^{\alpha_{i+1}} \\
 ={}& X_1^{\alpha_1}\cdots X_{i-1}^{\alpha_{i-1}}X_i^{\alpha_i}X_{i+1}^{\alpha_{i}}X_{i+2}^{\alpha_{i+2}}\cdots X_n^{\alpha_n} tT_i^{-1}X_{i+1}^{\alpha_{i+1}-\alpha_{i}} \\
 ={}& X_1^{\alpha_1}\cdots X_{i-1}^{\alpha_{i-1}}X_i^{\alpha_i}X_{i+1}^{\alpha_{i}}X_{i+2}^{\alpha_{i+2}}\cdots X_n^{\alpha_n} \\
 & \times \bigl( X_i^{\alpha_{i+1}-\alpha_i}tT_i^{-1} + (t-1)\bigl( X_{i+1}^{\alpha_{i+1}-\alpha_i} + \dots + X_{i+1}X_i^{\alpha_{i+1}-\alpha_i-1}\bigr) \bigr) \\
 ={}& X^{s_i(\alpha)}tT_i^{-1} + (t-1)\sum_{j=0}^{\alpha_{i+1}-\alpha_i-1} X^{\alpha + j(e_i - e_{i+1})}.
\end{align*}
 Lastly, from Definition~\ref{Bruhat def} it is clear that for all $0 \leq j \leq \alpha_{i+1}-\alpha_i-1$, $s_i(\alpha) \succ \alpha + j(e_i - e_{i+1})$.
\end{proof}

Using the above lemma, we are able to partially describe the expansions of the $F_{\tau}$ into the standard basis of $V_{\lambda}$.

\begin{Corollary}\label{triangular expansion}
For $\tau \in \PSYT(\lambda)$, each $F_{\tau}$ has a~triangular monomial expansion with respect to the Bruhat order on $\mathbb{Z}^n_{\geq 0}$ of the form
\smash{$ F_{\tau} = X^{w_{\tau}}\otimes f(\tau) + \sum_{\beta \prec w_{\tau}} X^{\beta}\otimes v_{\beta}
$}
 for some $v_{\beta} \in S_\lambda$ where $f(\tau) \in S_\lambda$ is given by the following recurrence relations:

\begin{itemize}
\itemsep=0pt
 \item If $\tau \in \SYT(\lambda)$, then $f(\tau) = e_{\tau}$.
 \item $f(\Psi(\tau)) = t^{-(n-1)}T_{n-1}\cdots T_1(f(\tau))$.
 \item If $w_{\tau}(i) < w_{\tau}(i+1)$, then $f(s_i(\tau)) = tT_i^{-1}f(\tau)$.
 \item If $w_{\tau}(i) = w_{\tau}(i+1)$ and $c_{\tau}(i)-c_{\tau}(i+1) > 1$, then
 \[
 f(s_i(\tau)) = \left( tT_i^{-1} + \frac{(t-1)t^{c_{\tau}(i+1)}}{t^{c_{\tau}(i)} - t^{c_{\tau}(i+1)}} \right)(f(\tau)).
 \]

\end{itemize}

\end{Corollary}
\begin{proof}
 We proceed by induction with respect to the partial ordering on $\PSYT(\lambda)$ defined in Definition~\ref{ordering and operations on tableaux defs}. At the same time, we verify the recurrence relations given for $f(\tau) \in S_{\lambda}$ as above.

 From Proposition~\ref{weight basis Proposition}, we know that if $\tau \in \SYT(\lambda)$, then $F_{\tau} = 1\otimes e_{\tau}$. Hence, $F_{\tau}$ trivially has a~triangular monomial expansion of the correct form in this case and that $f(\tau) = e_{\tau}$.
 In what follows, assume that for $\tau \in \PSYT(\lambda)$ we have that
$
 F_{\tau} = X^{w_{\tau}}\otimes f(\tau) + \smash{\sum_{\beta \prec w_{\tau}} X^{\beta}}\otimes v_{\beta}
$
 for some $v_{\beta} \in S_\lambda$.

 First, we see that
 \begin{align*}
 F_{\Psi(\tau)}& = q^{w_1(\tau)}X_n \pi_n^{-1}(F_{\tau}) = q^{w_1(\tau)}X_n \pi_n^{-1}X^{w_{\tau}}\otimes f(\tau) + \sum_{\beta \prec w_{\tau}} q^{w_1(\tau)}X_n \pi_n^{-1}X^{\beta}\otimes v_{\beta} \\
 & = q^{w_1(\tau)}q^{-w_1(\tau)} X^{\widetilde{\gamma}(w_{\tau})} \pi_n^{-1} \otimes f(\tau) + \sum_{\beta \prec w_{\tau}} q^{w_1(\tau)}q^{-\beta_1}X^{\widetilde{\gamma}(\beta)}\pi_n^{-1}\otimes v_{\beta} \\
 & = X^{\widetilde{\gamma}(w_{\tau})} \otimes \rho_n\bigl(\pi_n^{-1}\bigr)f(\tau) + \sum_{\beta \prec w_{\tau}} X^{\widetilde{\gamma}(\beta)} \otimes q^{w_1(\tau)-\beta_1} \rho_n\bigl(\pi_n^{-1}\bigr)v_{\beta} \\
 &= X^{\widetilde{\gamma}(w_{\tau})} \otimes t^{-(n-1)}T_{n-1}\cdots T_1f(\tau) + \sum_{\beta \prec w_{\tau}} X^{\widetilde{\gamma}(\beta)} \otimes q^{w_1(\tau)-\beta_1} t^{-(n-1)}T_{n-1}\cdots T_1v_{\beta}.
 \end{align*}

 From Lemma~\ref{gamma preserves Bruhat}, we know that if $\beta \prec w_{\tau}$, then $\widetilde{\gamma}(\beta) \prec \widetilde{\gamma}(w_{\tau})$. Therefore, we find that~$F_{\Psi(\tau)}$ has the expansion
 \[
 F_{\Psi(\tau)} = X^{\widetilde{\gamma}(w_{\tau})} \otimes t^{-(n-1)}T_{n-1}\cdots T_1f(\tau) + \sum_{\beta \prec \widetilde{\gamma}(\tau)} X^{\beta} \otimes v'_{\beta}
 \]
 for some $v'_{\beta} \in S_{\lambda}$. From this, we see that $f(\Psi(\tau)) = t^{-(n-1)}T_{n-1}\cdots T_1(f(\tau))$.

 Now suppose $s_i(\tau) > \tau$. From Proposition~\ref{weight basis Proposition}, we get
 \begin{align*}
 F_{s_i(\tau)} ={}& \left( tT_{i}^{-1} + \frac{(t-1)q^{w_{\tau}(i+1)}t^{c_{\tau}(i+1)}}{q^{w_{\tau}(i)}t^{c_{\tau}(i)}-q^{w_{\tau}(i+1)}t^{c_{\tau}(i+1)}}\right)(F_{\tau}) \\
 ={}& \left( tT_{i}^{-1} + \frac{(t-1)q^{w_{\tau}(i+1)}t^{c_{\tau}(i+1)}}{q^{w_{\tau}(i)}t^{c_{\tau}(i)}-q^{w_{\tau}(i+1)}t^{c_{\tau}(i+1)}}\right) \biggl(X^{w_{\tau}}\otimes f(\tau) + \sum_{\beta \prec w_{\tau}} X^{\beta}\otimes v_{\beta} \biggr) \\
 ={}& tT_{i}^{-1}\biggl(X^{w_{\tau}}\otimes f(\tau) + \sum_{\beta \prec w_{\tau}} X^{\beta}\otimes v_{\beta} \biggr) \\
 & + \left(\frac{(t-1)q^{w_{\tau}(i+1)}t^{c_{\tau}(i+1)}}{q^{w_{\tau}(i)}t^{c_{\tau}(i)}-q^{w_{\tau}(i+1)}t^{c_{\tau}(i+1)}}\right) \biggl(X^{w_{\tau}}\otimes f(\tau) + \sum_{\beta \prec w_{\tau}} X^{\beta}\otimes v_{\beta} \biggr).
 \end{align*}

 For any $\beta < w_{\tau}$, using Lemma~\ref{action of T inverse}, we find that
 \[
 tT_i^{-1}X^{\beta} \otimes v_{\beta} = \sum_{\beta' \prec s_i(w_{\tau})}X^{\beta'} \otimes u_{\beta',\beta}
 \]
 for some $u_{\beta',\beta} \in S_{\lambda}$; that is to say, each of the monomials $X^{\beta'}$ that appears in the standard basis expansion of $tT_i^{-1}X^{\beta} \otimes v_{\beta}$ must have $\beta'\prec s_i(w_{\tau})$.

 Assume $w_{\tau}(i)< w_{\tau}(i+1)$. By Lemma~\ref{action of T inverse}, we see
 \begin{align*}
 \bigl(tT_i^{-1}\bigr)X^{w_{\tau}}\otimes f(\tau) &= X^{s_i(w_{\tau})}\bigl(tT_i^{-1}\bigr)\otimes f(\tau) + \sum_{\beta \prec s_i(w_{\tau})} c_{\beta}X^{\beta}\otimes f(\tau) \\
 & = X^{s_i(w_{\tau})}\otimes \bigl(tT_i^{-1}\bigr)f(\tau) + \sum_{\beta \prec s_i(w_{\tau})} c_{\beta}X^{\beta}\otimes f(\tau).
 \end{align*}
 Therefore, $F_{s_i(\tau)}$ has the expansion
 \[
 F_{s_i(\tau)} = X^{s_i(w_{\tau})}\otimes tT_i^{-1}f(\tau) + \sum_{\beta \prec s_i(w_{\tau})} X^{\beta}\otimes v'_{\beta},
 \]
 where $v'_{\beta} \in S_{\lambda}$. Since $s_i(w_{\tau}) = w_{s_{i}(\tau)}$, we have
 \[
 F_{s_i(\tau)} = X^{w_{s_{i}(\tau)}}\otimes tT_i^{-1}f(\tau) + \sum_{\beta \prec w_{s_{i}(\tau)}} X^{\beta}\otimes v'_{\beta}
 \]
 and $f(s_i(\tau)) = tT_i^{-1}f(\tau)$.

 Now assume instead that $w_{\tau}(i) = w_{\tau}(i+1)$ and $c_{\tau}(i) - c_{\tau}(i+1) > 1$. Then ${T_iX^{w_{\tau}} = X^{w_{\tau}}T_i}$~so
 \begin{align*}
 F_{s_i(\tau)} = {}&\left( tT_{i}^{-1} + \frac{(t-1)q^{w_{\tau}(i+1)}t^{c_{\tau}(i+1)}}{q^{w_{\tau}(i)}t^{c_{\tau}(i)}-q^{w_{\tau}(i+1)}t^{c_{\tau}(i+1)}}\right)(F_{\tau}) \\
 = {}& \left( tT_{i}^{-1} + \frac{(t-1)t^{c_{\tau}(i+1)}}{t^{c_{\tau}(i)}-t^{c_{\tau}(i+1)}}\right)(F_{\tau}) \\
 = {}& \left( tT_{i}^{-1} + \frac{(t-1)t^{c_{\tau}(i+1)}}{t^{c_{\tau}(i)}-t^{c_{\tau}(i+1)}}\right) \biggl(X^{w_{\tau}}\otimes f(\tau) + \sum_{\beta \prec w_{\tau}} X^{\beta}\otimes v_{\beta} \biggr) \\
 = {}& \left( tT_{i}^{-1} + \frac{(t-1)t^{c_{\tau}(i+1)}}{t^{c_{\tau}(i)}-t^{c_{\tau}(i+1)}}\right)X^{w_{\tau}}\otimes f(\tau) \\
 &+ \left( tT_{i}^{-1} + \frac{(t-1)t^{c_{\tau}(i+1)}}{t^{c_{\tau}(i)}-t^{c_{\tau}(i+1)}}\right)\sum_{\beta \prec w_{\tau}} X^{\beta}\otimes v_{\beta} \\
 = {}& X^{w_{\tau}}\left( tT_{i}^{-1} + \frac{(t-1)t^{c_{\tau}(i+1)}}{t^{c_{\tau}(i)}-t^{c_{\tau}(i+1)}}\right)\otimes f(\tau) \\
 &+ \left( tT_{i}^{-1} + \frac{(t-1)t^{c_{\tau}(i+1)}}{t^{c_{\tau}(i)}-t^{c_{\tau}(i+1)}}\right)\sum_{\beta \prec w_{\tau}} X^{\beta}\otimes v_{\beta} \\
 = {}& X^{w_{\tau}} \otimes \left( tT_{i}^{-1} + \frac{(t-1)t^{c_{\tau}(i+1)}}{t^{c_{\tau}(i)}-t^{c_{\tau}(i+1)}}\right) f(\tau)\\
 &+ \left( tT_{i}^{-1} + \frac{(t-1)t^{c_{\tau}(i+1)}}{t^{c_{\tau}(i)}-t^{c_{\tau}(i+1)}}\right)\sum_{\beta \prec w_{\tau}} X^{\beta}\otimes v_{\beta}.
 \end{align*}
 Therefore, since $w_{\tau} = w_{s_i(\tau)}$ we find that
 \[
 F_{s_i(\tau)} = X^{w_{s_i(\tau)}}\otimes \left( tT_{i}^{-1} + \frac{(t-1)t^{c_{\tau}(i+1)}}{t^{c_{\tau}(i)}-t^{c_{\tau}(i+1)}}\right) f(\tau) + \sum_{\beta \prec w_{s_i(\tau)}} X^{\beta}\otimes v'_{\beta}
 \]
 for some $v'_{\beta} \in S_{\lambda}$ and
 \[
 f(s_i(\tau)) = \left( tT_i^{-1} + \frac{(t-1)t^{c_{\tau}(i+1)}}{t^{c_{\tau}(i)} - t^{c_{\tau}(i+1)}} \right)(f(\tau)).\tag*{\qed}
 \] \renewcommand{\qed}{}
\end{proof}

Using the $\zeta_i$ operators on $\PSYT(\lambda)$, we may compute $f(\Top(T))$ explicitly.

\begin{Proposition}\label{creation operators action}
 For $T \in \RYT(\lambda)$, we have that
 \[
 f(\Top(T)) = \mathscr{C}_1^{\nu(T)_1 - \nu(T)_2}\cdots \mathscr{C}_n^{\nu(T)_n}(e_{S(T)}),
 \]
 where for $1 \leq i \leq n$, we define
 $\mathscr{C}_i:= \left( t^{1-i} T_{i-1}\cdots T_1 \right)^{i}$.
\end{Proposition}

\begin{proof}
 Using the recurrence relations in Corollary~\ref{triangular expansion} for the elements $f(\tau)$ and Proposition~\ref{generating PSYT from SYT}, we see that for any $T\in \RSSYT(\lambda)$ since
 \[
 \Top(T) = \zeta_1^{\nu(T)_1 - \nu(T)_2}\cdots \zeta_n^{\nu(T)_n}(S(T))
 \]
 with each $\zeta_i:= (s_i\cdots s_{n-1}\Psi)^{i}$, we have a~similar expression for $f(\Top(T))$
 \[
 f(\Top(T))= \mathscr{C}_1^{\nu(T)_1 - \nu(T)_2}\cdots \mathscr{C}_n^{\nu(T)_n}(e_{S(T)}),
 \]
 where \smash{$\mathscr{C}_i :=\bigl( \bigl(tT_i^{-1}\bigr)\cdots \bigl(tT_{n-1}^{-1}\bigr) \bigl(t^{-(n-1)}T_{n-1}\cdots T_1\bigr) \bigr)^{i}=\bigl( t^{1-i} T_{i-1}\cdots T_1 \bigr)^{i}$} is obtained by replacing each $s_j$ and $\Psi$ in $\zeta_i$ with $tT_j^{-1}$ and $t^{-(n-1)}T_{n-1}\cdots T_1$, respectively. Importantly, when we apply $\zeta_i$ to any element of the form $\Top(T')$ we never perform any swaps $s_j(\tau) > \tau$ such that~${w_{\tau}(j) = w_{\tau}(j+1)}$ and hence never require the more complicated recurrence relation
 \[
 f(s_j(\tau)) = \left( tT_j^{-1} + \frac{(t-1)t^{c_{\tau}(j+1)}}{t^{c_{\tau}(j)} - t^{c_{\tau}(j+1)}} \right)(f(\tau)).\tag*{\qed}
 \] \renewcommand{\qed}{}
\end{proof}

The $\mathscr{C}_i$ operators can be identified concretely using the $\overline{\theta}_j$ elements of the finite Hecke algebra.

\begin{Lemma}\label{creation operators Lemma}
 For all $1 \leq i \leq n$,
$\mathscr{C}_i = A_i\cdots A_1$, where $A_j:= t^{-(j-1)}\overline{\theta}_j^{-1}$.
\end{Lemma}
\begin{proof}
 Let $0 \leq k \leq i-1$. We first show by induction that
 \[
 \bigl(t^{i-1}\overline{\theta}_{i}^{-1}\bigr)\cdots \bigl(t^{i-k-1}\overline{\theta}_{i-k}^{-1}\bigr) = (T_{i-1}\cdots T_1)^{k+1} (T_{k+1}\cdots T_{i-1})(T_k\cdots T_{i-2})\cdots (T_1\cdots T_{i-k-1}).
 \]

 To start, we see that for $k=0$ we have
 \[
 t^{i-1}\overline{\theta}_{i}^{-1} = T_{i-1}\cdots T_1^2 \cdots T_{i-1} = (T_{i-1}\cdots T_1)^{1} (T_1\cdots T_{i-1}).
 \]

 Now suppose that for $0 \leq k \leq i-2$ the formula above holds. Then
 \begin{gather*}
 \bigl(t^{i-1}\overline{\theta}_{i}^{-1}\bigr)\cdots \bigl(t^{i-(k+1)}\overline{\theta}_{i-(k+1)}^{-1}\bigr) \\
 \qquad = (T_{i-1}\cdots T_1)^{k+1} (T_{k+1}\cdots T_{i-1})(T_k\cdots T_{i-2})\cdots (T_1\cdots T_{i-k-1})\bigl(t^{i-(k+1)}\overline{\theta}_{i-(k+1)}^{-1}\bigr) \\
 \qquad = (T_{i-1}\cdots T_1)^{k+1} (T_{k+1}\cdots T_{i-1})(T_k\cdots T_{i-2})\cdots (T_1\cdots T_{i-k-1})\\
 \phantom{ \qquad = }{} \times\bigl(T_{i-k-2}\cdots T_1^2 \cdots T_{i-k-2}\bigr) \\
 \qquad = (T_{i-1}\cdots T_1)^{k+1} (T_{k+1}\cdots T_{i-1})(T_k\cdots T_{i-2})\cdots (T_2\cdots T_{i-k}) \\
 \phantom{ \qquad = }{} \times (T_1\cdots T_{i-k-1}) (T_{i-k-2}\cdots T_1)(T_1 \cdots T_{i-k-2}) \\
 \qquad = (T_{i-1}\cdots T_1)^{k+1} (T_{k+1}\cdots T_{i-1})(T_k\cdots T_{i-2})\cdots (T_2\cdots T_{i-k}) \\
 \phantom{ \qquad = }{} \times (T_{i-k-1}\cdots T_{2})(T_{1}\cdots T_{i-k-1})(T_1 \cdots T_{i-k-2}) \\
 \qquad = (T_{i-1}\cdots T_1)^{k+1} (T_{k+1}\cdots T_{i-1})(T_k\cdots T_{i-2})\cdots (T_2\cdots T_{i-k}) \\
 \phantom{ \qquad = }{} \times (T_{i-k-1}\cdots T_{1})(T_{2}\cdots T_{i-k-1})(T_1 \cdots T_{i-k-2}) \\
 \qquad = (T_{i-1}\cdots T_1)^{k+1} (T_{k+1}\cdots T_{i-1})(T_k\cdots T_{i-2})\cdots (T_3\cdots T_{i-k+1}) \\
 \phantom{ \qquad = }{} \times (T_2\cdots T_{i-k})(T_{i-k-1}\cdots T_{1})(T_{2}\cdots T_{i-k-1})(T_1 \cdots T_{i-k-2}) \\
 \qquad = (T_{i-1}\cdots T_1)^{k+1} (T_{k+1}\cdots T_{i-1})(T_k\cdots T_{i-2})\cdots (T_3\cdots T_{i-k+1}) \\
 \phantom{ \qquad = }{} \times (T_{i-k}\cdots T_{1})(T_{3}\cdots T_{i-k})(T_{2}\cdots T_{i-k-1})(T_1 \cdots T_{i-k-2}) \\
 \qquad = \cdots \\
 \qquad = (T_{i-1}\cdots T_1)^{k+1} (T_{i-1}\cdots T_1)(T_{k+2}\cdots T_{i-1})(T_{k+1}\cdots T_{i-2})\cdots (T_1\cdots T_{i-k-2}) \\
 \qquad = (T_{i-1}\cdots T_1)^{k+2}(T_{k+2}\cdots T_{i-1})(T_{k+1}\cdots T_{i-2})\cdots (T_1\cdots T_{i-k-2}).
 \end{gather*}

 By taking $k = i-1$, we find
\smash{$ \bigl(t^{i-1}\overline{\theta}_{i}^{-1}\bigr)\cdots \bigl(t^{0}\overline{\theta}_{1}^{-1}\bigr) = (T_{i-1}\cdots T_1)^{i}$}.
 Now we see that
 \begin{align*}
 \mathscr{C}_i& = \bigl( \bigl(tT_i^{-1}\bigr)\cdots \bigl(tT_{n-1}^{-1}\bigr) \bigl(t^{-(n-1)}T_{n-1}\cdots T_1\bigr) \bigr)^{i}
 = t^{-i(i-1)} (T_{i-1}\cdots T_1)^{i} \\
 & = t^{-i(i-1)} \bigl(t^{i-1}\overline{\theta}_{i}^{-1}\bigr)\cdots \bigl(t^{0}\overline{\theta}_{1}^{-1}\bigr)
 = t^{-2(i-1)-2(i-2)-\dots -2(1)-2(0)} \bigl(t^{i-1}\overline{\theta}_{i}^{-1}\bigr)\cdots \bigl(t^{0}\overline{\theta}_{1}^{-1}\bigr) \\
 & = \bigl(t^{-(i-1)}\overline{\theta}_{i}^{-1}\bigr)\cdots \bigl(t^{-0}\overline{\theta}_{1}^{-1}\bigr)
 = A_{i}\cdots A_1 ,
 \end{align*}
 where \smash{$A_j:= t^{-(j-1)}\overline{\theta}_j^{-1}$}.
\end{proof}

It will be important for our purposes later (Proposition~\ref{stability for locally maximal psyt}, for instance) that we are able to understand not just the $\Theta^{(n)}$-weight spaces spanned by each $F_{\tau}$ but rather exactly where this vector lies in that space. The next result determines a~specific scaling factor for the $F_{\Top(T)}$ in terms of the standard basis of $V_{\lambda}$. Putting the results of this section together gives the following.

\begin{Corollary}\label{triangular expansion of top element}
 For $T \in \RYT(\lambda)$, the triangular expansion of $F_{\Top(T)}$ has the form
 \[
 F_{\Top(T)}= t^{-b_T}X^{\nu(T)}\otimes e_{S(T)} + \sum_{\beta \prec \nu(T)} X^{\beta} \otimes v_{\beta}
 \]
 for some $v_{\beta} \in S_\lambda$.
\end{Corollary}
\begin{proof}
 First, notice that for $T \in \RYT(\lambda)$
 $w_{\Top(T)} = \nu(T)$. From Proposition~\ref{creation operators action} and Lemma~\ref{creation operator Lemma},
 \begin{align*}
 f(\Top(T)) ={}& \mathscr{C}_1^{\nu(T)_1 - \nu(T)_2}\cdots \mathscr{C}_n^{\nu(T)_n}(e_{S(T)}) \\
 ={}& A_1^{\nu(T)_1 - \nu(T)_2}(A_1A_2)^{\nu(T)_2-\nu(T)_3}\cdots (A_1\cdots A_n)^{\nu(T)_n}(e_{S(T)}) \\
 ={}& A_1^{(\nu(T)_1 - \nu(T)_2)+\cdots +(\nu(T)_{n-1}-\nu(T)_{n}) + \nu(T)_{n}}\cdots A_{n-1}^{(\nu(T)_{n-1}-\nu(T)_{n}) + \nu(T)_{n}}A_n^{\nu(T)_n}\\
 &\times (e_{S(T)}) \\
 ={}& \bigl(\overline{\theta}_1^{-1}\bigr)^{\nu(T)_1 }\cdots \bigl(\overline{\theta}_n^{-1}\bigr)^{\nu(T)_n} (e_{S(T)}) \\
 ={}& t^{-\nu(T)_1(c_{S(T)}(1)-(1-1))}\cdots t^{-\nu(T)_n(c_{S(T)}(n)-(n-1))} e_{S(T)} \\
 ={}& t^{-\sum_{i=1}^{n} \nu(T)_i( c_{S(T)}(i) + i-1)}e_{S(T)}
 = t^{-b_T} e_{S(T)}.
 \end{align*}
 Therefore, the leading term of $F_{\Top(T)}$ is
 \[
 X^{w_{\Top(T)}} \otimes f(\Top(T)) =t^{-b_T} X^{\nu(T)} \otimes e_{S(T)}.\tag*{\qed}
 \]\renewcommand{\qed}{}
\end{proof}

\subsection[Connecting maps between V\_lambda\^{}(n)]{Connecting maps between $\boldsymbol{V_{\lambda^{(n)}}}$}

Here we define connecting maps between the modules $V_{\lambda^{(n)}}$ which satisfy many exceptional properties. These connecting maps are certain lifts of the maps \smash{$\mathfrak{q}^{(n)}_{\lambda}$} from Lemma~\ref{finite hecke map} and will be instrumental in the construction of the Murnaghan-type representations of $\mathcal{E}^{+}$. We will show that the our connecting maps are well behaved with respect to both the standard basis and the non-symmetric vv Macdonald polynomial bases of $V_{\lambda}$.

\begin{Definition}\label{connecting maps def}
 Let $\lambda \in \Par$. For $n \geq n_{\lambda}$, define \smash{$\Phi^{(n)}_{\lambda}\colon V_{\lambda^{(n+1)}} \rightarrow V_{\lambda^{(n)}}$} as the $\mathbb{Q}(q,t)$-linear map given on any element $X^{\alpha}\otimes v \in V_{\lambda^{(n+1)}}$ by
 \[
 \Phi^{(n)}_{\lambda}(X^{\alpha}\otimes v) = \mathbbm{1}(\alpha_{n+1} = 0) X_1^{\alpha_1}\cdots X_n^{\alpha_n}\otimes \mathfrak{q}^{(n)}_{\lambda}(v).
 \]

\end{Definition}

 It is non-trivial and somewhat counterintuitive that the $\sD_n$-modules $V_{\lambda^{(n)}}$ should be related in a~direct way given that the algebras $\sD_n$ fail to include naturally as subalgebras (unlike the~$\sH_n$ and~$\sA_n$). However, the next proposition shows that the connecting maps \smash{$\Phi^{(n)}_{\lambda}$} satisfy multiple exceptional properties upon which the remainder of this paper heavily relies.

\begin{Proposition}\label{relations for connecting maps}
 The map \smash{$\Phi^{(n)}_{\lambda}$} satisfies the following relations:
 \begin{itemize}
\itemsep=0pt
 \item \smash{$\Phi^{(n)}_{\lambda}T_i = T_i\Phi^{(n)}_{\lambda}$} for $1\leq i \leq n-1$.
 \item \smash{$\Phi^{(n)}_{\lambda}X_i = X_i\Phi^{(n)}_{\lambda}$} for $1\leq i \leq n$.
 \item \smash{$\Phi^{(n)}_{\lambda}X_{n+1} = 0$}.
 \item \smash{$\Phi^{(n)}_{\lambda} t^{-n}\pi_{n+1}T_n = t^{-(n-1)}\pi_n \Phi^{(n)}_{\lambda}$}.
 \item \smash{$\Phi^{(n)}_{\lambda}\theta_i^{(n+1)} = \theta_i^{(n)}\Phi^{(n)}_{\lambda}$} for $1\leq i \leq n$.
 \item \smash{$\Phi^{(n)}_{\lambda}\bigl(\theta_{n+1}^{(n+1)} - t^{n-|\lambda|}\bigr) = 0$}.
 \end{itemize}
\end{Proposition}

\begin{proof}
 From Lemma~\ref{finite hecke map} and Definition~\ref{connecting maps def}, it follows immediately for all $1\leq i \leq n-1$ and~${1\leq j \leq n}$, that \smash{$\Phi^{(n)}_{\lambda}T_i = T_i\Phi^{(n)}_{\lambda}$}, \smash{$\Phi^{(n)}_{\lambda}X_j = X_j\Phi^{(n)}_{\lambda}$}, and \smash{$\Phi^{(n)}_{\lambda}X_{n+1} = 0$}.

 Let \smash{$X^{\alpha}\otimes v \in V_{\lambda^{(n+1)}}$}.
 By direct calculation, we find
 \begin{gather*}
 \Phi^{(n)}_{\lambda}t^{-n}\pi_{n+1}T_n\bigl(X_1^{\alpha_1}\cdots X_{n+1}^{\alpha_{n+1}}\otimes v\bigr) \\
 \qquad = \Phi^{(n)}_{\lambda}t^{-n}\pi_{n+1} X_1^{\alpha_1}\cdots X_{n-1}^{\alpha_{n-1}}T_n\bigl(X_n^{\alpha_n}X_{n+1}^{\alpha_{n+1}}\otimes v\bigr) \\
 \qquad = \Phi^{(n)}_{\lambda}t^{-n}X_2^{\alpha_1}\cdots X_n^{\alpha_{n-1}}\pi_{n+1}T_n\bigl(X_n^{\alpha_n}X_{n+1}^{\alpha_{n+1}} \otimes v\bigr) \\
 \qquad = t^{-n}X_2^{\alpha_1}\cdots X_n^{\alpha_{n-1}} \Phi^{(n)}_{\lambda} \pi_{n+1} T_n\bigl(X_n^{\alpha_n}X_{n+1}^{\alpha_{n+1}} \otimes v\bigr) \\
 \qquad = t^{-n}X_2^{\alpha_1}\cdots X_n^{\alpha_{n-1}} \\
 \phantom{ \qquad =}{} \times \Phi^{(n)}_{\lambda} \pi_{n+1} \left(X_{n+1}^{\alpha_n}X_n^{\alpha_{n+1}}T_n\otimes v + (1-t)X_n \frac{X_n^{\alpha_n}X_{n+1}^{\alpha_{n+1}}-X_{n+1}^{\alpha_n}X_{n}^{\alpha_{n+1}}}{X_n -X_{n+1}}\otimes v \right) \\
 \qquad = t^{-n}X_2^{\alpha_1}\cdots X_n^{\alpha_{n-1}} \Phi^{(n)}_{\lambda} \left( q^{\alpha_n}X_1^{\alpha_n}X_{n+1}^{\alpha_{n+1}} \pi_{n+1}T_n \otimes v \right. \\
\phantom{ \qquad =}{} \left. + (1-t)X_{n+1}\pi_{n+1}\frac{X_n^{\alpha_n}X_{n+1}^{\alpha_{n+1}}-X_{n+1}^{\alpha_n}X_{n}^{\alpha_{n+1}}}{X_n -X_{n+1}} \otimes v \right) \\
 \qquad = \mathbbm{1}(\alpha_{n+1}=0)t^{-n}q^{\alpha_n}X_1^{\alpha_n}X_2^{\alpha_1}\cdots X_n^{\alpha_{n-1}}\Phi^{(n)}_{\lambda}(1 \otimes \rho_{n+1}(\pi_{n+1}T_n)v) \\
 \qquad = \mathbbm{1}(\alpha_{n+1}=0)t^{-n}q^{\alpha_n}X_1^{\alpha_n}X_2^{\alpha_1}\cdots X_n^{\alpha_{n-1}}\Phi^{(n)}_{\lambda}\bigl(1 \otimes t^nT_1^{-1}\cdots T_{n-1}^{-1}v\bigr) \\
 \qquad = \mathbbm{1}(\alpha_{n+1}=0)q^{\alpha_n}X_1^{\alpha_n}X_2^{\alpha_1}\cdots X_n^{\alpha_{n-1}} \otimes T_1^{-1}\cdots T_{n-1}^{-1}\mathfrak{q}_{\lambda}^{(n)}(v).
\end{gather*}

On the other hand, we see
\begin{gather*}
 t^{-(n-1)}\pi_n \Phi^{(n)}_{\lambda}\bigl(X_1^{\alpha}\cdots X_{n+1}^{\alpha_{n+1}}\otimes v\bigr) \\
 \qquad =\mathbbm{1}(\alpha_{n+1}=0) t^{-(n-1)}\pi_n\bigl(X_1^{\alpha_1}\cdots X_n^{\alpha_n}\bigr) \otimes \mathfrak{q}_{\lambda}^{(n)}(v) \\
 \qquad = \mathbbm{1}(\alpha_{n+1} =0) t^{-(n-1)}q^{\alpha_n}X_1^{\alpha_n}X_2^{\alpha_1}\cdots X_n^{\alpha_{n-1}}\otimes \rho_n(\pi_n)\bigl(\mathfrak{q}_{\lambda}^{(n)}(v)\bigr) \\
 \qquad = \mathbbm{1}(\alpha_{n+1} =0)q^{\alpha_n}X_1^{\alpha_n}X_2^{\alpha_1}\cdots X_n^{\alpha_{n-1}}\otimes T_1^{-1}\cdots T_{n-1}^{-1}\mathfrak{q}_{\lambda}^{(n)}(v).
\end{gather*}

Therefore, \smash{$\Phi^{(n)}_{\lambda} t^{-n}\pi_{n+1}T_n = t^{-(n-1)}\pi_n \Phi^{(n)}_{\lambda}$} as desired.

Now, let $1\leq i \leq n$. We see that
\begin{align*}
 \Phi^{(n)}_{\lambda} \theta_i^{(n+1)}&= \Phi^{(n)}_{\lambda} t^{-(n-i+1)}T_{i-1}^{-1}\cdots T_1^{-1} \pi_{n+1} T_n\cdots T_i \\
 & = t^{i-1}T_{i-1}^{-1}\cdots T_1^{-1}\bigl(\Phi^{(n)}_{\lambda}t^{-n}\pi_{n+1} T_n\bigr) T_{n-1}\cdots T_i \\
 & = t^{i-1}T_{i-1}^{-1}\cdots T_1^{-1} \bigl(t^{-(n-1)}\pi_n \Phi^{(n)}_{\lambda}\bigr)T_{n-1}\cdots T_i \\
 & = t^{-(n-i)} T_{i-1}^{-1}\cdots T_1^{-1} \pi_n T_{n-1}\cdots T_i\Phi^{(n)}_{\lambda}
 = \theta_i^{(n)}\Phi^{(n)}_{\lambda}.
\end{align*}

Now, let \smash{$\alpha \in \mathbb{Z}^{n+1}_{\geq 0}$} and $\tau \in \SYT\bigl(\lambda^{(n+1)}\bigr)$. We find
\begin{align*}
 \Phi^{(n)}_{\lambda} \theta_{n+1}^{(n+1)}(X^{\alpha} \otimes e_{\tau}) &
 = \Phi^{(n)}_{\lambda}T_n^{-1}\cdots T_1^{-1} \pi_{n+1}(X^{\alpha} \otimes e_{\tau}) \\
 & = \Phi^{(n)}_{\lambda}T_n^{-1}\cdots T_1^{-1} q^{\alpha_{n+1}}X_1^{\alpha_{n+1}}X_2^{\alpha_1}\cdots X_{n+1}^{\alpha_n} \otimes \rho_{n+1}(\pi_{n+1})e_{\tau} \\
 & = q^{\alpha_{n+1}}\Phi^{(n)}_{\lambda}T_n^{-1}\cdots T_1^{-1}X_1^{\alpha_{n+1}}X_2^{\alpha_1}\cdots X_{n+1}^{\alpha_n} \bigl(1\otimes t^nT_1^{-1}\cdots T_n^{-1}(e_{\tau})\bigr).
\end{align*}
If $\alpha_{n+1} >0$, then this evaluates to $0$ since
\[
\Phi^{(n)}_{\lambda}T_n^{-1}\cdots T_1^{-1}X_1 = t^{-n}\Phi^{(n)}_{\lambda} X_{n+1}T_n \cdots T_1 =0.
\]
Hence,
\[
 \Phi^{(n)}_{\lambda} \theta_{n+1}^{(n+1)}(X^{\alpha} \otimes e_{\tau}) = \mathbbm{1}(\alpha_{n+1} =0)\Phi^{(n)}_{\lambda}T_n^{-1}\cdots T_1^{-1}X_2^{\alpha_1}\cdots X_{n+1}^{\alpha_n} \bigl(1\otimes t^nT_1^{-1}\cdots T_n^{-1}(e_{\tau})\bigr).
\]

By repeatedly applying Lemma~\ref{action of T inverse}, we see that as maps $V_{\lambda^{(n+1)}} \rightarrow V_{\lambda^{(n)}}$
\begin{gather*}
 \Phi^{(n)}_{\lambda}T_n^{-1}\cdots T_2^{-1} \bigl(T_{1}^{-1}X_2^{\alpha_1}\bigr) X_3^{\alpha_2}\cdots X_{n+1}^{\alpha_n} \\
 \qquad = \Phi^{(n)}_{\lambda}T_n^{-1}\cdots T_2^{-1}\left(X_{1}^{\alpha_1}T_1^{-1} + \bigl(1-t^{-1}\bigr)X_2 \frac{X_1^{\alpha_1}-X_2^{\alpha_1}}{X_1-X_2}\right)X_3^{\alpha_2}\cdots X_{n+1}^{\alpha_n} \\
 \qquad =X_{1}^{\alpha_1}\Phi^{(n)}_{\lambda}T_n^{-1}\cdots T_1^{-1}X_3^{\alpha_2}\cdots X_{n+1}^{\alpha_n} \\
 \phantom{\qquad =}{} + \bigl(1-t^{-1}\bigr)\Phi^{(n)}_{\lambda}T_n^{-1}\cdots T_2^{-1}X_2 \frac{X_1^{\alpha_1}-X_2^{\alpha_1}}{X_1-X_2}X_3^{\alpha_2}\cdots X_{n+1}^{\alpha_n} \\
 \qquad = X_{1}^{\alpha_1}\Phi^{(n)}_{\lambda}T_n^{-1}\cdots T_1^{-1}X_3^{\alpha_2}\cdots X_{n+1}^{\alpha_n} \\
 \phantom{\qquad =}{} + \bigl(1-t^{-1}\bigr)t^{-(n-2)}\Phi^{(n)}_{\lambda} X_{n+1} T_{n-1}\cdots T_2 \frac{X_1^{\alpha_1}-X_2^{\alpha_1}}{X_1-X_2}X_3^{\alpha_2}\cdots X_{n+1}^{\alpha_n} \\
 \qquad = X_{1}^{\alpha_1}\Phi^{(n)}_{\lambda}T_n^{-1}\cdots T_1^{-1}X_3^{\alpha_2}\cdots X_{n+1}^{\alpha_n} + 0 \\
 \qquad = X_{1}^{\alpha_1}\Phi^{(n)}_{\lambda}T_n^{-1}\cdots T_3^{-1}(T_2^{-1} X_3^{\alpha_2})T_1^{-1}X_4^{\alpha_3}\cdots X_{n+1}^{\alpha_n} \\
 \qquad = \cdots \\
 \qquad = X_{1}^{\alpha_1}X_2^{\alpha_2}\Phi^{(n)}_{\lambda}T_n^{-1}\cdots T_1^{-1}X_4^{\alpha_3}\cdots X_{n+1}^{\alpha_n} \\
 \qquad = \cdots \\
 \qquad = X_1^{\alpha_1}\cdots X_{n}^{\alpha_{n}}\Phi^{(n)}_{\lambda}T_n^{-1}\cdots T_1^{-1}.
\end{gather*}

As usual, let $\square_0$ denote the unique square of the skew diagram $\lambda^{(n+1)}/\lambda^{(n)}$. Returning to our main calculation now shows
\begin{gather*}
 \Phi^{(n)}_{\lambda} \theta_{n+1}^{(n+1)}(X^{\alpha} \otimes e_{\tau}) \\
 \qquad =\mathbbm{1}(\alpha_{n+1} =0)\Phi^{(n)}_{\lambda}T_n^{-1}\cdots T_1^{-1}X_2^{\alpha_1}\cdots X_{n+1}^{\alpha_n} \bigl(1\otimes t^nT_1^{-1}\cdots T_n^{-1}(e_{\tau})\bigr) \\
 \qquad = \mathbbm{1}(\alpha_{n+1} =0)X_1^{\alpha_1}\cdots X_{n}^{\alpha_{n}}\Phi^{(n)}_{\lambda}T_n^{-1}\cdots T_1^{-1} \bigl(1\otimes t^nT_1^{-1}\cdots T_n^{-1}(e_{\tau})\bigr) \\
 \qquad = \mathbbm{1}(\alpha_{n+1} =0)X_1^{\alpha_1}\cdots X_{n}^{\alpha_{n}}\Phi^{(n)}_{\lambda}\bigl(1 \otimes t^nT_n^{-1}\cdots T_1^{-1}T_1^{-1}\cdots T_n^{-1}(e_{\tau})\bigr) \\
 \qquad = \mathbbm{1}(\alpha_{n+1} =0)X_1^{\alpha_1}\cdots X_{n}^{\alpha_{n}}\Phi^{(n)}_{\lambda}\bigl(1 \otimes \overline{\theta}_{n+1}^{(n+1)}(e_{\tau})\bigr) \\
 \qquad = \mathbbm{1}(\alpha_{n+1} =0)X_1^{\alpha_1}\cdots X_{n}^{\alpha_{n}}\Phi^{(n)}_{\lambda}\bigl(1 \otimes t^{c_{\tau}(n+1)}e_{\tau}\bigr) \\
 \qquad = t^{c_{\tau}(n+1)}\mathbbm{1}(\alpha_{n+1} =0)X_1^{\alpha_1}\cdots X_{n}^{\alpha_{n}}\otimes \mathfrak{q}_{\lambda}^{(n)}(e_{\tau}) \\
 \qquad = t^{c(\square_0)}\mathbbm{1}(\alpha_{n+1} =0)\mathbbm{1}( \tau(\square_0) = n+1)X_1^{\alpha_1}\cdots X_{n}^{\alpha_{n}}\otimes e_{\tau|{\lambda^{(n)}}} \\
 \qquad = t^{n-|\lambda|} \mathbbm{1}(\alpha_{n+1} =0)\mathbbm{1}( \tau(\square_0) = n+1)X_1^{\alpha_1}\cdots X_{n}^{\alpha_{n}}\otimes e_{\tau|{\lambda^{(n)}}} = \Phi^{(n)}_{\lambda}\bigl(t^{n-|\lambda|}X^{\alpha}\otimes e_{\tau}\bigr).
\end{gather*}

Therefore,
\smash{$
\Phi^{(n)}_{\lambda}\bigl( \theta_{n+1}^{(n+1)} - t^{n-|\lambda|}\bigr) = 0$}.
\end{proof}

Here we utilize the exceptional properties in Proposition~\ref{relations for connecting maps} to show that the $F_{\tau}$ behave well with respect to the connecting maps \smash{$\Phi^{(n)}_{\lambda}$} from Definition~\ref{connecting maps def}.

\begin{Proposition}\label{stability for locally maximal psyt}
 Let $n \geq n_{\lambda}$ and \smash{$\square_0 = \lambda^{(n+1)}/\lambda^{(n)}$}. For \smash{$\tau \in \PSYT\bigl(\lambda^{(n+1)}\bigr)$}, we have
 \[
 \Phi^{(n)}_{\lambda}(F_{\tau}) := \begin{cases}
 F_{\tau|_{\lambda^{(n)}}}, & \tau(\square_0) = n+1,\\
 0, & \tau(\square_0) \neq n+1.
 \end{cases}
 \]

\end{Proposition}

\begin{proof}
 First, we deal with the case when $\tau(\square_0) = n+1$. Let $T \in \RYT\bigl(\lambda^{(n)}\bigr)$ and let \smash{$T' \in \RYT\bigl(\lambda^{(n+1)}\bigr)$} with $T'(\square_0) = 0$ and \smash{$T'|_{\lambda^{(n)}}=T|_{\lambda^{(n)}}$}. By looking at the eigenvalues of \smash{$\theta_1^{(n+1)},\dots, \theta_n^{(n+1)}$} on $F_{\Top(T')}$ and the eigenvalues of \smash{$\theta_1^{(n)},\dots, \theta_n^{(n)}$} on $F_{\Top(T)}$, we see that \smash{$\Phi^{(n)}_{\lambda}(F_{\Top(T')}) = \beta F_{\Top(T)}$} for some scalar $\beta$. We now show that $\beta = 1$. From Corollary~\ref{triangular expansion of top element}, we know that
 \[
 F_{\Top(T')} = t^{-b_{T'}}X^{\nu(T')}\otimes e_{S(T')} + \sum_{\beta \prec \nu(T')} X^{\beta} \otimes v'_{\beta}
 \]
 and
 \[
 F_{\Top(T)} = t^{-b_{T}}X^{\nu(T)}\otimes e_{S(T)} + \sum_{\beta \prec \nu(T)} X^{\beta} \otimes v_{\beta}
 \]
 for some \smash{$v_{\beta} \in S_{\lambda^{(n)}}$} and \smash{$v'_{\beta} \in S_{\lambda^{(n+1)}}$}. Since $T'(\square_0) = 0$ and \smash{$T'|_{\lambda^{(n)}}=T|_{\lambda^{(n)}}$}, it follows that $b_{T'} = b_{T}$, $\nu(T') = \nu(T)*0$, and \smash{$\mathfrak{q}_{\lambda}^{(n)}(e_{S(T')})= e_{S(T)}$}. Therefore,
 \[
 \Phi^{(n)}_{\lambda}\bigl(t^{-b_{T'}}X^{\nu(T')}\otimes e_{S(T')}\bigr) = t^{-b_{T}}X^{\nu(T)}\otimes e_{S(T)}.
 \]
 Now if $\beta \prec \nu(T')$, then
 \smash{$
 \Phi^{(n)}_{\lambda}\bigl(X^{\beta}\otimes v'_{\beta}\bigr) = \mathbbm{1}(\beta_{n+1}=0)X_1^{\beta_1}\cdots X_n^{\beta_n} \otimes \mathfrak{q}_{\lambda}^{(n)}(v'_{\beta})$} cannot be of the form \smash{$X^{\nu(T)}\otimes w$} for any \smash{$w \in S_{\lambda^{(n)}}$}. As such, the coefficient of \smash{$X^{\nu(T)}\otimes e_{S(T)}$} in the standard basis expansion of \smash{$\Phi^{(n)}_{\lambda}(F_{\Top(T')})$} is $t^{-b_T}$. Since this agrees with the same coefficient in the expansion of $F_{\Top(T)}$, we know that $\beta = 1$ and thus \smash{$\Phi^{(n)}_{\lambda}(F_{\Top(T')}) =F_{\Top(T)}$}.

 Consider any \smash{$\tau' \in \PSYT\bigl(\lambda^{(n+1)}\bigr)$} with $\tau'(\square_0) = n+1$. Let
 \[
 T' := \mathfrak{p}_{\lambda^{(n+1)}}(\tau) \in \RYT\bigl(\lambda^{(n+1)}\bigr).
 \]
 Then $T'(\square_0) = 0$ so if we set \smash{$T:= T'|_{\lambda^{(n)}}$}, we have that \smash{$\Phi^{(n)}_{\lambda}(F_{\Top(T')}) = F_{\Top(T)}$}. Write $\tau:= \smash{ \tau'|_{\lambda^{(n)}}}$. As seen before, there exists a~sequence $\tau < s_{i_1}(\tau) < \dots< s_{i_r}\cdots s_{i_1}(\tau) = \Top(T)$. Since $\tau'(\square_0) = n+1$, we see that $\tau' < s_{i_1}(\tau') < \dots< s_{i_r}\cdots s_{i_1}(\tau') = \Top(T')$ as well. For each $1\leq j \leq r$, we consider using the intertwiner operators from Proposition~\ref{weight basis Proposition} to obtain~\smash{$F_{s_{i_j}s_{i_{j-1}}\cdots s_{i_1}(\tau)}$} from \smash{$F_{s_{i_{j-1}}\cdots s_{i_1}(\tau)}$}. We have that
 \begin{gather*}
 \left( tT_{i_j}^{-1} + \frac{(t-1)q^{w_{s_{i_{j-1}}\cdots s_{i_1}(\tau)}(i_j+1)}t^{c_{s_{i_{j-1}}\cdots s_{i_1}(\tau)}(i_j+1)}}{q^{w_{s_{i_{j-1}}\cdots s_{i_1}(\tau)}(i_j)}t^{c_{s_{i_{j-1}}\cdots s_{i_1}(\tau)}(i_j)}-q^{w_{s_{i_{j-1}}\cdots s_{i_1}(\tau)}(i_j+1)}t^{c_{s_{i_{j-1}}\cdots s_{i_1}(\tau)}(i_j+1)}}\right)\\ \qquad\times(F_{s_{i_{j-1}}\cdots s_{i_1}(\tau)})
 = F_{s_{i_j}(s_{i_{j-1}}\cdots s_{i_1}(\tau))}.
 \end{gather*}

 Now the same exact formula holds with $\tau$ replaced by $\tau'$. Importantly, we have that
 \begin{gather*}
 w_{s_{i_{j-1}}\cdots s_{i_1}(\tau)}(i_j+1) = w_{s_{i_{j-1}}\cdots s_{i_1}(\tau')}(i_j+1),\\ c_{s_{i_{j-1}}\cdots s_{i_1}(\tau)}(i_j+1) = c_{s_{i_{j-1}}\cdots s_{i_1}(\tau')}(i_j+1).\end{gather*}
 Therefore, we may write
 \begin{gather*}
 D_j(F_{s_{i_{j-1}}\cdots s_{i_1}(\tau)}) = F_{s_{i_j}s_{i_{j-1}}\cdots s_{i_1}(\tau)}
 \qquad \text{and}\qquad
 D_j(F_{s_{i_{j-1}}\cdots s_{i_1}(\tau')}) = F_{s_{i_j}s_{i_{j-1}}\cdots s_{i_1}(\tau')}
\end{gather*}
 for \smash{$D_j \in \sA_n \subset \sA_{(n,1)}$} of the form $D_j = T_{i_j} + \alpha_j$ where $\alpha_j \in \mathbb{Q}(q,t)$. Here we have used \smash{$tT_{i_j}^{-1} = T_{i_j} + t-1$}. By using the quadratic relation for $T_{i_j}$, we may locally invert the operator~$D_j$ in the sense that there exists operators $C_j \in \sA_n$ with
 \begin{gather*}
 F_{s_{i_{j-1}}\cdots s_{i_1}(\tau)} = C_j(F_{s_{i_j}s_{i_{j-1}}\cdots s_{i_1}(\tau)})
 \qquad \text{and}\qquad
 F_{s_{i_{j-1}}\cdots s_{i_1}(\tau')} = C_j(F_{s_{i_j}s_{i_{j-1}}\cdots s_{i_1}(\tau')}).
 \end{gather*}
 Therefore, if we assume that \smash{$\Phi^{(n)}_{\lambda}(F_{s_{i_j}s_{i_{j-1}}\cdots s_{i_1}(\tau')})= F_{s_{i_j}s_{i_{j-1}}\cdots s_{i_1}(\tau)}$}, then
 \begin{align*}
 \Phi^{(n)}_{\lambda}(F_{s_{i_{j-1}}\cdots s_{i_1}(\tau')})&= \Phi^{(n)}_{\lambda}(C_j(F_{s_{i_j}s_{i_{j-1}}\cdots s_{i_1}(\tau')}))
 = C_j\Phi^{(n)}_{\lambda}(F_{s_{i_j}s_{i_{j-1}}\cdots s_{i_1}(\tau')}) \\
 & = C_jF_{s_{i_j}s_{i_{j-1}}\cdots s_{i_1}(\tau)}
 = F_{s_{i_{j-1}}\cdots s_{i_1}(\tau)}.
 \end{align*}

Thus by induction, since we know \smash{$\Phi^{(n)}_{\lambda}(F_{\Top(T')}) = F_{\Top(T)}$}, it follows that \smash{$\Phi^{(n)}_{\lambda}(F_{\tau'}) = F_{\tau}$}.

 Lastly, we consider the case of $\tau(\square_0) \neq n+1$. Then $\tau(\square_0) = iq^a$ with either $i \neq n+1$ or~${a \geq 0}$. If $a >0$, then $\tau = \Psi(\tau')$ for some $\tau'$ and thus from Proposition~\ref{weight basis Proposition}, we know $X_{n+1}$ divides $F_{\tau}$. Since $\Phi^{(n)}_{\lambda}X_{n+1} = 0$ it follows that $\Phi^{(n)}_{\lambda}(F_{\tau}) =0$. Suppose $a = 0$ and $i \neq n+1$. Notice for any $m \geq n_{\lambda}$ that the largest power of $t$ occurring in the \smash{$\Theta^{(m)}$}-weight of any $F_{\tau'}$ with~${\tau' \in \PSYT(\lambda^{(m)})}$ is exactly $t^{m-|\lambda|-1}$. Since $i \neq n+1$, we know that if \smash{$\Phi^{(n)}_{\lambda}(F_{\tau}) = \beta F_{\tau'}$} for some nonzero scalar $\beta$ and $\tau' \in \PSYT\bigl(\lambda^{(n+1)}\bigr)$, then the maximal power of $t$ occurring in the $\Theta^{(n)}$-weight of $F_{\tau'}$ is $t^{n-|\lambda|}$ coming from
 \[
 \theta_i(F_{\tau'}) = \theta_i\bigl( \Phi^{(n)}_{\lambda}(F_{\tau})\bigr)= \Phi^{(n)}_{\lambda}(\theta_i(F_{\tau})))= \Phi^{(n)}_{\lambda}\bigl( t^{c(\square_0)}F_{\tau}\bigr)= \Phi^{(n)}_{\lambda}\bigl( t^{n-|\lambda|}F_{\tau}\bigr)= t^{n-|\lambda|} F_{\tau'}.
 \]
 Thus, \smash{$\Phi^{(n)}_{\lambda}(F_{\tau})$} cannot be a~$\Theta^{(n)}$-weight vector in $V_{\lambda^{(n)}}$ and so $\beta=0$.
\end{proof}

The maps \smash{$\Phi^{(n)}_{\lambda}$} possess another important stability property.

\begin{Proposition}\label{connecting maps and symmetrization of thetas}
 For all $\ell \in \mathbb{Z}\setminus \{0\}$ and $n \geq n_{\lambda}$,
 \[
 \Phi^{(n)}_{\lambda}\left( \sum_{j=1}^{n+1}\bigl(\theta^{(n+1)}_j\bigr)^{\ell} - \sum_{\square \in \lambda^{(n+1)}} t^{\ell c(\square)} \right) = \left( \sum_{j=1}^{n}\bigl(\theta^{(n)}_j\bigr)^{\ell} - \sum_{\square \in \lambda^{(n)}} t^{\ell c(\square)} \right) \Phi^{(n)}_{\lambda}.
 \]

\end{Proposition}

\begin{proof}
 Let $\ell \in \mathbb{Z}\setminus \{0\}$ and $n \geq n_{\lambda}$. As usual, let $\square_0$ denote the unique square of the skew diagram $\lambda^{(n+1)}/\lambda^{(n)}$. Directly from Proposition~\ref{relations for connecting maps}, we see
 \begin{gather*}
 \Phi^{(n)}_{\lambda}\left( \sum_{j=1}^{n+1}\bigl(\theta^{(n+1)}_j\bigr)^{\ell} - \sum_{\square \in \lambda^{(n+1)}} t^{\ell c(\square)} \right) \\
 \qquad = \Phi^{(n)}_{\lambda}\left( \sum_{j=1}^{n}\bigl(\theta^{(n+1)}_j\bigr)^{\ell} - \sum_{\square \in \lambda^{(n)}} t^{\ell c(\square)} \right) + \Phi^{(n)}_{\lambda}\bigl(\bigl(\theta^{(n+1)}_{n+1}\bigr)^{\ell} - t^{\ell c(\square_0)}\bigr) \\
 \qquad = \left( \sum_{j=1}^{n}\bigl(\theta^{(n)}_j\bigr)^{\ell} - \sum_{\square \in \lambda^{(n)}} t^{\ell c(\square)} \right)\Phi^{(n)}_{\lambda} + \Phi^{(n)}_{\lambda}\bigl(\bigl(\theta^{(n+1)}_{n+1}\bigr)^{\ell} - t^{\ell (n-|\lambda|)}\bigr).
 \end{gather*}
 It follows from the relation \smash{$\Phi^{(n)}_{\lambda}\bigl(\theta^{(n+1)}_{n+1} - t^{n-|\lambda|}\bigr) = 0$} and the fact that \smash{$\theta^{(n+1)}_{n+1}$} is invertible on~\smash{$V_{\lambda^{(n+1)}}$} that
 \smash{$
 \Phi^{(n)}_{\lambda}\bigl(\bigl(\theta^{(n+1)}_{n+1}\bigr)^{\ell} - t^{\ell (n-|\lambda|)}\bigr) = 0$}.
 Therefore,
 \[
 \Phi^{(n)}_{\lambda}\left( \sum_{j=1}^{n+1}\bigl(\theta^{(n+1)}_j\bigr)^{\ell} - \sum_{\square \in \lambda^{(n+1)}} t^{\ell c(\square)} \right) =\left( \sum_{j=1}^{n}\bigl(\theta^{(n)}_j\bigr)^{\ell} - \sum_{\square \in \lambda^{(n)}} t^{\ell c(\square)} \right)\Phi^{(n)}_{\lambda}.\tag*{\qed}
 \]\renewcommand{\qed}{}
\end{proof}

\section{Positive EHA representations from Young diagrams}\label{Positive EHA Representations from Young Diagrams}

In this section, we construct the Murnaghan-type representations \smash{$\widetilde{W}_{\lambda}$} of $\mathcal{E}^+$ along with specified bases $\MacD_{T}$ which generalize the symmetric Macdonald functions. The representations \smash{$\widetilde{W}_{\lambda}$} are shown to be mutually non-isomorphic and for $\lambda = \varnothing$ the representation $\widetilde{W}_{\varnothing}$ may be identified with the representation of $\mathcal{E}^+$ on the space of symmetric functions $\Lambda_{q,t}$. The action of $P_{0,1} \in \mathcal{E}^+$ gives an analogue of the Macdonald operator $\Delta$ for each of the spaces \smash{$\widetilde{W}_{\lambda}$}. We give explicit formulas for the spectrum of $\Delta$ on each \smash{$\widetilde{W}_{\lambda}$} and recover the known formula in the $\lambda = \varnothing$ case. To accomplish these results, we consider the spherical modules $W_{\lambda^{(n)}}$ corresponding to the $\sD_n$-modules $V_{\lambda^{(n)}}$ as $n$ varies. We show, after restricting the connecting maps \smash{$\Phi_{\lambda}^{(n)}$} to the $W_{\lambda^{(n)}}$, they are equivariant with respect to the $\sD_n^{\mathrm{sph}}$ actions in the sense of Definition~\ref{EHA def}. At the same time, we construct the bases $\MacD_{T}$ by further studying the combinatorics underpinning the $\sD_n$-modules $V_{\lambda^{(n)}}$.

\subsection[The D\_n\^{}sph-modules W\_lambda\^{}(n)]{The $\boldsymbol{\sD_n^{\mathrm{sph}}}$-modules $\boldsymbol{W_{\lambda^{(n)}}}$}

We now consider the spherical DAHA modules and symmetric vv polynomials corresponding to the positive DAHA modules $V_{\lambda}$ and the non-symmetric vv polynomials $F_{\tau}$ considered in the prior sections.
 \begin{Definition}\label{sph DAHA reps Definition}
 For $\lambda \in \Par$ with $|\lambda| = n$, define the $\sD_n^{\mathrm{sph}}$-module $W_{\lambda}:= \epsilon^{(n)}(V_{\lambda})$.
 \end{Definition}

 The $F_{\tau}$ expansions of any symmetrized element of any $\sA_n$-submodule $U_{T}$ satisfy a~simple set of recurrence relations.

\begin{Lemma}\label{coefficient relations for sym element}
 Let $T \in \RSSYT(\lambda)$ and $v \in \epsilon^{(n)}(U_T)$. Suppose that $v$ has the following expansion into the $F_{\tau}$ basis
 \smash{$
 v = \sum_{\tau \in \PSYT(\lambda;T)} \kappa_{\tau}F_{\tau}$}.

 Then for each $\tau \in \PSYT(\lambda;T)$ with $1\leq i \leq n-1$ such that $s_i(\tau) > \tau$, we have the relation
 \[
 \kappa_{s_i(\tau)} = \left( \frac{q^{w_{\tau}(i)}t^{c_{\tau}(i)} - q^{w_{\tau}(i+1)}t^{c_{\tau}(i+1)}}{q^{w_{\tau}(i)}t^{c_{\tau}(i)} - q^{w_{\tau}(i+1)}t^{c_{\tau}(i+1)+1}} \right) \kappa_{\tau}.
 \]
 As a~consequence, if $\kappa_{\Top(T)} \neq 0$, then each coefficient $\kappa_{\tau}$ is also nonzero.
\end{Lemma}

\begin{proof}
 Let $\tau \in \PSYT(\lambda;T)$ and $1 \leq i \leq n-1$ with $s_i(\tau) > \tau$. Note that $\mathbb{Q}(q,t)\{F_{\tau}, F_{s_i(\tau)}\}$ is a~2-dimensional submodule for $\mathbb{Q}(q,t)[T_i]$. The $T_i$-invariant subspace of $\mathbb{Q}(q,t)\{F_{\tau}, F_{s_i(\tau)}\}$ is given by \smash{$\mathbb{Q}(q,t)\bigl(1+tT_i^{-1}\bigr)F_{\tau}$}. From Proposition~\ref{weight basis Proposition}, we find
 \begin{align*}
 \bigl(1+tT_{i}^{-1}\bigr)F_{\tau}& = F_{\tau} + tT_i^{-1}F_{\tau}
 = F_{\tau} + F_{s_i(\tau)} + \frac{(1-t)q^{w_{\tau}(i+1)}t^{c_{\tau}(i+1)}}{q^{w_{\tau}(i)}t^{c_{\tau}(i)}-q^{w_{\tau}(i+1)}t^{c_{\tau}(i+1)}} F_{\tau} \\
 & = F_{s_i(\tau)} + \frac{q^{w_{\tau}(i)}t^{c_{\tau}(i)}-q^{w_{\tau}(i+1)}t^{c_{\tau}(i+1)+1}}{q^{w_{\tau}(i)}t^{c_{\tau}(i)}-q^{w_{\tau}(i+1)}t^{c_{\tau}(i+1)}}F_{\tau}.
 \end{align*}
 Since \smash{$v = \sum_{\tau \in \PSYT(\lambda;T)} \kappa_{\tau}F_{\tau}$} is $T_i$-invariant then we know that in particular
 $\kappa_{\tau}F_{\tau}\allowbreak + \smash{\kappa_{s_i(\tau)}F_{s_i(\tau)}}$ is also $T_i$-invariant and therefore must be a~scalar multiple of $\bigl(1+tT_i^{-1}\bigr)F_{\tau}$. Therefore,
 \[
 \kappa_{\tau}F_{\tau} + \kappa_{s_i(\tau)}F_{s_i(\tau)} = \kappa_{s_i(\tau)}F_{s_i(\tau)} + \frac{q^{w_{\tau}(i)}t^{c_{\tau}(i)}-q^{w_{\tau}(i+1)}t^{c_{\tau}(i+1)+1}}{q^{w_{\tau}(i)}t^{c_{\tau}(i)}-q^{w_{\tau}(i+1)}t^{c_{\tau}(i+1)}}\kappa_{s_i(\tau)}F_{\tau}
 \]
 and so
 \[
 \kappa_{s_i(\tau)} = \left( \frac{q^{w_{\tau}(i)}t^{c_{\tau}(i)} - q^{w_{\tau}(i+1)}t^{c_{\tau}(i+1)}}{q^{w_{\tau}(i)}t^{c_{\tau}(i)} - q^{w_{\tau}(i+1)}t^{c_{\tau}(i+1)+1}} \right) \kappa_{\tau}.\tag*{\qed}
 \]\renewcommand{\qed}{}
\end{proof}

Using the recurrence relations in Lemma~\ref{coefficient relations for sym element} and the irreducibility of each of the $\sA_n$-sub\-modules $U_{T}$ of $V_{\lambda}$, we may determine which $T \in \RYT(\lambda)$ have a~non-zero space of $\sH_n$\nobreakdash-invariants $\epsilon^{(n)}(U_T)$.

\begin{Proposition}\label{sym AHA submodules}
 For $\lambda \in \Par$ with $|\lambda| = n$ and $T \in \RYT(\lambda)$,
 \[
 \mathrm{dim}_{\mathbb{Q}(q,t)} \epsilon^{(n)}(U_T) =\begin{cases}
 1 ,& T \in \RSSYT(\lambda),\\
 0, & T \notin \RSSYT(\lambda).
 \end{cases}
 \]

 \end{Proposition}

 \begin{proof}
 By Proposition~\ref{Mackey decomposition}, each $\sA_n$-module $U_T$ is irreducible with simple $\Theta^{(n)}$ spectrum. This implies that \smash{$\mathrm{dim}_{\mathbb{Q}(q,t)} \epsilon^{(n)}(U_T) \leq 1$} for any $T \in \RYT(\lambda)$. Further, we have that $\epsilon^{(n)}(U_T)$ is zero if for any $\Theta^{(n)}$-weight vector $v$ in $U_T$, $\epsilon^{(n)}(v)$ is zero. If $T \in \RYT(\lambda) \setminus \RSSYT(\lambda)$, then there exists a~pair of boxes $\square_1,\square_2 \in \lambda$ with $\square_1$ directly above $\square_2$ such that $T(\square_1) = T(\square_2) = a$. Hence, $\Top(T)(\square_1) = iq^{a}$ and $\Top(T)(\square_2) = (i+1)q^a$ for some $1\leq i \leq n-1$. Then $T_i(F_{\Top(T)}) = -tF_{\Top(T)}$, which implies that \smash{$\epsilon^{(n)}(F_{\Top(T)}) = 0$}. Thus \smash{$\epsilon^{(n)}(U_T) = 0$}.

 Alternatively, now suppose $T \in \RSSYT(\lambda)$. Following Lemma~\ref{coefficient relations for sym element}, we construct a~vector~${v \in U_{T}}$ of the form
 \[
 v = \sum_{\tau \in \PSYT(\lambda;T)} \kappa_{\tau}F_{\tau},
 \]
 where $\kappa_{\Top(T)} = 1$ and if $s_i(\tau) > \tau$, then
 \[
 \kappa_{s_i(\tau)} = \left( \frac{q^{w_{\tau}(i)}t^{c_{\tau}(i)} - q^{w_{\tau}(i+1)}t^{c_{\tau}(i+1)}}{q^{w_{\tau}(i)}t^{c_{\tau}(i)} - q^{w_{\tau}(i+1)}t^{c_{\tau}(i+1)+1}} \right) \kappa_{\tau}.
 \]
 These coefficients $\kappa_{\tau}$ have the property that if $s_i(\tau) > \tau$, then
 \[
 T_i(\kappa_{\tau}F_{\tau} + \kappa_{s_i(\tau)}F_{s_i(\tau)}) = \kappa_{\tau}F_{\tau} + \kappa_{s_i(\tau)}F_{s_i(\tau)}.
 \]
 By construction, $v \neq 0$ since \smash{$\frac{q^{w_{\tau}(i)}t^{c_{\tau}(i)} - q^{w_{\tau}(i+1)}t^{c_{\tau}(i+1)}}{q^{w_{\tau}(i)}t^{c_{\tau}(i)} - q^{w_{\tau}(i+1)}t^{c_{\tau}(i+1)+1}} \neq 0$} whenever $s_i(\tau) > \tau$. We show that~${T_i(v) = v}$ for all $1\leq i \leq n-1$ and thus \smash{$\epsilon^{(n)}(U_T) \neq 0$}.

 We find that
 \begin{align*}
 T_i(v)
 ={}& \sum_{\tau \in \PSYT(\lambda;T)} \kappa_{\tau} T_i(F_{\tau}) \\
 ={}& \sum_{\substack{(\tau, s_i(\tau))\PSYT(\lambda)^2 \\ s_i(\tau) > \tau}}T_i\left(\kappa_{\tau} (F_{\tau}) + \kappa_{s_i(\tau)} (F_{s_i(\tau)}) \right) \\
 & + \sum_{\substack{\tau \in \PSYT(\lambda) \\ i,i+1 ~\text{same row of}~ \tau }} \kappa_{\tau} T_i(F_{\tau}) + \sum_{\substack{\tau \in \PSYT(\lambda) \\ i,i+1 ~\text{same column of}~ \tau }} \kappa_{\tau} T_i(F_{\tau}) \\
 ={}& \sum_{\substack{(\tau, s_i(\tau))\PSYT(\lambda)^2 \\ s_i(\tau) > \tau}}\left(\kappa_{\tau} (F_{\tau}) + \kappa_{s_i(\tau)} (F_{s_i(\tau)}) \right) \\
 & + \sum_{\substack{\tau \in \PSYT(\lambda) \\ i,i+1 ~\text{same row of}~ \tau }} \kappa_{\tau} F_{\tau} + \sum_{\substack{\tau \in \PSYT(\lambda) \\ i,i+1 ~\text{same column of}~ \tau }} (-t) \kappa_{\tau} F_{\tau}.
 \end{align*}

 Thus,
 \[
 T_i(v) - v = \sum_{\substack{\tau \in \PSYT(\lambda) \\ i,i+1 ~\text{same column of}~ \tau }} (1+t)\kappa_{\tau} F_{\tau}.
 \]

 Lastly, since $T \in \RSSYT(\lambda)$ there cannot be any $\tau \in \PSYT(\lambda;T)$ with $i$, $i+1$ occurring in the same column as necessarily this would imply that $T$ would have redundant values in those boxes contradicting the fact that $T$ is reverse semi-standard. Hence, the above sum vanishes and we find $T_i(v) = v$.
 \end{proof}

Finally, we are able to define the symmetric vv Macdonald polynomials following the conventions of this paper. By combining Lemma~\ref{coefficient relations for sym element} and Proposition~\ref{sym AHA submodules}, we are able for each $T \in \RSSYT(\lambda)$ to uniquely specify a~non-zero vector in $\epsilon^{(n)}(U_T)$.

\begin{Definition}\label{Macdonald polynomial def}
 Let $T \in \RSSYT(\lambda)$. Define \smash{$P_{T} \in \epsilon^{(n)}(U_T)$} to be the unique element of the form
 \smash{$
 P_T = F_{\Top(T)} + \sum_{y} \kappa_{y}F_y$},
 where the sum above ranges over $y \in \PSYT(\lambda)$ with $\mathfrak{p}_{\lambda}(y) = T$ and $y < \Top(T)$.
\end{Definition}

Now we are able to explicitly compute the $F_{\tau}$ expansion of each $P_{T}$ using the recurrence relations found in Lemma~\ref{coefficient relations for sym element}.

\begin{Corollary}\label{expansion of sym into nonsym}
 For all $T\in \RSSYT(\lambda)$,
 \[
 P_T = \sum_{\tau \in \PSYT(\lambda;T)} \prod_{(\square_1,\square_2) \in \Inv(\tau)} \left( \frac{q^{T(\square_1)}t^{c(\square_1) +1} - q^{T(\square_2)}t^{c(\square_2) }}{q^{T(\square_1)}t^{c(\square_1)} - q^{T(\square_2)}t^{c(\square_2)}} \right) F_{\tau}.
 \]

\end{Corollary}
\begin{proof}
 For $\tau \in \PSYT(\lambda;T)$, let
 \[
 \kappa_{\tau}= \prod_{(\square_1,\square_2) \in \Inv(\tau)} \left( \frac{q^{T(\square_1)}t^{c(\square_1) +1} - q^{T(\square_2)}t^{c(\square_2) }}{q^{T(\square_1)}t^{c(\square_1)} - q^{T(\square_2)}t^{c(\square_2)}} \right).
 \]

 {\samepage From Lemma~\ref{coefficient relations for sym element}, it suffices to show that
 \begin{itemize}
\itemsep=0pt
 \item $\kappa_{\Top(T)} = 1$.
 \item If $s_i(\tau) > \tau$, then \smash{$\kappa_{s_i(\tau)} = \bigl( \frac{q^{w_{\tau}(i)}t^{c_{\tau}(i)} - q^{w_{\tau}(i+1)}t^{c_{\tau}(i+1)}}{q^{w_{\tau}(i)}t^{c_{\tau}(i)} - q^{w_{\tau}(i+1)}t^{c_{\tau}(i+1)+1}} \bigr) \kappa_{\tau}$}.
 \end{itemize}}

 It is easy to see that $\Inv(\Top(T)) = \varnothing$ so $\kappa_{\Top(T)} = 1$. Suppose $s_i(\tau) > \tau$. Let $\square^{(i)},\square^{(i+1)} \in \lambda$ denote the boxes of $\lambda$ with \smash{$\tau(\square^{(i)})= i q^a$} and \smash{$\tau(\square^{(i+1)}) = (i+1)q^b$} for some $a,b \geq 0$. It is straightforward to check that
$
 \Inv(s_i(\tau)) = \big\{\bigl(\square^{(i)},\square^{(i+1)}\bigr)\big\} \sqcup \Inv(\tau)$.

 Therefore,
 \begin{align*}
 \kappa_{s_i(\tau)} ={}& \prod_{(\square_1,\square_2) \in \Inv(s_i(\tau))} \left( \frac{q^{T(\square_1)}t^{c(\square_1) +1} - q^{T(\square_2)}t^{c(\square_2) }}{q^{T(\square_1)}t^{c(\square_1)} - q^{T(\square_2)}t^{c(\square_2)}} \right) \\
 = {}&\left( \frac{q^{T(\square^{(i)})}t^{c(\square^{(i)}) +1} - q^{T(\square^{(i+1)})}t^{c(\square^{(i+1)}) }}{q^{T(\square^{(i)})}t^{c(\square^{(i)})} - q^{T(\square^{(i+1)})}t^{c(\square^{(i+1)})}} \right) \\
 &\times\prod_{(\square_1,\square_2) \in \Inv(\tau)} \left( \frac{q^{T(\square_1)}t^{c(\square_1) +1} - q^{T(\square_2)}t^{c(\square_2) }}{q^{T(\square_1)}t^{c(\square_1)} - q^{T(\square_2)}t^{c(\square_2)}} \right) \\
 ={}& \left( \frac{q^{w_{\tau}(i)}t^{c_{\tau}(i)} - q^{w_{\tau}(i+1)}t^{c_{\tau}(i+1)}}{q^{w_{\tau}(i)}t^{c_{\tau}(i)} - q^{w_{\tau}(i+1)}t^{c_{\tau}(i+1)+1}} \right) \kappa_{\tau}.\tag*{\qed}
 \end{align*}\renewcommand{\qed}{}
\end{proof}

We look now at the action of the special spherical DAHA elements $P_{0,\ell}^{(n)}$.

\begin{Proposition}\label{finite rank gen macd function properties}
 Let $|\lambda| = n$. The set $\{P_T\mid T \in \RSSYT(\lambda)\}$ is a~\smash{$\mathbb{Q}(q,t)\big[\theta_1^{\pm 1},\dots, \theta_n^{\pm 1}\big]^{\mathfrak{S}_n}$}-weight basis for $W_{\lambda}$. Further, for $\ell \in \mathbb{Z}\setminus\{0\}$
 \[
 P_{0,\ell}^{(n)}(P_T) = \left( \sum_{\square \in \lambda} q^{\ell T(\square)}t^{\ell c(\square)} \right) P_T.
 \]

 Consequently, $P_{0,1}^{(n)}$ acts on $W_{\lambda}$ with simple spectrum
 \[
 \left\{ \sum_{\square \in \lambda} q^{T(\square)}t^{c(\square)} \mid T \in \RSSYT(\lambda) \right\}.
 \]

\end{Proposition}

\begin{proof}
 It follows directly from Propositions~\ref{Mackey decomposition} and~\ref{sym AHA submodules}, that the set $\{P_T\mid T \in \RSSYT(\lambda)\}$ is a~linear basis for $W_{\lambda}$. We need to show that the $P_T$ are \smash{$\mathbb{Q}(q,t)\big[\theta_1^{\pm 1},\dots, \theta_n^{\pm 1}\big]^{\mathfrak{S}_n}$}-weight vectors.
 Let $T \in \RSSYT(\lambda)$. Then from Definition~\ref{Macdonald polynomial def}, we know that
\smash{$
 P_T = \beta \epsilon^{(n)}(F_{\Top(T)})$}
 for some nonzero scalar $\beta$ depending on $T$. Then for any $\ell \in \mathbb{Z} \setminus\{0\}$, we have that
 \begin{align*}
 \left(\sum_{j=1}^{n} (\theta_j^{(n)})^{\ell}\right)(P_T) &
 = \left(\sum_{j=1}^{n} (\theta_j^{(n)})^{\ell}\right)(\beta \epsilon^{(n)}(F_{\Top(T)})) = \beta \epsilon^{(n)}\left( \left(\sum_{j=1}^{n} (\theta_j^{(n)})^{\ell}\right)(F_{\Top(T)})\right) \\
 & = \beta \epsilon^{(n)}\left( \left(\sum_{j=1}^{n} q^{\ell w_{\Top(T)}(j)}t^{\ell c_{\Top(T)}(j)}\right)F_{\Top(T)} \right) \\
 & = \beta \epsilon^{(n)}\left( \left(\sum_{\square \in \lambda} q^{\ell T(\square)}t^{\ell c(\square)}\right)F_{\Top(T)} \right) \\
 &= \left(\sum_{\square \in \lambda} q^{\ell T(\square)}t^{\ell c(\square)}\right) \beta \epsilon^{(n)}(F_{\Top(T)})
 = \left(\sum_{\square \in \lambda} q^{\ell T(\square)}t^{\ell c(\square)}\right) P_T.
 \end{align*}

 Hence, $P_T$ is a~\smash{$\mathbb{Q}(q,t)\big[\theta_1^{\pm 1},\dots, \theta_n^{\pm 1}\big]^{\mathfrak{S}_n}$}-weight vector. Now let $S \in \RSSYT(\lambda)$ and suppose that
 \[
 \sum_{\square \in \lambda} q^{T(\square)}t^{c(\square)} = \sum_{\square \in \lambda} q^{S(\square)}t^{c(\square)}.
 \]

 Fix any $d \in \mathbb{Z}$.
 Since $q$ and $t$ are algebraically independent over $\mathbb{Q}$,
 \[
 \sum_{\substack{\square \in \lambda \\ c(\square) = d}} q^{T(\square)} = \sum_{\substack{\square \in \lambda \\ c(\square) = d}} q^{S(\square)}.
 \]

 Since the labeling $T$ is reverse semi-standard, the values of $T(\square)$ for $\square \in \lambda$ with $c(\square) = d$ are all distinct and strictly decreasing down the $d$-diagonal. Of course, the same is true for $S$. Therefore, the values of $T$ and $S$ agree along the $d$-diagonal of $\lambda$. As $d \in \mathbb{Z}$ was general, it follows that $T = S$. Thus the spectrum of the operator \smash{$P_{0,1}^{(n)}$} on $W_{\lambda}$ is simple.
\end{proof}

As mentioned previously, the non-symmetric vv Macdonald polynomials do not agree with those of Dunkl--Luque. To see this simply note that the $F_{\tau}$ in this paper are $\Theta^{(n)}$-weight vectors but the non-symmetric vv Macdonald polynomials of Dunkl--Luque are weight vectors for the standard Cherednik elements \smash{$\xi^{(n)}_i := t^{n-i+1}T_{i-1}\cdots T_1 \pi_n T_{n-1}^{-1}\cdots T_i^{-1}$} for $1\leq i\leq n$. Although similar, the $\theta$'s and the $\xi$'s are not the same and in fact do not have the same spectrum of~$V_{\lambda}$. However, we are able to show that, once symmetrized, the symmetrized vv polynomials agree.\looseness=1

\begin{Corollary}\label{Corollary recovering sym vv Macd}
 The symmetric vv Macdonald polynomials of Dunkl--Luque {\rm\cite[\emph{Theorem} 5.27]{DL_2011}} agree with the $P_T$ of this paper up to nonzero scalars.
\end{Corollary}

\begin{proof}
 The $\sD_n$-modules $V_{\lambda}$ in this paper are isomorphic (after aligning conventions) to the $\sD_n$-modules $\mathcal{M}_{\lambda}$ in Dunkl--Luque's paper. Dunkl and Luque characterize the symmetric vector valued Macdonald polynomials as eigenvectors with distinct eigenvalues for the operator~\smash{$\xi^{(n)}_1+\cdots + \xi^{(n)}_n$} acting on \smash{$\epsilon^{(n)}(\mathcal{M}_{\lambda})$}. Here \smash{$\xi^{(n)}_i$} are the standard Cherednik elements given in our conventions by \smash{$\xi^{(n)}_i = t^{n-i+1}T_{i-1}\cdots T_1 \pi_n T_{n-1}^{-1}\cdots T_i^{-1}$}. A~simple calculation shows that the spherical DAHA elements \smash{$\epsilon^{(n)}\bigl(\xi^{(n)}_1+\cdots + \xi^{(n)}_n\bigr)\epsilon^{(n)}$} and \smash{$\epsilon^{(n)}\bigl(\theta^{(n)}_1+\cdots + \theta^{(n)}_n\bigr)\epsilon^{(n)}$} are both nonzero scalar multiples of \smash{$\epsilon^{(n)}\pi_n\epsilon^{(n)}$}. Namely,
 \begin{align*}
 \epsilon^{(n)}\bigl(\theta^{(n)}_1+\cdots + \theta^{(n)}_n\bigr)\epsilon^{(n)} &= \epsilon^{(n)}\bigl(t^{-(n-1)}\pi_nT_{n-1}\cdots T_1+\cdots + T_{n-1}^{-1}\cdots T_{1}^{-1}\bigr)\epsilon^{(n)}\\
 &= \frac{1-t^{-n}}{1-t^{-1}}\epsilon^{(n)} \pi_n \epsilon^{(n)}
 \end{align*}
 and similarly,
 \begin{align*}
 \epsilon^{(n)}\bigl(\xi^{(n)}_1+\cdots + \xi^{(n)}_n\bigr)\epsilon^{(n)} &= \epsilon^{(n)}\bigl(t^n \pi_n T_{n-1}^{-1}\cdots T_1^{-1}+\cdots +tT_{n-1}\cdots T_1\pi_n\bigr)\epsilon^{(n)}\\
 & = t\frac{1-t^n}{1-t}\epsilon^{(n)}\pi_n\epsilon^{(n)}.
 \end{align*}
 Since the spectrum of $\epsilon^{(n)}\pi_n\epsilon^{(n)}$ acting on $W_{\lambda}$ is simple, it follows that the $P_T$ are eigenvectors for \smash{$\epsilon^{(n)}\bigl(\xi^{(n)}_1+\cdots + \xi^{(n)}_n\bigr)\epsilon^{(n)}$} and hence agree with the symmetric vector valued Macdonald polynomials of Dunkl--Luque up to re-normalization.
\end{proof}

\subsection[Limit of the W\_lambda\^{}(n)]{Limit of the $\boldsymbol{W_{\lambda^{(n)}}}$}

Here we will leverage the exceptional properties of the connecting maps $\Phi_{\lambda}^{(n)}$ and the polynomials~$F_{\tau}$ to define limits of the spherical modules $W_{\lambda^{(n)}}$. To start, we identify a~special stability property for the $P_{T}$ elements.

\begin{Corollary}\label{stability for MacD poly}
 For $T \in \RSSYT\bigl(\lambda^{(n)}\bigr)$, let $T' \in \RSSYT\bigl(\lambda^{(n+1)}\bigr)$ be such that $T(\square) = T'(\square)$ for $ \square \in \lambda^{(n)}$ and $T'(\square_0) = 0$ for $\square_0 \in \lambda^{(n+1)}/\lambda^{(n)}$.
 Then,
\smash{$ \Phi^{(n)}_{\lambda}(P_{T'})= P_T$}.
\end{Corollary}

\begin{proof}
 Note that restriction from \smash{$\lambda^{(n+1)}$} to \smash{$\lambda^{(n)}$} identifies \smash{$\PSYT\bigl(\lambda^{(n)};T\bigr)$} as the subset of~$\smash{\tau \in \PSYT\bigl(\lambda^{(n+1)};T'\bigr)}$ with $\tau(\square_0) = n+1$. Thus, by using Proposition~\ref{stability for locally maximal psyt} in conjunction with Corollary~\ref{expansion of sym into nonsym}, we find that
 \begin{align*}
 \Phi^{(n)}_{\lambda}(P_{T'}) &
 = \sum_{\tau \in \PSYT(\lambda^{(n+1)};T')} \prod_{(\square_1,\square_2) \in \Inv(\tau)} \left( \frac{q^{T(\square_1)}t^{c(\square_1) +1} - q^{T(\square_2)}t^{c(\square_2) }}{q^{T(\square_1)}t^{c(\square_1)} - q^{T(\square_2)}t^{c(\square_2)}} \right) \Phi^{(n)}_{\lambda}(F_{\tau}) \\
 & = \sum_{\substack{\tau \in \PSYT(\lambda^{(n+1)};T')\\ \tau(\square_0) = n+1}} \prod_{(\square_1,\square_2) \in \Inv(\tau)} \left( \frac{q^{T(\square_1)}t^{c(\square_1) +1} - q^{T(\square_2)}t^{c(\square_2) }}{q^{T(\square_1)}t^{c(\square_1)} - q^{T(\square_2)}t^{c(\square_2)}} \right)F_{\tau|_{\lambda^{(n)}}} \\
 & = \sum_{\substack{\tau \in \PSYT(\lambda^{(n+1)};T')\\ \tau(\square_0) = n+1}} \prod_{(\square_1,\square_2) \in \Inv(\tau|_{\lambda^{(n)}})} \left( \frac{q^{T(\square_1)}t^{c(\square_1) +1} - q^{T(\square_2)}t^{c(\square_2) }}{q^{T(\square_1)}t^{c(\square_1)} - q^{T(\square_2)}t^{c(\square_2)}} \right)F_{\tau|_{\lambda^{(n)}}} \\
 & = \sum_{\tau \in \PSYT(\lambda^{(n)};T)} \prod_{(\square_1,\square_2) \in \Inv(\tau)} \left( \frac{q^{T(\square_1)}t^{c(\square_1) +1} - q^{T(\square_2)}t^{c(\square_2) }}{q^{T(\square_1)}t^{c(\square_1)} - q^{T(\square_2)}t^{c(\square_2)}} \right) F_{\tau}
 = P_T. \tag*{\qed}
 \end{align*} \renewcommand{\qed}{}
\end{proof}

We will use the above result to define limits of the polynomials $P_{T}$. First, we must define the required associated combinatorics.

\begin{Definition}\label{stable limit Definition1}
 Let $\lambda \in \Par$. Define the infinite diagram \smash{$\lambda^{(\infty)}:= \bigcup_{n \geq n_{\lambda}} \lambda^{(n)}$}. Define $\Omega(\lambda)$ to be the set of all labelings \smash{$T\colon \lambda^{(\infty)} \rightarrow \mathbb{Z}_{\geq 0}$} such that
 \begin{itemize}
\itemsep=0pt
 \item $|\{ \square \in \lambda^{(\infty)} \mid T(\square) \neq 0\}| < \infty$,
 \item $T$ decreases weakly left to right across rows,
 \item $T$ decreases strictly top to bottom down columns.
 \end{itemize}

 For $T \in \Omega(\lambda)$, we define the degree of $T$ as \smash{$\deg(T) := \sum_{\square \in \lambda^{(\infty)}} T(\square)$}. Define the rank of $T$, $\rk(T)$, to be the minimal $n \geq n_{\lambda}$ such that \smash{$T|_{\lambda^{(\infty)}\setminus\lambda^{(n)}} = 0$}.
\end{Definition}

\begin{Example}\label{Example of asymptotic diagram and labeling}
 For $\lambda = (3,2,1)$, below is an example of a~labeling in $\Omega(\lambda)$
 \[
 T= \ytableausetup{centertableaux, boxframe= normal, boxsize= 2.25em}\begin{ytableau}
5&3&3&2&1&0&0&0& \none [\dots]\\
 3& 2 & 0 \\
 1& 1 & \none \\
 0& \none & \none \\
 \end{ytableau} \in \Omega(\lambda).
 \]

 Note that $T$ is strictly decreasing down columns but only weakly decreasing across rows. The degree of $T$ is $\deg(T):= 5+3+1+3+2+1+3+2+1 = 21$ and the rank of $T$ is $\rk(T) = 11$.

\end{Example}

We may now define the main objects of study in this paper.

\begin{Definition}\label{stable limit Definition}
 Define the space \smash{$W_{\lambda}^{(\infty)}$} to be the inverse limit \smash{$\varprojlim W_{\lambda^{(n)}}$} with respect to the maps \smash{$\Phi_{\lambda}^{(n)}$}. Let \smash{$\widetilde{W}_{\lambda}$} be the subspace of all bounded $X$-degree elements of \smash{$W_{\lambda}^{(\infty)}$}. For $T \in \Omega_{\lambda}$, define the \textit{\textit{generalized Macdonald function}}
\smash{$\MacD_T:= \lim_{n} P_{T|_{\lambda^{(n)}}} \in \widetilde{W}_{\lambda}$}.
\end{Definition}

{\samepage\begin{Remark}
 The degree of each $\MacD_T$ is given simply as
 \smash{$
 \mathrm{deg}(\MacD_T) = \deg(T) = \sum_{\square \in \lambda^{(\infty)}} T(\square)$}.
 It is clear from definition that the set of all $\MacD_T$ for $T \in \Omega(\lambda)$ gives a~$\mathbb{Q}(q,t)$-basis of \smash{$\widetilde{W}_{\lambda}$}.
\end{Remark}

Using the stability of the action of the \smash{$P_{r,0}^{(n)}$} and \smash{$P_{0,\ell}^{(n)}$} operators, we may define the following operators.}

\begin{Definition}\label{MacD operator Definition}
 For $r \geq 1$, define the operator \smash{$p_r^{\bullet}\colon W_{\lambda}^{(\infty)} \rightarrow W_{\lambda}^{(\infty)}$} as the limit of multiplication operators \smash{$p_r^{\bullet}:= \lim_n \bigl(X_1^{r}+\dots + X_n^r\bigr)$}.
 For $\ell \in \mathbb{Z}\setminus \{0\}$, define the operator \smash{$\widetilde{\Delta}_r^{(\infty)}\colon W_{\lambda}^{(\infty)} \rightarrow W_{\lambda}^{(\infty)}$} to be the limit
 \[
 \widetilde{\Delta}_{\ell}^{(\infty)}:= \lim_n \Biggl(\sum_{j=1}^{n}\bigl(\theta^{(n)}_j\bigr)^{\ell} - \sum_{\square \in \lambda^{(n)}} t^{\ell c(\square)} \Biggr).
 \]
\end{Definition}

Proposition~\ref{relations for connecting maps} guarantees that for any $r\geq 1$, $p_r^{\bullet}$ defines a~well-defined operator on \smash{$\widetilde{W}_{\lambda}$} which restricts to an operator on \smash{$\widetilde{W}_{\lambda}$} since $p_r^{\bullet}$ increases $X$-degree by $r$ uniformly. A~simple calculation shows the following.

\begin{Lemma}\label{action of MacD operator}
 For all $\ell \in \mathbb{Z}\setminus \{0\}$ and $T \in \Omega(\lambda)$,
 \[
 \widetilde{\Delta}_{\ell}^{(\infty)}( \MacD_T)= \biggl( \sum_{\square \in \lambda^{(\infty)}} \bigl(q^{\ell T(\square)} -1\bigr)t^{\ell c(\square)} \biggr) \MacD_T.
 \]

\end{Lemma}

\begin{proof}
 Let $\ell \in \mathbb{Z}\setminus \{0\}$ and $T \in \Omega(\lambda)$. Then
 \begin{align*}
 \widetilde{\Delta}_{\ell}^{(\infty)}( \MacD_T) &= \lim_{n} \left(\sum_{j=1}^{n}\bigl(\theta^{(n)}_j\bigr)^{\ell} - \sum_{\square \in \lambda^{(n)}} t^{\ell c(\square)} \right) \bigl( \lim_{n} P_{T|_{\lambda^{(n)}}} \bigr) \\
 &= \lim_{n} \left(\sum_{j=1}^{n}\bigl(\theta^{(n)}_j\bigr)^{\ell} - \sum_{\square \in \lambda^{(n)}} t^{\ell c(\square)} \right) (P_{T|_{\lambda^{(n)}}}) \\
 &= \lim_{n} \biggl(\sum_{\square \in \lambda^{(n)}}q^{\ell T(\square)}t^{\ell c(\square)} -\sum_{\square \in \lambda^{(n)}}t^{\ell c(\square)} \biggr) P_{T|_{\lambda^{(n)}}} \\
 &= \lim_{n} \biggl(\sum_{\square \in \lambda^{(n)}} \bigl(q^{\ell T(\square)} -1\bigr)t^{\ell c(\square)} \biggr) P_{T|_{\lambda^{(n)}}}.
 \end{align*}
 Importantly, for $n \geq \rk(T)$
 \[
 \sum_{\square \in \lambda^{(n)}} \bigl(q^{\ell T(\square)} -1\bigr)t^{\ell c(\square)} = \sum_{\square \in \lambda^{(\infty)}} \bigl(q^{\ell T(\square)} -1\bigr)t^{\ell c(\square)}.
 \]
 Therefore,
 \begin{align*}
 \lim_{n} \biggl(\sum_{\square \in \lambda^{(n)}} \bigl(q^{\ell T(\square)} -1\bigr)t^{\ell c(\square)} \biggr) P_{T|_{\lambda^{(n)}}}&
 = \lim_{n} \biggl(\sum_{\square \in \lambda^{(\infty)}} \bigl(q^{\ell T(\square)} -1\bigr)t^{\ell c(\square)} \biggr) P_{T|_{\lambda^{(n)}}} \\
 & = \biggl(\sum_{\square \in \lambda^{(\infty)}} \bigl(q^{\ell T(\square)} -1\bigr)t^{\ell c(\square)} \biggr) \lim_{n} P_{T|_{\lambda^{(n)}}} \\
 & = \biggl(\sum_{\square \in \lambda^{(\infty)}} \bigl(q^{\ell T(\square)} -1\bigr)t^{\ell c(\square)} \biggr)\MacD_T. \tag*{\qed}
 \end{align*} \renewcommand{\qed}{}
\end{proof}

\begin{Corollary}
 For $\ell \in \mathbb{Z}\setminus \{0\}$, the operator $\widetilde{\Delta}_{\ell}^{(\infty)}$ restricts to an operator on \smash{$\widetilde{W}_{\lambda}$}.
\end{Corollary}

\begin{proof}
 Let $\ell \in \mathbb{Z} \setminus \{0\}$. We know that the set $\{ \MacD_{T} \mid T \in \Omega(\lambda)\}$ is a~basis for \smash{$\widetilde{W}_{\lambda}$}. From Lemma~\ref{action of MacD operator}, we further know that \smash{$\widetilde{\Delta}_{\ell}^{(\infty)}$} acts diagonally on this basis. Therefore, \smash{$\widetilde{\Delta}_{\ell}^{(\infty)}$} restricts to an operator on \smash{$\widetilde{W}_{\lambda}$}.
\end{proof}

\begin{Example}
 For $T \in \Omega(3,2,1)$ as in Example~\ref{Example of asymptotic diagram and labeling},
 \begin{align*}
 \widetilde{\Delta}_{1}^{(\infty)}(\MacD_T) ={}& \bigl( \bigl(q^5-1\bigr) + \bigl(q^3-1\bigr)\bigl(t^{-1}+t^1+t^2\bigr) + \bigl(q^2-1\bigr)\bigl(t^0+t^3\bigr)\\
 &+ (q-1)\bigl(t^{-2}+t^{-1}+t^4\bigr)\bigr)\MacD_T.
 \end{align*}

\end{Example}

\subsection[Positive elliptic Hall algebra action on W\_lambda]{Positive elliptic Hall algebra action on $\boldsymbol{\widetilde{W}_{\lambda}}$}

Combining every result of this paper thus far, we are able to define a~novel family of positive EHA representations.

\begin{Theorem} \label{positive EHA module}
 For $\lambda \in \Par$, \smash{$\widetilde{W}_{\lambda}$} is a~graded $\sE^{+}$-module with action determined for $\ell \in \mathbb{Z}\setminus \{0\}$ and $r>0$ by
 \begin{itemize}
\itemsep=0pt
 \item $P_{r,0} \rightarrow q^{r} p_{r}^{\bullet}$,
 \item $P_{0,\ell} \rightarrow \widetilde{\Delta}_{\ell}^{(\infty)}$.
 \end{itemize}
 Further, \smash{$\widetilde{W}_{\lambda}$} is spanned by a~basis of eigenvectors $\{ \MacD_{T}\}_{T \in \Omega(\lambda)}$ with distinct eigenvalues for the Macdonald operator \smash{$\Delta = \widetilde{\Delta}_{1}^{(\infty)}$}.
\end{Theorem}
\begin{proof}
It suffices to establish that the map $\sE^{+} \rightarrow \End_{\mathbb{Q}(q,t)}\bigl(\widetilde{W}_{\lambda}\bigr)$ satisfies the generating relations of $\sE^{+}$. Any such relation is a~non-commutative polynomial expression in $\sE^{+}$ of the form~${ F(P_{0,-r},\dots, P_{0,-1}P_{0,1},\dots,P_{0,r},P_{1,0},\dots,P_{s,0}) = 0}
$
for some $r > 0$ and $s>0$. By an argument of Schiffmann--Vasserot \cite[Lemma 1.3]{SV}, there are automorphisms \smash{$\Gamma^{(n)}$} of \smash{$\sD_n^{\mathrm{sph}}$} such that for all $\ell \in \mathbb{Z}\setminus \{0\}$ and $s>0$,
\smash{$\Gamma^{(n)}\big(P^{(n)}_{0,\ell}\big)= P^{(n)}_{0,\ell} - \sum_{\square \in \lambda^{(n)}} t^{\ell c(\square)}
$}
 and \smash{$\Gamma^{(n)}\bigl(P^{(n)}_{s,0}\bigr) = P^{(n)}_{s,0}$}. By applying the canonical quotient maps \smash{$\Pi_n\colon \widetilde{W}_{\lambda} \rightarrow W_{\lambda^{(n)}}$}, we see using Corollary~\ref{connecting maps and symmetrization of thetas} that as maps
\begin{gather*}
\Pi_nF(P_{0,-r},\dots, P_{0,-1},P_{0,1},\dots,P_{0,r},P_{1,0},\dots,P_{s,0}) \\
\qquad=F\bigl(\Gamma^{(n)}\bigl(P^{(n)}_{0,-r}\bigr),\dots,\Gamma^{(n)}\bigl(P^{(n)}_{0,-1}\bigr),\Gamma^{(n)}\bigl(P^{(n)}_{0,1}\bigr), \dots, \Gamma^{(n)}\bigl(P^{(n)}_{0,r}\bigr),\Gamma^{(n)}\bigl(P^{(n)}_{1,0}\bigr),\\
\phantom{\qquad=}{}\dots,\Gamma^{(n)}\bigl(P^{(n)}_{s,0}\bigr)\bigr)\Pi_n \\
\qquad= \Gamma^{(n)}\bigl(F\bigl(P^{(n)}_{0,-r},\dots, P^{(n)}_{0,-1},P^{(n)}_{0,1},\dots,P^{(n)}_{0,r},P^{(n)}_{1,0},\dots,P^{(n)}_{s,0}\bigr)\bigr) \Pi_n = 0.
\end{gather*}

As this holds for all $n \geq n_{\lambda}$, it follows that $F(P_{0,-r},\dots, P_{0,-1},P_{0,1},\dots,P_{0,r},P_{1,0},\dots,P_{s,0}) = 0$ in \smash{$\End_{\mathbb{Q}(q,t)}\bigl(\widetilde{W}_{\lambda}\bigr)$} as desired. The last statement regarding the spectrum of $\Delta$ follows directly from Proposition~\ref{finite rank gen macd function properties} and Corollary~\ref{stability for MacD poly}.
\end{proof}

\begin{Remark}
 For $T \in \Omega(\lambda)$ and $\ell \in \mathbb{Z} \setminus \{0\}$,
 \[
 P_{0,\ell}(\MacD_T) = \biggl( \sum_{\square \in \lambda^{(\infty)}} \bigl(q^{\ell T(\square)}-1 \bigr) t^{\ell c(\square)} \biggr) \MacD_T.
 \]

\end{Remark}

For $\lambda = \varnothing$, we recover the polynomial representation of $\mathcal{E}^{+}$ on the space of symmetric functions $\Lambda_{q,t}$.

\begin{Proposition}\label{Proposition recover usual Macd}
 $\widetilde{W}_{\varnothing} = \Lambda_{q,t}$ is the standard representation of $\sE^{+}$. In this case, the sets $\Omega(\varnothing)$ and $\Par$ are in bijection and for all $\mu \in \Par$, $\MacD_{\mu} = P_{\mu}\bigl(x;q^{-1},t\bigr)$ up to some nonzero scalar.
\end{Proposition}
\begin{proof}
 We begin by noticing that for all $n\geq 1$, $\varnothing^{(n)} = (n)$ corresponds to the trivial representation $S_{(n)}$ of the finite Hecke algebra $\sH_n$. Thus $\rho_n^{*}(S_{(n)})$ is the trivial representation of the affine Hecke algebra $\sA_n$ and so when we induce from $\sA_n$ to $\sD_n$, we find that the $\sD_n$-module $V_{(n)}$ is the polynomial representation of $\sD_n$ (see \cite[Theorem 3.2.1]{CherednikBook}). As such, $W_{(n)}$ is just the symmetric polynomial representation of \smash{$\sD_n^{\mathrm{sph}}$}. In this case, the connecting maps \smash{$\Phi_{\varnothing}^{(n)}\colon V_{(n+1)} \rightarrow V_{(n)}$} are the usual projections $f(x_1,\dots,x_n,x_{n+1}) \mapsto f(x_1,\dots,x_n,0)$ on symmetric polynomials.
 Therefore, when we limit, we find that $\widetilde{W}_{\varnothing}$ is the space of symmetric functions~$\Lambda_{q,t}$ along with the usual action of $\mathcal{E}^{+}$ \cite[Corollary 1.5]{SV}. The Macdonald element $P_{0,1}$ acts by the classic Macdonald operator $\Delta$ on $\Lambda_{q,t}$ and has 1-dimensional weight spaces indexed by integer partitions~${\mu \in \Par}$ each spanned by the symmetric Macdonald polynomial $P_{\mu}\bigl(x;q^{-1},t\bigr)$. Furthermore, for all~${\mu \in \Par}$,
 \smash{$\Delta\bigl(P_{\mu}\bigl(x;q^{-1},t\bigr)\bigr) = \bigl( \sum_{1\leq i \leq \ell(\mu)}\bigl(q^{\mu_i}-1\bigr)t^{i-1}\bigr) P_{\mu}\bigl(x;q^{-1},t\bigr)$}. The sets $\Omega(\varnothing)$ and $\Par$ are in canonical bijection. Namely, given $T \in \Omega(\varnothing)$ the labels of $T$ are weakly decreasing to the right along~\smash{$\varnothing^{(\infty)}$} and thus define an integer partition $\mu \in \Par$. By Lemma~\ref{action of MacD operator}, for $\mu \in \Omega(\varnothing)= \Par$, $\MacD_{\mu}$ has the same $\Delta$-weight as $P_{\mu}\bigl(x;q^{-1},t\bigr)$. Therefore, $\MacD_{\mu} = P_{\mu}\bigl(x;q^{-1},t\bigr)$ up to a~non-zero scalar.
\end{proof}

We identify a~special element of each \smash{$\widetilde{W}_{\lambda}$}.

\begin{Definition}
 For any $\lambda \in \Par$, define the labeling $T^{\min}_{\lambda}$ of $\lambda^{(\infty)}$ by
$ T^{\min}_{\lambda}(\square)= \#\big\{ \square' \in \lambda^{(\infty)} \mid \square' ~\text{strictly below} ~\square\big\}$.
\end{Definition}

\begin{Lemma}\label{lowest labeling}
The labeling $T^{\min}_{\lambda}$ is the unique element of $\Omega(\lambda)$ of lowest degree $d_{\lambda}:= \sum_{i \geq 1} i\lambda_i$.
\end{Lemma}

\begin{proof}
 It is immediate that since $\lambda$ is a~partition, $T^{\min}_{\lambda} \in \Omega(\lambda)$. Further, by construction each entry of $T^{\min}_{\lambda}$ is chosen minimally in that for any $T \in \Omega(\lambda)$ and $\square \in \lambda^{(\infty)}$, $T^{\min}_{\lambda}(\square) \leq T(\square)$. To see this, simply note that if $T \in \Omega(\lambda)$, $\square \in \lambda^{(\infty)}$, and $\square'$ is the box directly below $\square$, then~${T(\square)> T(\square')}$.
 Hence, $T(\square)$ must be at least as large as the number of boxes strictly below $\square$. Therefore, $T^{\min}_{\lambda}$ has the minimal degree among all elements of $\Omega(\lambda)$. Lastly, the number of boxes $\square \in \lambda^{(\infty)}$ with $T^{\min}_{\lambda}(\square) = i$ is $\lambda_i$ so \smash{$\deg(T^{\min}_{\lambda}) = d_{\lambda}$} as defined above.
\end{proof}

Here we prove that the modules \smash{$\widetilde{W}_{\lambda}$} are mutually non-isomorphic. Although intuitive, it is not immediately clear apriori that the $\Delta$-spectra of the modules \smash{$\widetilde{W}_{\lambda}$} will be distinct. However, it is straightforward to verify that the $\Delta$-weights of the minimal degree generalized Macdonald functions $\MacD_{T^{\min}_{\lambda}}$ are distinct for distinct $\lambda$.

\begin{Proposition}
 For $\lambda, \mu \in \Par$ distinct, \smash{$\widetilde{W}_{\lambda} \ncong \widetilde{W}_{\mu}$} as graded $\mathcal{E}^{+}$-modules.
\end{Proposition}

\begin{proof}
 Let $\lambda, \mu \in \Par$ and suppose that \smash{$f\colon\widetilde{W}_{\lambda} \rightarrow \widetilde{W}_{\mu}$} is a~graded $\mathcal{E}^{+}$-module isomorphism. Then by Lemma~\ref{lowest labeling} we know that
 \smash{$
 f(\MacD_{T^{\min}_{\lambda}}) = \alpha \MacD_{T^{\min}_{\mu}}
$}
 for some nonzero scalar $\alpha \in \mathbb{Q}(q,t)$. Further,
 \begin{align*}
 P_{0,1}(f(\MacD_{T^{\min}_{\lambda}})) &= f(P_{0,1}(\MacD_{T^{\min}_{\lambda}}))
 = f\biggl(\biggl(\sum_{\square\in \lambda^{(\infty)}}\bigl(q^{T^{\min}_{\lambda}(\square)}-1\bigr)t^{c(\square)}\biggr)\MacD_{T^{\min}_{\lambda}}\biggr) \\
 & = \biggl(\sum_{\square\in \lambda^{(\infty)}}\bigl(q^{T^{\min}_{\lambda}(\square)}-1\bigr)t^{c(\square)}\biggr)f(\MacD_{T^{\min}_{\lambda}}) \\
 & = \biggl(\sum_{\square\in \lambda^{(\infty)}}\bigl(q^{T^{\min}_{\lambda}(\square)}-1\bigr)t^{c(\square)}\biggr) \alpha \MacD_{T^{\min}_{\mu}}.
 \end{align*}

 On the other hand,
 \[
 P_{0,1}(\alpha \MacD_{T^{\min}_{\mu}}) = \biggl(\sum_{\square\in \mu^{(\infty)}}\bigl(q^{T^{\min}_{\mu}(\square)}-1\bigr)t^{c(\square)}\biggr)\alpha \MacD_{T^{\min}_{\mu}}.
 \]

 By assumption, $\alpha \neq 0$ so
 \[
 \sum_{\square\in \lambda^{(\infty)}}\bigl(q^{T^{\min}_{\lambda}(\square)}-1\bigr)t^{c(\square)} = \sum_{\square\in \mu^{(\infty)}}\bigl(q^{T^{\min}_{\mu}(\square)}-1\bigr)t^{c(\square)}.
 \]

 This gives
 \[
 \sum_{\square\in \lambda^{(n_{\lambda})}}\bigl(q^{T^{\min}_{\lambda}(\square)}-1\bigr)t^{c(\square)} = \sum_{\square\in \mu^{(n_{\mu})}}\bigl(q^{T^{\min}_{\mu}(\square)}-1\bigr)t^{c(\square)},
 \]
 which after limiting $q \rightarrow 0$ gives
\smash{$
 \sum_{\square \in \lambda}t^{c(\square)} = \sum_{\square \in \mu}t^{c(\square)}$}.
 By comparing the coefficient of $t^d$ for all $d \in \mathbb{Z}$ on both sides of the above equality, we see that $\lambda$ and $\mu$ have the same number of boxes on each diagonal and are therefore equal.
\end{proof}

\section{Pieri rule for generalized Macdonald functions}\label{Pieri Rule section}

Recall that we have defined (see Definition~\ref{MacD operator Definition}) multiplication operators \smash{$p_r^{\bullet}\colon \widetilde{W}_{\lambda}\rightarrow \widetilde{W}_{\lambda}$} for all~${r\geq 1}$. We may similarly consider the operators \smash{$e_r^{\bullet}\colon \widetilde{W}_{\lambda}\rightarrow \widetilde{W}_{\lambda}$} for $r\geq 1$ defined by $e_r^{\bullet} := \lim_{n} e_r(X_1,\dots,X_n) $ where $e_r(u_1,\dots,u_n):=\sum_{1<i_1<\dots< i_r \leq n} u_{i_1}\cdots u_{i_r}$ is the $r$-th elementary symmetric polynomial for any variables $u_1,\dots,u_n$.
 Proposition~\ref{relations for connecting maps} shows that these operators are well defined. The goal of this section is to derive and utilize an explicit combinatorial formula for the action of the multiplication operators $e_r^{\bullet}$ on \smash{$\widetilde{W}_{\lambda}$}. We investigate the $e_1^{\bullet}$ matrix coefficients in the generalized Macdonald basis, which we refer to as the $e_1$-Pieri coefficients, in more depth and show that they satisfy a~simple non-vanishing condition. We use this non-vanishing to prove that the modules~\smash{$\widetilde{W}_{\lambda}$} are cyclic.

\subsection{Pieri rule preliminaries}

First, we establish some useful lemmas.
 We know apriori that for any $\tau \in \PSYT(\lambda;\allowbreak T)$ for $T \in \RSSYT(\lambda)$, $\epsilon^{(n)}(F_{\tau})= (*)P_{T}$ up to some scalar $(*)$. To find our Pieri rule, we~must compute these scalars $(*)$ explicitly. First, we~understand the Hecke-symmetrizations of polynomials of the form $F_{\Min(T)}$.

\begin{Lemma}\label{sym of min term}
 For $T \in \RYT(\lambda)$,
 \[
 \epsilon^{(n)}(F_{\Min(T)}) = \frac{[\mu(T)]_{t}!}{[n]_{t}!} \sum_{\sigma \in \mathfrak{S}_n/\mathfrak{S}_{\mu(T)}} t^{ ({n \choose 2} - {\mu(T) \choose 2}) - \ell(\sigma)}T_{\sigma}(F_{\Min(T)}).
 \]

\end{Lemma}

\begin{proof}
 The result follows from the following simple calculation:
 \begin{align*}
 \epsilon^{(n)}(F_{\Min(T)}) &
 = \frac{1}{[n]_{t}!}\sum_{\sigma \in \mathfrak{S}_n} t^{{n\choose 2} - \ell(\sigma)} T_{\sigma} (F_{\Min(T)}) \\
 & = \frac{1}{[n]_{t}!}\sum_{\sigma \in \mathfrak{S}_n/\mathfrak{S}_{\mu(T)}} \sum_{\gamma \in \mathfrak{S}_{\mu(T)}} t^{{n\choose 2} - \ell(\sigma\gamma)} T_{\sigma \gamma} (F_{\Min(T)}) \\
 & = \frac{1}{[n]_{t}!}\sum_{\sigma \in \mathfrak{S}_n/\mathfrak{S}_{\mu(T)}} \sum_{\gamma \in \mathfrak{S}_{\mu(T)}} t^{{n\choose 2} - \ell(\sigma) - \ell(\gamma)} T_{\sigma} T_{\gamma} (F_{\Min(T)}) \\
 & = \frac{1}{[n]_{t}!}\sum_{\sigma \in \mathfrak{S}_n/\mathfrak{S}_{\mu(T)}} \sum_{\gamma \in \mathfrak{S}_{\mu(T)}} t^{({n\choose 2}- {\mu(T) \choose 2})- \ell(\sigma)} t^{ {\mu(T) \choose 2} - \ell(\sigma) } T_{\sigma} (F_{\Min(T)}) \\
 & = \frac{1}{[n]_{t}!}\sum_{\sigma \in \mathfrak{S}_n/\mathfrak{S}_{\mu(T)}} t^{({n\choose 2}- {\mu(T) \choose 2})- \ell(\sigma)} T_{\sigma} (F_{\Min(T)}) \biggl( \sum_{\gamma \in \mathfrak{S}_{\mu(T)}}t^{ {\mu(T) \choose 2} - \ell(\sigma) }\biggr) \\
 & = \frac{[\mu(T)]_{t}!}{[n]_{t}!} \sum_{\sigma \in \mathfrak{S}_n/\mathfrak{S}_{\mu(T)}} t^{ ({n \choose 2} - {\mu(T) \choose 2}) - \ell(\sigma)}T_{\sigma}(F_{\Min(T)}).\tag*{\qed}
 \end{align*} \renewcommand{\qed}{}
\end{proof}

The action of the $T_i$-operators on polynomials of the form $F_{\Min(T)}$ may be partially understood as follows.

\begin{Lemma}\label{triangularity for min coset reps}
For $\RSSYT(\lambda)$ and $\sigma \in \mathfrak{S}_n/\mathfrak{S}_{\mu(T)}$,
\[
T_{\sigma}(F_{\Min(T)}) = F_{\sigma(\Min(T))} + \sum_{\tau < \sigma(\Min(T))} \eta_{\tau} F_{\tau}
\]
 for some scalars $\eta_{\tau}$.
\end{Lemma}

\begin{proof}
 We proceed by induction using the fact that $\PSYT(\lambda;T)$ is isomorphic to $\mathfrak{S}_n/\mathfrak{S}_{\mu(T)}$ as posets which we saw in Remark~\ref{PSYT orbits using Bruhat}. Certainly, the statement holds trivially for $\tau = \Min(T)$. Take some $\sigma(\Min(T)) = \tau \in \PSYT(\lambda;T)$ with $s_i(\tau) > \tau$ and suppose that
 \[
 T_{\sigma}(F_{\Min(T)}) = F_{\sigma(\Min(T))} + \sum_{\tau' < \sigma(\Min(T))} \eta_{\tau'} F_{\tau'}
 \]
 for some scalars $\eta_{\tau}$.
 Then using Proposition~\ref{weight basis Proposition},
 \begin{align*}
 T_{s_i\sigma}(F_{\Min(T)})
 ={}& T_iT_{\sigma}(F_{\Min(T)})
 = T_iF_{\tau} + \sum_{\tau' < \sigma(\Min(T))} \eta_{\tau'} T_iF_{\tau'} \\
 ={}& F_{s_i(\tau)} + \frac{(1-t)q^{w_{\tau}(i)}t^{c_{\tau}(i)}}{q^{w_{\tau}(i)}t^{c_{\tau}(i)} - q^{w_{\tau}(i+1)}t^{c_{\tau}(i+1)}} F_{\tau} \\
 & + \sum_{\tau' < \sigma(\Min(T))} \eta_{\tau'} \left(F_{s_i(\tau')} + \frac{(1-t)q^{w_{\tau'}(i)}t^{c_{\tau'}(i)}}{q^{w_{\tau'}(i)}t^{c_{\tau'}(i)} - q^{w_{\tau'}(i+1)}t^{c_{\tau'}(i+1)}} F_{\tau'} \right) \\
 ={}& F_{(s_i\sigma)(\Min(T))} + \sum_{\tau' < (s_i\sigma)(\Min(T))}\eta'_{\tau'} F_{\tau'}. \tag*{\qed}
 \end{align*} \renewcommand{\qed}{}
\end{proof}

The above lemmas may now be used to compute the symmetrization of each $F_{\tau}$ in terms of the $P_{T}$ basis.

\begin{Proposition}\label{sym Macd coeff for symmetrization of lowest term}
 For $T \in \RSSYT(\lambda)$,
 \[
 \epsilon^{(n)}(F_{\Min(T)}) = \frac{[\mu(T)]_{t}!}{[n]_{t}!}P_T.
 \]

\end{Proposition}

\begin{proof}
 Recall from Definition~\ref{Macdonald polynomial def} that the coefficient of $F_{\Top(T)}$ in $P_T$ is $1$. We know from the proof of Proposition~\ref{sym AHA submodules} that since $T \in \RSSYT(\lambda)$,
 \smash{$
 \epsilon^{(n)}(F_{\min(T)}) = \alpha P_T$}
 for some nonzero scalar $\alpha$. Let $\sigma_0$ denote the longest element of \smash{$\mathfrak{S}_n/\mathfrak{S}_{\mu(T)}$}. Note that $\sigma_0(\Min(T)) = \Top(T)$. We now use Lemmas~\ref{sym of min term} and~\ref{triangularity for min coset reps} to compute the coefficient of $F_{\Top(T)}$ in \smash{$\epsilon^{(n)}(F_{\min(T)})$} determining~$\alpha$
 \begin{align*}
 \epsilon^{(n)}(F_{\Min(T)}) &
 = \frac{[\mu(T)]_{t}!}{[n]_{t}!} \sum_{\sigma \in \mathfrak{S}_n/\mathfrak{S}_{\mu(T)}} t^{ ({n \choose 2} - {\mu(T) \choose 2}) - \ell(\sigma)}T_{\sigma}(F_{\Min(T)}) \\
 & = \frac{[\mu(T)]_{t}!}{[n]_{t}!} \sum_{\sigma \in \mathfrak{S}_n/\mathfrak{S}_{\mu(T)}} t^{ ({n \choose 2} - {\mu(T) \choose 2}) - \ell(\sigma)} \biggl( F_{\sigma(\Min(T))} + \sum_{\tau < \sigma(\Min(T))} \eta_{\tau}^{\sigma} F_{\tau} \biggr) \\
 & = \frac{[\mu(T)]_{t}!}{[n]_{t}!} F_{\sigma_0(\Min(T))}t^{ ({n \choose 2} - {\mu(T) \choose 2}) - \ell(\sigma_0)} + \sum_{\tau < \sigma_0(\Min(T))} \eta_{\tau}' F_{\tau} \\
 & = \frac{[\mu(T)]_{t}!}{[n]_{t}!}F_{\Top(T)} + \sum_{\tau < \Top(T)} \eta_{\tau}' F_{\tau}.
 \end{align*}
 Therefore, $\alpha = \frac{[\mu(T)]_{t}!}{[n]_{t}!}$.
\end{proof}

Now we compute the scaling factors between each $\epsilon^{(n)}(F_{\tau})$ and $\epsilon^{(n)}(F_{\Top(T)})$.

\begin{Lemma}\label{symmetrization coefficients from top element}
 For $\tau \in \PSYT(\lambda)$ with $\mathfrak{p}_{\lambda}(\tau) = T \in \RSSYT(\lambda)$,
 \[
 \epsilon^{(n)}(F_{\tau}) = \prod_{(\square_1,\square_2) \in \Inv(\tau)} \left( \frac{q^{T(\square_1)}t^{c(\square_1)} - q^{T(\square_2)}t^{c(\square_2) }}{q^{T(\square_1)}t^{c(\square_1)} - q^{T(\square_2)}t^{c(\square_2) +1}} \right) \epsilon^{(n)}(F_{\Top(T)}).
 \]

\end{Lemma}

\begin{proof}
 Let $T \in \RSSYT(\lambda)$ and $\tau \in \PSYT(\lambda;T)$ with $s_i(\tau) > \tau$.
 Then using Proposition~\ref{weight basis Proposition}, we see
 \begin{align*}
 \epsilon^{(n)}(F_{s_i(\tau)})& = \epsilon^{(n)}\left( \left(T_i +\frac{(t-1)q^{w_{\tau}(i)}t^{c_{\tau}(i)}}{q^{w_{\tau}(i)}t^{c_{\tau}(i)} - q^{w_{\tau}(i+1)}t^{c_{\tau}(i+1)}} \right)F_{\tau}\right) \\
 & = \left(1+ \frac{(t-1)q^{w_{\tau}(i)}t^{c_{\tau}(i)}}{q^{w_{\tau}(i)}t^{c_{\tau}(i)} - q^{w_{\tau}(i+1)}t^{c_{\tau}(i+1)}} \right) \epsilon^{(n)}(F_{\tau}) \\
 &= \left( \frac{q^{w_{\tau}(i+1)}t^{c_{\tau}(i+1)}- q^{w_{\tau}(i)}t^{c_{\tau}(i)+1}}{q^{w_{\tau}(i+1)}t^{c_{\tau}(i+1)} - q^{w_{\tau}(i)}t^{c_{\tau}(i)}} \right) \epsilon^{(n)}(F_{\tau}).
 \end{align*}
 Using an induction argument nearly identical to the proof of Corollary~\ref{expansion of sym into nonsym}, we see that for any~${\tau \in \PSYT(\lambda;T)}$
 \[
 \epsilon^{(n)}(F_{\tau}) = \prod_{(\square_1,\square_2) \in \Inv(\tau)} \left( \frac{q^{T(\square_1)}t^{c(\square_1)} - q^{T(\square_2)}t^{c(\square_2) }}{q^{T(\square_1)}t^{c(\square_1)} - q^{T(\square_2)}t^{c(\square_2) +1}} \right) \epsilon^{(n)}(F_{\Top(T)}).\tag*{\qed}
 \] \renewcommand{\qed}{}
\end{proof}

Finally, we may compute each $\epsilon^{(n)}(F_{\tau})$ explicitly in terms of the $P_{T}$.

\begin{Corollary}\label{symmetrization coefficients in general}
 For $\mathfrak{p}_{\lambda}(\tau) = T \in \RSSYT(\lambda)$,
 \[
 \epsilon^{(n)}(F_{\tau}) = K_{T}(q,t)\prod_{(\square_1,\square_2) \in \Inv(\tau)} \left( \frac{q^{T(\square_1)}t^{c(\square_1)} - q^{T(\square_2)}t^{c(\square_2) }}{q^{T(\square_1)}t^{c(\square_1)} - q^{T(\square_2)}t^{c(\square_2) +1}} \right) P_T,
 \]
 where
 \[
 K_T(q,t):= \frac{[\mu(T)]_{t}!}{[n]_{t}!} \prod_{(\square_1,\square_2) \in \Inv(\Min(T))} \left( \frac{q^{T(\square_1)}t^{c(\square_1)} - q^{T(\square_2)}t^{c(\square_2) +1 }}{q^{T(\square_1)}t^{c(\square_1)} - q^{T(\square_2)}t^{c(\square_2)}} \right).
 \]

\end{Corollary}

\begin{proof}
 We begin by noting that from Lemma~\ref{symmetrization coefficients from top element} applied to $\Min(T)$
 \[
 \epsilon^{(n)}(F_{\Min}) = \prod_{(\square_1,\square_2) \in \Inv(\tau)} \left( \frac{q^{T(\square_1)}t^{c(\square_1)} - q^{T(\square_2)}t^{c(\square_2) }}{q^{T(\square_1)}t^{c(\square_1)} - q^{T(\square_2)}t^{c(\square_2) +1}} \right) \epsilon^{(n)}(F_{\Top(T)}).
 \]
 But from Proposition~\ref{sym Macd coeff for symmetrization of lowest term}, we know that $\epsilon^{(n)}(F_{\Min}) = \frac{[\mu(T)]_{t}!}{[n]_{t}!} P_T$ so
 \[
 \prod_{(\square_1,\square_2) \in \Inv(\tau)} \left( \frac{q^{T(\square_1)}t^{c(\square_1)} - q^{T(\square_2)}t^{c(\square_2) }}{q^{T(\square_1)}t^{c(\square_1)} - q^{T(\square_2)}t^{c(\square_2) +1}} \right) \epsilon^{(n)}(F_{\Top(T)}) = \frac{[\mu(T)]_{t}!}{[n]_{t}!} P_T.
 \]

 Thus
\smash{$
 \epsilon^{(n)}(F_{\Top(T)}) = K_T(q,t) P_T
$}
 as defined in the corollary statement above.

 Lastly, we can now use Lemma~\ref{symmetrization coefficients from top element} to finish the proof.
\end{proof}

The last lemma of this section relates the action of $e_r^{\bullet}$ to the action of $\gamma_n^{r}$ on symmetrized elements.

\begin{Lemma}\label{multiplication in terms of gamma}
For $1\leq r \leq n$,
\[\epsilon^{(n)}e_r(X_1,\dots, X_n)\epsilon^{(n)} = t^{-( (n-1)+\cdots + (n-r))}e_r\bigl(1,\dots,t^{n-1}\bigr) \epsilon^{(n)} \gamma_n^r \epsilon^{(n)}.\]
\end{Lemma}
\begin{proof}
 First, we show by induction that for $1 \leq r \leq n$,
 \[
 \gamma_n^{r} = t^{(n-1)+\cdots + (n-r)} \bigl(T_{n-1}^{-1}\cdots T_1^{-1}\bigr)\bigl(T_{n-1}^{-1}\cdots T_2^{-1}T_1\bigr)\cdots \bigl(T_{n-1}^{-1}\cdots T_r^{-1}T_{r-1}\cdots T_1\bigr) X_1\cdots X_r.
 \]
 For $r=1$, we see that
$
 \gamma_n = X_nT_{n-1}\cdots T_1 = t^{n-1}T_{n-1}^{-1}\cdots T_1^{-1}X_1$.
 Suppose this equation holds for some $1\leq r \leq n-1$. Then we have
 \begin{align*}
 \gamma_n^{r+1} = {}&\gamma_n^{r} \gamma_n
 = t^{(n-1)+\cdots + (n-r)} \bigl(T_{n-1}^{-1}\cdots T_1^{-1}\bigr)\bigl(T_{n-1}^{-1}\cdots T_2^{-1}T_1\bigr)\cdots \bigl(T_{n-1}^{-1}\cdots T_r^{-1}T_{r-1}\cdots T_1\bigr) \\
 &\times X_1\cdots X_r t^{n-1}T_{n-1}^{-1}\cdots T_1^{-1}X_1 \\
 = {}& t^{(n-1)+\cdots + (n-r)} t^{n-1}\bigl(T_{n-1}^{-1}\cdots T_1^{-1}\bigr)\bigl(T_{n-1}^{-1}\cdots T_2^{-1}T_1\bigr)\cdots \bigl(T_{n-1}^{-1}\cdots T_r^{-1}T_{r-1}\cdots T_1\bigr) \\
 & \times T_{n-1}^{-1}\cdots T_{r+1}^{-1} X_1\cdots X_r T_r^{-1}\cdots T_1^{-1}X_1.
 \end{align*}
 A~simple calculation verifies that
$
 X_1\cdots X_r T_r^{-1}\cdots T_1^{-1} = t^{-r} T_r\cdots T_1 X_2\cdots X_{r+1}$.
 Therefore,
 \begin{align*}
 \gamma_n^{r+1} ={}& t^{(n-1)+\cdots + (n-r)} t^{n-1}\bigl(T_{n-1}^{-1}\cdots T_1^{-1}\bigr)\bigl(T_{n-1}^{-1}\cdots T_2^{-1}T_1\bigr)\cdots \bigl(T_{n-1}^{-1}\cdots T_r^{-1}T_{r-1}\cdots T_1\bigr) \\
 & \times T_{n-1}^{-1}\cdots T_{r+1}^{-1} t^{-r} T_r\cdots T_1 X_1X_2\cdots X_{r+1} \\
 ={}& t^{(n-1)+\cdots + (n-r) +(n-(r+1))} \bigl(T_{n-1}^{-1}\cdots T_1^{-1}\bigr)\bigl(T_{n-1}^{-1}\cdots T_2^{-1}T_1\bigr) \cdots \\
 & \times \bigl(T_{n-1}^{-1}\cdots T_r^{-1}T_{r-1}\cdots T_1\bigr) \bigl(T_{n-1}^{-1}\cdots T_{r+1}^{-1} T_r\cdots T_1\bigr) X_1\cdots X_{r+1} ,
 \end{align*}
which is of the correct form.

 We see that for any $1\leq r \leq n$,
 \begin{align*}
 \epsilon^{(n)} \gamma_n^{r} \epsilon^{(n)}
 ={}& \epsilon^{(n)}t^{(n-1)+\cdots + (n-r)} \bigl(T_{n-1}^{-1}\cdots T_1^{-1}\bigr)\bigl(T_{n-1}^{-1}\cdots T_2^{-1}T_1\bigr)\cdots\\
 & \times \bigl(T_{n-1}^{-1}\cdots T_r^{-1}T_{r-1}\cdots T_1\bigr) X_1\cdots X_r \epsilon^{(n)} \\
 ={}& t^{(n-1)+\cdots + (n-r)} \epsilon^{(n)} X_1\cdots X_r \epsilon^{(n)}.
 \end{align*}
Suppose that $1 = i_0 \leq i_1 < \dots < i_r \leq i_{r+1} = n $ with $i_j < i_{j+1}-1$ for some $0 \leq j \leq r$. Then
\begin{align*}
 X_{i_1}\cdots X_{i_{j-1}}X_{i_j +1}X_{i_{j+1}}X_{i_{j+2}}\cdots X_{i_r}
 &= X_{i_1}\cdots X_{i_{j-1}} \bigl(t^{-1}T_{i_j}^{-1}X_{i_j}T_{i_j}^{-1}\bigr)X_{i_{j+1}}X_{i_{j+2}}\cdots X_{i_r} \\
 &= tT_{i_j}^{-1} X_{i_1}\cdots X_{i_{j-1}} X_{i_j} X_{i_{j+1}}X_{i_{j+2}}\cdots X_{i_r}T_{i_j}^{-1} \\
 &= tT_{i_j}^{-1}X_{i_1}\cdots X_{i_r} T_{i_j}^{-1} ,
\end{align*}
which shows that
\[
\epsilon^{(n)}X_{i_1}\cdots X_{i_{j-1}}X_{i_j +1}X_{i_{j+1}}X_{i_{j+2}}\cdots X_{i_r} \epsilon^{(n)} = t\epsilon^{(n)}X_{i_1}\cdots X_{i_r}\epsilon^{(n)}.
\]

It follows that for any $1\leq i_1 < \dots < i_r \leq n$,
\[
\epsilon^{(n)}X_{i_1}\cdots X_{i_r}\epsilon^{(n)} = t^{ (i_r -r)+\cdots + (i_1 -1)}\epsilon^{(n)} X_1\cdots X_r\epsilon^{(n)}.
\]

Now we see
\begin{align*}
 \epsilon^{(n)} e_r(X_1,\dots, X_n) \epsilon^{(n)}
 &= \epsilon^{(n)} \biggl(\sum_{1\leq i_1< \dots < i_r \leq n} X_{i_1}\cdots X_{i_r} \biggr) \epsilon^{(n)} \\
 & = \sum_{1\leq i_1< \dots < i_r \leq n} \epsilon^{(n)} X_{i_1}\cdots X_{i_r} \epsilon^{(n)} \\
 & = \sum_{1\leq i_1< \dots < i_r \leq n} t^{ (i_r -r)+\cdots + (i_1 -1)}\epsilon^{(n)} X_1\cdots X_r\epsilon^{(n)} \\
 &= \sum_{1\leq i_1< \dots < i_r \leq n} t^{ (i_r -r)+\cdots + (i_1 -1)} t^{-( (n-1)+\cdots + (n-r))}\epsilon^{(n)} \gamma_n^{r} \epsilon^{(n)} \\
 &= t^{-( (n-1)+\cdots + (n-r))}\biggl(\sum_{1\leq i_1< \dots < i_r \leq n} t^{ (i_1 -1)+\cdots + (i_r -r)}\biggr) \epsilon^{(n)} \gamma_n^{r} \epsilon^{(n)} \\
 & = t^{-( (n-1)+\cdots + (n-r))} e_r\bigl(1,\dots, t^{n-1}\bigr) \epsilon^{(n)} \gamma_n^{r} \epsilon^{(n)}. \tag*{\qed}
\end{align*} \renewcommand{\qed}{}
\end{proof}

\subsection{Pieri rule}

Using the above lemmas, we may derive an explicit formula for the action of $e_r$ on the symmetric vv Macdonald polynomials in the finite variable situation. We then use the stability of the $P_T$ to derive a~similar formula for the $\MacD_{T}$.

\begin{Theorem}\label{Pieri Rule for Finite vars}
 For $T \in \RSSYT(\lambda)$ and $1 \leq r \leq n$, we have the expansion
 \[
 e_r(X_1,\dots, X_n)P_T = \sum_{S} d^{(r)}_{S,T} P_{S},
 \]
 where for $S \in \RSSYT(\lambda)$ we have
 \begin{gather*}
 \frac{d^{(r)}_{S,T}}{t^{ {r\choose 2}}e_r\bigl(1,\dots, t^{n-1}\bigr)K_{S}(q,t)} \\
 \qquad = \sum_{\substack{\tau \in \PSYT(\lambda;T) \\ \text{s.t.}\\ \Psi^r(\tau) \in \PSYT(\lambda;S)}}
 t^{c_{\tau}(1)+\cdots + c_{\tau}(r)}
 \prod_{(\square_1,\square_2) \in \Inv(\tau)} \left( \frac{q^{T(\square_1)}t^{c(\square_1) +1} - q^{T(\square_2)}t^{c(\square_2)}}{q^{T(\square_1)}t^{c(\square_1)} - q^{T(\square_2)}t^{c(\square_2)}} \right) \\
 \phantom{ \qquad = }{} \times \prod_{(\square_1,\square_2) \in \Inv(\Psi^r(\tau))} \left( \frac{q^{S(\square_1)}t^{c(\square_1)} - q^{S(\square_2)}t^{c(\square_2)}}{q^{S(\square_1)}t^{c(\square_1)} - q^{S(\square_2)}t^{c(\square_2) +1}} \right).
 \end{gather*}
Here $S$ ranges over all $S \in \RSSYT(\lambda)$ which can obtained from $T$ by increasing by~$1$ the labels of~$T$ in~$r$ distinct boxes of~$\lambda$.
\end{Theorem}

\begin{proof}
Using Lemma~\ref{expansion of sym into nonsym} and Remark~\ref{action of gamma}, we find
 \begin{gather*}
 e_r(X_1,\dots, X_n)P_T \\
 \qquad = \epsilon^{(n)}e_r(X_1,\dots, X_n)\epsilon^{(n)}(P_T)
 = t^{-((n-1)+\cdots +(n-r))}e_r\bigl(1,\dots, t^{n-1}\bigr)\epsilon^{(n)}\gamma_n^{r} \epsilon^{(n)}(P_T) \\
 \qquad = t^{-((n-1)+\cdots +(n-r))}e_r\bigl(1,\dots, t^{n-1}\bigr)\epsilon^{(n)}\gamma_n^{r}(P_T) \\
 \qquad = t^{-((n-1)+\cdots +(n-r))}e_r\bigl(1,\dots, t^{n-1}\bigr) \\
 \phantom{\qquad =}{} \times \epsilon^{(n)}\gamma_n^{r} \sum_{\tau \in \PSYT(\lambda;T)} \prod_{(\square_1,\square_2) \in \Inv(\tau)} \left( \frac{q^{T(\square_1)}t^{c(\square_1) +1} - q^{T(\square_2)}t^{c(\square_2) }}{q^{T(\square_1)}t^{c(\square_1)} - q^{T(\square_2)}t^{c(\square_2)}} \right) F_{\tau} \\
 \qquad =t^{-((n-1)+\cdots +(n-r))}e_r\bigl(1,\dots, t^{n-1}\bigr) \\
 \phantom{\qquad =}{} \times \sum_{\tau \in \PSYT(\lambda;T)} \prod_{(\square_1,\square_2) \in \Inv(\tau)} \left( \frac{q^{T(\square_1)}t^{c(\square_1) +1} - q^{T(\square_2)}t^{c(\square_2) }}{q^{T(\square_1)}t^{c(\square_1)} - q^{T(\square_2)}t^{c(\square_2)}} \right) \epsilon^{(n)}\gamma_n^{r}(F_{\tau}) \\
 \qquad =t^{r \choose 2}e_r\bigl(1,\dots, t^{n-1}\bigr)\sum_{\tau \in \PSYT(\lambda;T)} \\
 \phantom{\qquad =}{} \times \prod_{(\square_1,\square_2) \in \Inv(\tau)} \left( \frac{q^{T(\square_1)}t^{c(\square_1) +1} - q^{T(\square_2)}t^{c(\square_2) }}{q^{T(\square_1)}t^{c(\square_1)} - q^{T(\square_2)}t^{c(\square_2)}} \right) \epsilon^{(n)}t^{c_{\tau}(1)+\cdots + c_{\tau}(r)} F_{\Psi^{r}(\tau)} \\
 \qquad = t^{r \choose 2}e_r\bigl(1,\dots, t^{n-1}\bigr) \sum_{\tau \in \PSYT(\lambda;T)}t^{c_{\tau}(1)+\cdots + c_{\tau}(r)} \\
 \phantom{\qquad =}{} \times \prod_{(\square_1,\square_2) \in \Inv(\tau)} \left( \frac{q^{T(\square_1)}t^{c(\square_1) +1} - q^{T(\square_2)}t^{c(\square_2) }}{q^{T(\square_1)}t^{c(\square_1)} - q^{T(\square_2)}t^{c(\square_2)}} \right) \epsilon^{(n)}( F_{\Psi^{r}(\tau)}).
 \end{gather*}

 From Corollary~\ref{symmetrization coefficients in general},
\begin{align*}
 \epsilon^{(n)}( F_{\Psi^{r}(\tau)})
 ={}&
 \mathbbm{1}\left(\mathfrak{p}_{\lambda}(\Psi^{r}(\tau)) \in \RSSYT(\lambda)\right) K_{\mathfrak{p}_{\lambda}(\Psi^{r}(\tau))}(q,t) \\
 &\times \prod_{(\square_1,\square_2) \in \Inv(\Psi^{r}(\tau))} \left( \frac{q^{T(\square_1)}t^{c(\square_1) } - q^{T(\square_2)}t^{c(\square_2) }}{q^{T(\square_1)}t^{c(\square_1)} - q^{T(\square_2)}t^{c(\square_2)+1}} \right) P_{\mathfrak{p}_{\lambda}(\Psi^{r}(\tau))}.
\end{align*}

 Hence, by collecting coefficients around each $P_S$ for $S \in \RSSYT(\lambda)$ we see that
 \[
 e_r(X_1,\dots, X_n)P_T = \sum_{S} d^{(r)}_{S,T} P_{S},
 \]
 where \smash{$d^{(r)}_{S,T}$} are as given in the theorem statement above.

 Lastly, if $\tau \in \PSYT(\lambda;T)$, then the boxes of $\lambda$ containing the labels $1,\dots, r$ (with some powers of $q$ given by $T$) are exactly those boxes $\square \in \lambda$ with $\mathfrak{p}_{\lambda}(\Psi^r(\tau))(\square) = T(\square) +1$. Thus if~${S= \mathfrak{p}_{\lambda}(\Psi^r(\tau)) \in \RSSYT(\lambda)}$, then we may obtain $S$ from $T$ by adding the value $1$ to $r$ boxes of $T$ with at most one $1$ added to each box.
\end{proof}

We define the $e_r$-Pieri coefficients as follows.

\begin{Definition}\label{Pieri coeff Definition}
 For $T \in \Omega(\lambda)$ and $r \geq 1$, define \smash{$\mathfrak{d}^{(r)}_{S,T} \in \mathbb{Q}(q,t)$} by
 \[
 e_r^{\bullet}\MacD_T = \sum_{S \in \Omega(\lambda)} \mathfrak{d}^{(r)}_{S,T} \MacD_{S}.
 \]

\end{Definition}

\begin{Remark}
 Note that from Theorem~\ref{Pieri Rule for Finite vars}, it is clear for $T \in \Omega(\lambda)$ and $ r\geq 1$ that each $S \in \Omega(\lambda)$ such that \smash{$\mathfrak{d}^{(r)}_{S,T} \neq 0$} will necessarily be obtained from $T$ by adding $r$ $1$'s to the boxes of $T$ with at most one $1$ being added to each box. As such, the set of such~$S$ is finite. Further, any such~$S$ has $\rk(S) \leq \rk(T) + r$.
\end{Remark}

As an immediate consequence of Theorem~\ref{Pieri Rule for Finite vars} and the definition of $\MacD_{T}$ from Definition~\ref{stable limit Definition}, we obtain the following result.

\begin{Corollary}[Pieri rule]\label{Pieri Rule}
 Let $S,T \in \Omega(\lambda)$ and $r \geq 1$. For all $n \geq \rk(T)+r$,
 \[
 \mathfrak{d}^{(r)}_{S,T} = d^{(r)}_{S|_{\lambda^{(n)}},T|_{\lambda^{(n)}}}.
 \]

\end{Corollary}

\subsection[Non-vanishing for e\_1-Pieri coefficients]{Non-vanishing for $\boldsymbol{e_1}$-Pieri coefficients}

In this subsection, we prove that if $S$, $T$ satisfy a~simple combinatorial relationship, then $\mathfrak{d}_{S,T}^{(1)} \neq 0$. This will be instrumental in the proof that the modules \smash{$\widetilde{W}_{\lambda}$} are cyclic.

\begin{Definition}\label{raising definition}
 Let $\lambda \in \Par$ and $T \in \RSSYT(\lambda)$. A~box $\square_0$ in $\lambda$ is $T$-\textit{\textit{raisable}} if the labeling~$S$ defined by
 \[
 S(\square) = \begin{cases}
 T(\square), & \square \neq \square_0, \\
 T(\square)+1 ,& \square = \square_0
 \end{cases}
 \]
 is also in $\RSSYT(\lambda)$. We say that $S$ is obtained by raising the box $\square_0$ of $T$. Further, we say that $\square_0$ is a~$S$-\textit{\textit{lowerable}} box in $\lambda$.

 We write $T \uparrow S$ if $S$ may be obtained by raising one box of $T$.
\end{Definition}

\begin{Remark}
 We may define a~partial order $\sqsubseteq$ on the set $\RSSYT(\lambda)$ simply by
 \[
 T \sqsubseteq S \leftrightarrow \forall ~\square \in \lambda^{(\infty)}, \qquad T(\square) \leq S(\square).
 \]
 Then the relation $T \uparrow S$ defined in Definition~\ref{raising definition} is simply the cover relation of $\sqsubseteq$.
 Lastly, we may extend the definitions of raisable/lowerable boxes and of the relation $T \uparrow S$ to $\Omega(\lambda)$ analogously in the obvious way.
\end{Remark}

We require the following lemmas regarding the combinatorics of inversion sets. These lemmas are instrumental in proving the non-vanishing of the $e_1$-Pieri coefficients.

\begin{Lemma}\label{gaps between contents}
 Let $\tau \in \PSYT(\lambda;T)$ for $T\in \RSSYT(\lambda)$. If $(\square_1,\square_2)\in \Inv(\tau)$, with $T(\square_1)=T(\square_2)$, then $c(\square_2)-c(\square_1) \geq 2$.
\end{Lemma}

\begin{proof}
 Since $T \in \RSSYT(\lambda)$, for all $n \geq 0$ the set of boxes $\{\square \in \lambda \mid T(\square) = n\}$ is a~skew-diagram consisting of a~union of disjoint horizontal strips. Suppose $(\square_1,\square_2) \in \Inv(\tau)$ with $T(\square_1)=T(\square_2) = n$. Then, $\square_1$ and $\square_2$ cannot be in the same horizontal strip component of~${\{\square \in \lambda \mid T(\square) = n\}}$. Further, $\square_1$ must be to the left of $\square_2$. Hence, $c(\square_2)-c(\square_1) \geq 2$.
\end{proof}

\begin{Lemma}\label{minimal inversion set for raisable labeling}
 Let $T\in \RSSYT(\lambda)$. Given a~$T$-raisable box of $\lambda$, $\square_0$, there exists a~unique $\tau \in \PSYT(\lambda;T)$ such that
 \begin{itemize}
\itemsep=0pt
 \item $\tau(\square_0)= 1q^{a}$ for some $a \geq 0$,
 \item $\inv(\tau) = S(T)(\square_0) -1$.
 \end{itemize}
\end{Lemma}

\begin{proof}
 Since the count $\inv(\tau) = S(T)(\square_0) -1$ is tight, there exists at most one such labeling. We may simply define $\tau \in \PSYT$ by labeling the boxes $\square \in \lambda$ with $\square <_{T} \square_0$ with the labels~${\{2,\dots, S(T)(\square_0) -1 \}}$ following the box ordering $S(T)$. Then fill the boxes $\square >_T \square_0$ with the values $\{ S(T)(\square_0),\dots, n\}$ following the box ordering $S(T)$ where $n = |\lambda$|. Thus, $\tau$ has exactly $S(T)(\square_0) -1$ inversion pairs.
\end{proof}

\begin{Lemma}\label{inversion sets and raising}
 Let $T,T' \in \RSSYT(\lambda)$ with $T\uparrow T'$. Let $\square_0$ be the box of $\lambda$ on which $T$ and $T'$ differ.
 Let $\tau \in \PSYT(\lambda;T)$ with $\Psi(\tau) \in \PSYT(\lambda;T')$. If $\square_1,\square_2 \in \lambda$ with $\square_0 \notin \{\square_1,\square_2\}$, then $(\square_1,\square_2) \in \Inv(\tau)$ if and only if $(\square_1,\square_2) \in \Inv(\Psi(\tau))$. Furthermore, we have
 \begin{itemize}
\itemsep=0pt
 \item $\{(\square_1,\square_2) \in \Inv(\tau)| \square_0 \in \{\square_1,\square_2\} \} = \{ (\square,\square_0) \mid \square <_{T} \square_0\}$ and
 \item $\{(\square_1,\square_2) \in \Inv(\Psi(\tau)) \mid \square_0 \in \{\square_1,\square_2\}\} = \{ (\square_0,\square) \mid \square_0 <_{T'} \square\}$.
 \end{itemize}
\end{Lemma}

\begin{proof}
By assumptions, \smash{$\tau(\square_0) = 1q^{T(\square_0)}$} and \smash{$\Psi(\tau)(\square_0) = nq^{T'(\square_0)} = nq^{T(\square_0)+1}$}, where $n = |\lambda| $. Therefore, the relative orderings on the boxes of $\lambda \setminus \{\square_0\}$ determined by $S(T)$ and $S(T')$ are the same. Likewise, suppose we take we take two boxes $\square_1,\square_2 \in \lambda$ with $\square_0 \notin \{\square_1,\square_2\}$ and $\tau(\square_1) < \tau(\square_2)$. If $T(\square_1) \neq T(\square_2)$, then
\[
T'(\square_1) = T(\square_1) = T(\square_2) = T'(\square_2)
\]
 and so $\Psi(\tau)(\square_1) < \Psi(\tau)(\square_2)$. Suppose instead that $T(\square_1) = T(\square_2)$ and $\tau(\square_1) < \tau(\square_2)$, say $\tau(\square_j) = i_jq^{a}$ with $1< i_1 <i_2$. Then $\Psi(\tau)(\square_1) = (i_1-1)q^{a}, \Psi(\tau)(\square_2) = (i_2-1)q^{a}$, and $a = T'(\square_1) = T'(\square_2)$ so $\Psi(\tau)(\square_1) < \Psi(\tau)(\square_2)$. Thus, $(\square_1,\square_2) \in \Inv(\tau)$ if and only if $(\square_1,\square_2) \in \Inv(\Psi(\tau))$.

Let $(\square_1,\square_2) \in \Inv(\tau)$ with $\square_0 \in \{\square_1,\square_2\}$. Since $\tau(\square_0) = 1q^{T(\square_0)}$, we must have that $\square_2 = \square_0$ and $\square_1 <_{T} \square_0$. Similarly, let $(\square_1,\square_2) \in \Inv(\Psi(\tau))$ with $\square_0 \in \{\square_1,\square_2\}$. Since $\Psi(\tau)(\square_0) = nq^{T(\square_0)}$, we must have that $\square_1 = \square_0$ and $\square_0 <_{T'} \square_2$.
\end{proof}

Putting these lemmas together we may show the following.

\begin{Theorem}\label{non-vanishing of e_1 Pieri}
 If $\lambda \in \Par$ and $T,T' \in \RSSYT(\lambda)$ with $T \uparrow T'$, then $d_{T',T}^{(1)} \neq 0$.
\end{Theorem}
\begin{proof}
 Let $\square_0$ be the $T$-raisable box on which $T$ and $T'$ differ. From Theorem~\ref{Pieri Rule for Finite vars}, we see that
 \begin{align*}
 \frac{d^{(1)}_{T',T}}{\bigl(\frac{1-t^n}{1-t}\bigr)K_{T'}(q,t)} ={}&
 \sum_{\substack{\tau \in \PSYT(\lambda;T) \\ \text{s.t.}\\ \Psi(\tau) \in \PSYT(\lambda;T')}}
 t^{c_{\tau}(1)}
 \prod_{(\square_1,\square_2) \in \Inv(\tau)} \left( \frac{q^{T(\square_1)}t^{c(\square_1) +1} - q^{T(\square_2)}t^{c(\square_2)}}{q^{T(\square_1)}t^{c(\square_1)} - q^{T(\square_2)}t^{c(\square_2)}} \right) \\
 &\times \prod_{(\square_1,\square_2) \in \Inv(\Psi(\tau))} \left( \frac{q^{T'(\square_1)}t^{c(\square_1)} - q^{T'(\square_2)}t^{c(\square_2)}}{q^{T'(\square_1)}t^{c(\square_1)} - q^{T'(\square_2)}t^{c(\square_2) +1}} \right).
\end{align*}
Therefore, it suffices to show that the sum on the right-hand side of the above equation is nonzero. If $\tau \in \PSYT(\lambda;T)$ with $\Psi(\tau) \in \PSYT(\lambda;T')$, then $c_{\tau}(1)= c(\square_0)$. Hence, we may factor out the term $t^{c_{\tau}(1)}= t^{c(\square_0)}$ outside the sum. From Lemma~\ref{inversion sets and raising}, we have the following for any $\tau \in \PSYT(\lambda;T)$ with $\Psi(\tau) \in \PSYT(\lambda;T')$:
\begin{gather*}
 \prod_{(\square_1,\square_2) \in \Inv(\tau)} \left( \frac{q^{T(\square_1)}t^{c(\square_1) +1} - q^{T(\square_2)}t^{c(\square_2)}}{q^{T(\square_1)}t^{c(\square_1)} - q^{T(\square_2)}t^{c(\square_2)}} \right)\\
 \qquad \times
 \prod_{(\square_1,\square_2) \in \Inv(\Psi(\tau))} \left( \frac{q^{T'(\square_1)}t^{c(\square_1)} - q^{T'(\square_2)}t^{c(\square_2)}}{q^{T'(\square_1)}t^{c(\square_1)} - q^{T'(\square_2)}t^{c(\square_2) +1}} \right) \\
 \phantom{ \qquad \times}{} = \prod_{\square <_{T} \square_0} \left( \frac{q^{T(\square)}t^{c(\square) +1} - q^{T(\square_0)}t^{c(\square_0)}}{q^{T(\square)}t^{c(\square)} - q^{T(\square_0)}t^{c(\square_0)}} \right) \prod_{\square_0 <_{T'} \square} \left( \frac{q^{T(\square_0)+1}t^{c(\square_0)} - q^{T(\square)}t^{c(\square)}}{q^{T(\square_0)+1}t^{c(\square_0)} - q^{T(\square)}t^{c(\square)+1}} \right) \\
 \phantom{ \qquad\times =}{} \times \prod_{\substack{(\square_1,\square_2) \in \Inv(\tau)\\ \square_i \neq \square_0}} \left( \frac{q^{T(\square_1)}t^{c(\square_1)+1} - q^{T(\square_2)}t^{c(\square_2)}}{q^{T(\square_1)}t^{c(\square_1)} - q^{T(\square_2)}t^{c(\square_2)+1}} \right).
\end{gather*}

The first two products above are nonzero and do not depend on $\tau$ and can therefore be brought outside the summation
\[
 \sum_{\substack{\tau \in \PSYT(\lambda;T) \\ \text{s.t.}\\ \Psi(\tau) \in \PSYT(\lambda;T')}}.
\]
Hence, it suffices to show that
\[
\sum_{\substack{\tau \in \PSYT(\lambda;T) \\ \text{s.t.}\\ \Psi(\tau) \in \PSYT(\lambda;T')}} \prod_{\substack{(\square_1,\square_2) \in \Inv(\tau)\\ \square_i \neq \square_0}} \left( \frac{q^{T(\square_1)}t^{c(\square_1)+1} - q^{T(\square_2)}t^{c(\square_2)}}{q^{T(\square_1)}t^{c(\square_1)} - q^{T(\square_2)}t^{c(\square_2)+1}} \right) \neq 0.
\]
Notice that we can rewrite the above product terms in the following way:
\begin{gather*}
 \prod_{\substack{(\square_1,\square_2) \in \Inv(\tau)\\ \square_i \neq \square_0}} \left( \frac{q^{T(\square_1)}t^{c(\square_1)+1} - q^{T(\square_2)}t^{c(\square_2)}}{q^{T(\square_1)}t^{c(\square_1)} - q^{T(\square_2)}t^{c(\square_2)+1}} \right) \\
 \qquad = t^{\inv(\tau) -S(T)(\square_0)+1} \prod_{\substack{(\square_1,\square_2) \in \Inv(\tau)\\ \square_i \neq \square_0}} \left( \frac{1 - q^{T(\square_2)-T(\square_1)}t^{c(\square_2)-c(\square_1)-1}}{1 - q^{T(\square_2)-T(\square_1)}t^{c(\square_2)-c(\square_1)+1}} \right) .
\end{gather*}

Therefore,
\begin{gather*}
\sum_{\substack{\tau \in \PSYT(\lambda;T) \\ \text{s.t.}\\ \Psi(\tau) \in \PSYT(\lambda;T')}} \prod_{\substack{(\square_1,\square_2) \in \Inv(\tau)\\ \square_i \neq \square_0}} \left( \frac{q^{T(\square_1)}t^{c(\square_1)+1} - q^{T(\square_2)}t^{c(\square_2)}}{q^{T(\square_1)}t^{c(\square_1)} - q^{T(\square_2)}t^{c(\square_2)+1}} \right) \\
\qquad= \sum_{\substack{\tau \in \PSYT(\lambda;T) \\ \text{s.t.}\\ \Psi(\tau) \in \PSYT(\lambda;T')}} t^{\inv(\tau) -S(T)(\square_0)+1} \prod_{\substack{(\square_1,\square_2) \in \Inv(\tau)\\ \square_i \neq \square_0}} \left( \frac{1 - q^{T(\square_2)-T(\square_1)}t^{c(\square_2)-c(\square_1)-1}}{1 - q^{T(\square_2)-T(\square_1)}t^{c(\square_2)-c(\square_1)+1}} \right) .
\end{gather*}

We have by definition for any inversion pair $(\square_1,\square_2)$ that $T(\square_2)- T(\square_1) \leq 0$. Therefore, by limiting $q \rightarrow \infty$ we see that
\begin{gather*}
\lim_{q\rightarrow \infty} \sum_{\substack{\tau \in \PSYT(\lambda;T) \\ \text{s.t.}\\ \Psi(\tau) \in \PSYT(\lambda;T')}} t^{\inv(\tau) -S(T)(\square_0)+1} \prod_{\substack{(\square_1,\square_2) \in \Inv(\tau)\\ \square_i \neq \square_0}} \left( \frac{1 - q^{T(\square_2)-T(\square_1)}t^{c(\square_2)-c(\square_1)-1}}{1 - q^{T(\square_2)-T(\square_1)}t^{c(\square_2)-c(\square_1)+1}} \right) \\
\qquad= \sum_{\substack{\tau \in \PSYT(\lambda;T) \\ \text{s.t.}\\ \Psi(\tau) \in \PSYT(\lambda;T')}} t^{\inv(\tau) -S(T)(\square_0)+1} \prod_{\substack{(\square_1,\square_2) \in \Inv(\tau)\\ \square_i \neq \square_0 \\ T(\square_1)=T(\square_2)}} \left( \frac{1 - t^{c(\square_2)-c(\square_1)-1}}{1 -t^{c(\square_2)-c(\square_1)+1}} \right).
\end{gather*}

By Lemma~\ref{gaps between contents}, we see that for each of the inversion pairs $(\square_1,\square_2) \in \Inv(\tau)$ for $\tau \in \PSYT(\lambda;T)$ with $\Psi(\tau) \in \PSYT(\lambda;T')$ and $T(\square_1) = T(\square_2)$ that $c(\square_2)-c(\square_1)-1 \geq 1$. Therefore, if we limit $t \rightarrow 0$,
\begin{gather*}
 \lim_{t \rightarrow 0} \sum_{\substack{\tau \in \PSYT(\lambda;T) \\ \text{s.t.}\\ \Psi(\tau) \in \PSYT(\lambda;T')}} t^{\inv(\tau) -S(T)(\square_0)+1} \prod_{\substack{(\square_1,\square_2) \in \Inv(\tau)\\ \square_i \neq \square_0 \\ T(\square_1)=T(\square_2)}} \left( \frac{1 - t^{c(\square_2)-c(\square_1)-1}}{1 -t^{c(\square_2)-c(\square_1)+1}} \right) \\
 \qquad = \sum_{\substack{\tau \in \PSYT(\lambda;T) \\ \text{s.t.}\\ \Psi(\tau) \in \PSYT(\lambda;T')}} \mathbbm{1}\left(\inv(\tau) = S(T)(\square_0)-1\right) \lim_{t \rightarrow 0} \prod_{\substack{(\square_1,\square_2) \in \Inv(\tau)\\ \square_i \neq \square_0 \\ T(\square_1)=T(\square_2)}} \left( \frac{1 - t^{c(\square_2)-c(\square_1)-1}}{1 -t^{c(\square_2)-c(\square_1)+1}} \right) \\
 \qquad = \sum_{\substack{\tau \in \PSYT(\lambda;T) \\ \text{s.t.}\\ \Psi(\tau) \in \PSYT(\lambda;T')}} \mathbbm{1}\left(\inv(\tau) = S(T)(\square_0)-1\right) \prod_{\substack{(\square_1,\square_2) \in \Inv(\tau)\\ \square_i \neq \square_0 \\ T(\square_1)=T(\square_2)}} (1) \\
 \qquad = \#\{ \tau \in \PSYT(\lambda;T) \mid \Psi(\tau) \in \PSYT(\lambda;T'), \inv(\tau)=S(T)(\square_0)-1\}.
\end{gather*}

By Lemma~\ref{minimal inversion set for raisable labeling}, $\#\{ \tau \in \PSYT(\lambda;T) \mid \Psi(\tau) \in \PSYT(\lambda;T'), \inv(\tau)=S(T)(\square_0)-1\} = 1$ which in particular is not $0$. Therefore, $d_{T',T}^{(1)} \neq 0$.
\end{proof}

Using stability, we find the following.

\begin{Corollary}\label{non-vanishing of e_1 Pieri for limit MacD}
 If $\lambda \in \Par$ and $T,T' \in \Omega(\lambda)$ with $T\uparrow T'$, then $\mathfrak{d}_{T'.T}^{(1)} \neq 0$.
\end{Corollary}

\begin{proof}
 From Corollary~\ref{Pieri Rule}, we know that for all $n \geq \rk(T)+1$,
 \[
 \mathfrak{d}_{T',T}^{(1)} = d_{T'|_{\lambda^{(n)}},T|_{\lambda^{(n)}}}^{(1)}.
 \]
 Since $T'$ is obtained from $T$ by increasing the value of a~single box of $T$ by $1$, we know that the same must be true for \smash{$T'|_{\lambda^{(n)}}$} and \smash{$T|_{\lambda^{(n)}}$} for all $n \geq \rk(T)+1$. Therefore, from Theorem~\ref{non-vanishing of e_1 Pieri} we conclude that \smash{$\mathfrak{d}_{T',T}^{(1)} = d_{T'|_{\lambda^{(n)}},T|_{\lambda^{(n)}}}^{(1)} \neq 0$}.
\end{proof}

The non-vanishing of the $e_1$-Pieri coefficients is sufficient to prove that the \smash{$\widetilde{W}_{\lambda}$} are cyclic $\sE^{+}$-modules.

\begin{Corollary}\label{Murnaghan type modules are cyclic +EHA reps}
 For $\lambda \in \Par$, \smash{$\widetilde{W}_{\lambda}$} is a~cyclic $\mathcal{E}^{+}$-module.
\end{Corollary}

\begin{proof}
 We show that \smash{$\MacD_{T^{\min}_{\lambda}}$} is a~cyclic vector for \smash{$\widetilde{W}_{\lambda}$}, i.e.,
 \smash{$\mathcal{E}^{+} \MacD_{T^{\min}_{\lambda}} = \widetilde{W}_{\lambda}$}. It suffices to show that for every $T \in \Omega(\lambda)$ there exists some $A \in \mathcal{E}^{+}$ with \smash{$A(\MacD_{T^{\min}_{\lambda}}) = \MacD_{T}$}. Notice that given any $T \in \Omega(\lambda)$, we may choose any lowerable box $\square_1$ of $T$ and obtain a~labeling $T_1 \in \Omega(\lambda)$ by subtracting the value of $1$ from $\square_1$ in the labeling $T$. Continuing this process, yields a~sequence of labelings $T_1,T_2,\dots$ with $T_{i+1} \uparrow T_i$ which must eventually terminate as $\deg(T_i) = \deg(T)-i$. It is easy to verify that the only element of $\Omega(\lambda)$ without any lowerable boxes is $T^{\min}_{\lambda}$ so the sequence $T_1,T_2,\dots$ must end at $T^{\min}_{\lambda}$. Reversing this process shows that any $T \in \Omega(\lambda)$ may be obtained from $T^{\min}_{\lambda}$ by a~sequence $T^{\min}_{\lambda} = T_1,\dots, T_n= T$ with $T_i\uparrow T_{i+1}$. Hence, by induction it suffices to show that if $T\uparrow T'$, then there exists $A \in \mathcal{E}^{+}$ such that $A(\MacD_{T}) = \MacD_{T'}$.

 Let $T,T' \in \Omega(\lambda)$ with $T\uparrow T'$. Consider the element $A \in \mathcal{E}^{+}$ defined by
 \[
 A:= \prod_{\substack{T\uparrow S \\ S \neq T'}} \left( \frac{P_{0,1}-\sum_{\square \in \lambda^{(\infty)}}(q^{S(\square)}-1)t^{c(\square)}}{\sum_{\square \in \lambda^{(\infty)}}(q^{T'(\square)}-q^{S(\square)})t^{c(\square)}} \right).
 \]

 The denominator of the above product is nonzero since $P_{0,1}$ acts with simple spectrum on~\smash{$\widetilde{W}_{\lambda}$}. Further, as mentioned before the set of $S \in \Omega(\lambda)$ with $T \uparrow S$ is finite so the above product is finite. We have that for $T \uparrow V$
 \begin{align*}
 A(\MacD_V)&= \prod_{\substack{T\uparrow S \\ S \neq T'}} \left( \frac{P_{0,1}-\sum_{\square \in \lambda^{(\infty)}}(q^{S(\square)}-1)t^{c(\square)}}{\sum_{\square \in \lambda^{(\infty)}}(q^{T'(\square)}-q^{S(\square)})t^{c(\square)}} \right)(\MacD_V) \\
 & = \prod_{\substack{T\uparrow S \\ S \neq T'}} \left( \frac{\sum_{\square \in \lambda^{(\infty)}}(q^{V(\square)}-q^{S(\square)})t^{c(\square)}}{\sum_{\square \in \lambda^{(\infty)}}(q^{T'(\square)}-q^{S(\square)})t^{c(\square)}} \right) \MacD_V
 = \delta_{V,T'} \MacD_V.
 \end{align*}
From Corollary~\ref{non-vanishing of e_1 Pieri for limit MacD}, we know that $\smash{\mathfrak{d}^{(1)}_{T',T} \neq 0}$. Therefore, we may consider the ele\-ment~$\smash{A' \in \mathcal{E}^{+}}$ defined by
 \[
 A':= \frac{q^{-1}}{\mathfrak{d}^{(1)}_{T',T}} A P_{1,0}.
 \]

 We find that
 \begin{align*}
 A'(\MacD_T)& = \frac{q^{-1}}{\mathfrak{d}^{(1)}_{T',T}} A P_{1,0} (\MacD_T)
 = \frac{q^{-1}}{\mathfrak{d}^{(1)}_{T',T}} A~qe_1^{\bullet}(\MacD_T)
 = \frac{1}{\mathfrak{d}^{(1)}_{T',T}}A\biggl( \sum_{T\uparrow S} \mathfrak{d}^{(1)}_{S,T} \MacD_S\biggr) \\
 & = \frac{1}{\mathfrak{d}^{(1)}_{T',T}}\sum_{T\uparrow S} \mathfrak{d}^{(1)}_{S,T} \delta_{S,T'}\MacD_S
 = \MacD_{T'}. \tag*{\qed}
 \end{align*} \renewcommand{\qed}{}
\end{proof}

\subsection*{Acknowledgements}

The author would like to thank their graduate advisor Monica Vazirani for her consistent guidance. The author would also like to thank Erik Carlsson, Daniel Orr, and Eugene Gorsky for helpful conversations about the elliptic Hall algebra and the geometry of Hilbert schemes. The author was supported during this work by the 2023 UC Davis Dean's Summer Research Fellowship. The author would also like to thank the reviewers for their helpful and highly detailed feedback which greatly improved the quality of this paper.

\pdfbookmark[1]{References}{ref}
\LastPageEnding

\end{document}